\magnification=1300

  \font\srzec=msbm10  \font\srsap=msbm7  \font\srcin=msbm5

\def\footnoterule{\kern-3pt \hrule width 0truemm \kern2.6pt}
\hsize=12.7cm
\vsize=19cm
\normalbaselineskip=17pt
\parindent=24pt



\font\sm=cmcsc10

\font\doisperm=cmr12

\font\cuczec=msam10

\font\cucsap=msam7


\def\patr{\hbox{\cuczec\char3 }}

\mathchardef\d="0064
\def\e{\mathord{\rm e}}
\def\im{\mathord{\rm i}}
\def\Im{\mathop{\rm Im}}
\def\Re{\mathop{\rm Re}}
\def\Id{\mathop{\rm Id}}
\def\Range{\mathop{\rm Range}}
\def\Dom{\mathop{\rm Dom}}
\def\Sym{\mathop{\rm Sym}}
\def\PSLtwoZ{\mathop{\rm PSL}\left(2,\Z\right)}
\def\PSLtwoR{\mathop{\rm PSL}\left(2,\R\right)}
\def\twoone{{\rm II}_{1}}
\def\twoinfty{{\rm II}_{\infty}}

\def\ggeq{\geq}

\def\sup{\mathop{\rm sup}}
\def\min{\mathop{\rm min}}
\def\max{\mathop{\rm max}}

\def\area{\mathop{\rm area}}
\def\arg{\mathop{\rm arg}}
\def\const{\mathop{\rm const}}
\countdef\nrpag=27
\global\nrpag=1
\global\advance\nrpag by 1
\def\buildrelt#1\over#2\over#3{\mathrel{\mathop{\kern-1pt#2}\limits^{#1}_{#3}}}


\nopagenumbers
\pageno=1
\def\lastheadline{\vbox{\eightit\relax }}
\def\rightheadline{\vbox{\ninerm\noindent\number\nrpag\hfill 
\leftskip0pt \rightskip0pt
  Non-commutative Markov processes in
  free factors \hfill\folio\vskip4pt\hrule}\global\advance\nrpag by 1}
\def\leftheadline{\vbox{\ninerm\noindent\folio\hfill \leftskip0pt 
\rightskip0pt Florin
  R\u adulescu \hfill\number\nrpag\vskip4pt\hrule}\global\advance\nrpag by 1}
\headline{\ifnum\pageno=1{\lastheadline}
\else{\ifodd\pageno\rightheadline\else\leftheadline\fi}\fi}


\font\ninerm=cmr9  \font\eightrm=cmr8  \font\sixrm=cmr6
\font\ninei=cmmi9  \font\eighti=cmmi8  \font\sixi=cmmi6
\font\ninesy=cmsy9 \font\eightsy=cmsy8 \font\sixsy=cmsy6
\font\ninebf=cmbx9 \font\eightbf=cmbx8 
\font\ninett=cmtt9 \font\eighttt=cmtt8
\font\nineit=cmti9 \font\eightit=cmti8
\font\ninesl=cmsl9 \font\eightsl=cmsl8
\skewchar\ninei='177  \skewchar\eighti='177  \skewchar\sixi='177
\skewchar\ninesy='60  \skewchar\eightsy='60  \skewchar\sixsy='60
\hyphenchar\ninett=-1 \hyphenchar\eighttt=-1 \hyphenchar\tentt=-1

\font\srzec=msbm10  \font\srsap=msbm7  \font\srcin=msbm5

\font\caligzec=cmsy10   
\font\caligsap=cmsy7
                        \font\caligcin=cmsy5
\font\goticzec=eufm10   
\font\goticsap=eufm7
                        \font\goticcin=eufm5


\def\tenpoint{\def\rm{\fam0\tenrm}
   \textfont0=\tenrm \scriptfont0=\sevenrm \scriptscriptfont0=\fiverm
   \textfont1=\teni \scriptfont1=\seveni \scriptscriptfont1=\fivei
   \textfont2=\tensy \scriptfont2=\sevensy \scriptscriptfont2=\fivesy
   \textfont3=\tenex \scriptfont3=\tenex \scriptscriptfont3=\tenex
   \textfont8=\caligzec \scriptfont8=\caligsap \scriptscriptfont8=\caligcin
   \textfont9=\goticzec \scriptfont9=\goticsap \scriptscriptfont9=\goticcin
   \textfont10=\srzec \scriptfont10=\srsap \scriptscriptfont10=\srcin
   \textfont11=\cuczec \scriptfont11=\cucsap
\textfont\itfam=\tenit \def\it{\fam\itfam\tenit}
   \textfont\slfam=\tensl \def\sl{\fam\slfam\tensl}
   \textfont\ttfam=\tentt \def\tt{\fam\ttfam\tentt}
   \textfont\bffam=\tenbf \def\scriptfont\bffam=\sevenbf
   \scriptscriptfont\bffam=\fivebf \def\bf{\fam\bffam\tenbf}
   \normalbaselineskip=14pt

\mathchardef\ac="0841
\mathchardef\bc="0842
\mathchardef\cc="0843
\mathchardef\dc="0844
\mathchardef\ec="0845
\mathchardef\fc="0846
\mathchardef\gc="0847
\mathchardef\hc="0848
\mathchardef\ic="0849
\mathchardef\jc="084A
\mathchardef\kc="084B
\mathchardef\lc="084C
\mathchardef\mc="084D
\mathchardef\nc="084E
\mathchardef\oc="084F
\mathchardef\pc="0850
\mathchardef\qc="0851
\mathchardef\rc="0852
\mathchardef\sc="0853
\mathchardef\tc="0854
\mathchardef\uc="0855
\mathchardef\vc="0856
\mathchardef\wc="0857
\mathchardef\xc="0858
\mathchardef\yc="0859
\mathchardef\zc="085A
\mathchardef\ag="0941
\mathchardef\bg="0942
\mathchardef\cg="0943
\mathchardef\dg="0944
\mathchardef\eg="0945
\mathchardef\fg="0946
\mathchardef\Gg="0947
\mathchardef\hg="0948
\mathchardef\ig="0949
\mathchardef\jg="094A
\mathchardef\kg="094B
\mathchardef\lg="094C
\mathchardef\mg="094D
\mathchardef\ng="094E
\mathchardef\og="094F
\mathchardef\pg="0950
\mathchardef\qg="0951
\mathchardef\rg="0952
\mathchardef\sg="0953
\mathchardef\tg="0954
\mathchardef\ug="0955
\mathchardef\vg="0956
\mathchardef\wg="0957
\mathchardef\xg="0958
\mathchardef\yg="0959
\mathchardef\zg="095A
\mathchardef\C="0A43
\mathchardef\D="0A44
\mathchardef\H="0A48
\mathchardef\N="0A4E
\mathchardef\Q="0A51
\mathchardef\R="0A52
\mathchardef\Z="0A5A

\mathchardef\leq="3B36
\mathchardef\geq="3B3E

   \setbox\strutbox=\hbox{\vrule height8.5pt depth3.5pt width0pt}
      \normalbaselines\rm}


\def\ninepoint{\def\rm{\fam0\ninerm}
   \textfont0=\ninerm 
   \textfont1=\ninei  
   \textfont2=\ninesy 
   \textfont3=\tenex  
   \textfont\itfam=\nineit  \def\it{\fam\itfam\nineit}
   \textfont\slfam=\ninesl  \def\sl{\fam\slfam\ninesl}
   \textfont\ttfam=\ninett  \def\tt{\fam\ttfam\ninett}
   \textfont\bffam=\ninebf 
    \scriptscriptfont\bffam=\fivebf  \def\bf{\fam\bffam\ninebf}
   \normalbaselineskip=13pt
   \font\sm=cmcsc9
   \let\sc=\sevenrm  \let\big=\ninebig  \normalbaselines\rm}


\def\eightpoint{\def\rm{\fam0\eightrm}
   \textfont0=\eightrm 
   \textfont1=\eighti  
   \textfont2=\eightsy 
   \textfont3=\tenex   
   \textfont\itfam=\eightit \def\it{\fam\itfam\eightit}
   \textfont\slfam=\eightsl \def\sl{\fam\itfam\eightsl}
   \textfont\ttfam=\eighttt \def\tt{\fam\ttfam\eighttt}
   \textfont\bffam=\eightbf 
   \def\bf{\fam\bffam\eightbf}
   \normalbaselineskip=11pt
   \setbox\strutbox=\hbox{\vrule height7pt depth2pt  width0pt}
   \font\sm=cmcsc8
   \let\sc=\sixrm  \let\big=\eightbig \normalbaselines\rm}

\def\backslash{/}
\def\zxi{[(\overline{z}-\xi)\backslash (-2\im)]}
\def\zueta{[(\overline{z}-\eta)\backslash (-2\im)]}
\def\zizj{[(\overline {z}_i-z_j)\backslash (-2\im)]}
\def\zwz{[(\overline{z}-w)\backslash (-2\im)]}
\def\etxi{[(\overline{\eta}-\xi)\backslash (-2\im)]}

\def\qed{\llap{$\patr$}\llap{$\patr$}}

\tenpoint
\centerline{\doisperm NON-COMMUTATIVE MARKOV PROCESSES}
\centerline{\doisperm IN  FREE GROUP FACTORS, RELATED TO}
\centerline{\doisperm BEREZIN'S QUANTIZATION AND AUTOMORPHIC FORMS}\vskip12pt
\centerline{\ninerm FLORIN R\u ADULESCU}\vskip20pt

{\ninepoint

\leftskip8mm

\rightskip8mm

\noindent {\nineit Abstract}. In this paper we use the description of free
group factors as the von Neumann algebras of Berezin's deformation of 
the upper half-plane,
modulo $\PSLtwoZ$.

The derivative, in the deformation parameter, of the product in the
corresponding algebras,  is a positive Hochschild $2$-cocycle, defined on
a dense subalgebra.
By analyzing the structure of the cocycle we prove that there is a
generator
$\cal L$ for a quantum dynamical semigroup that implements the cocycle on
a strongly dense subalgebra.

For $x$ in the dense subalgebra, ${\cal L} (x)$ is the (diffusion) operator
$${\cal L}(x)=\Lambda (x)
-{1\over2}\{T,x\},$$
where $\Lambda$ is the pointwise (Schur) multiplication
operator with a symbol function related to the logarithm of the automorphic
form
$\Delta$. The operator
$T$ is positive and affiliated with the algebra $\ac_t$ and $T$
  corresponds to $\Lambda(1)$, in a
sense to be made precise in the paper.  After a suitable 
normalization, corresponding to a
principal-value type method adapted for  $\twoone$ factors, $\Lambda$ becomes
(completely) positive on a union of weakly dense subalgebras. 
Moreover, the 2-cyclic cohomology
cocycle associated to the  deformation may be expressed in terms of $\Lambda$.\vskip6pt}\vskip24pt
\centerline{{\ninebf INTRODUCTION}}\vskip12pt


In this paper we analyze the structure of the positive Hochschild cocycle that
determines Berezin's deformation [Be] of the upper half-plane $\H$, modulo
$\PSLtwoZ$.

As described in [Ra], the algebras ${\cal A}_{t,t>1}$ in the 
deformation are $\twoone$ factors
(free group factors, by [Dy] and [Ra], based on [Vo]
whose elements are (reproducing) kernels $k$ that are functions on 
$\H\times\H$,
analytic in the second variable and antianalytic in the first variable,
diagonally $\PSLtwoZ$-invariant and subject to boundedness conditions
(see [Ra]).

The product $k*_tl$ of two such kernels is the convolution product
$$(k*_tl)(\overline z,\xi)=c_t\int_{\H}k(\overline z,\eta)l(\overline \eta,\xi)
[\overline z,\eta,\overline \eta,\xi]^t\,\d \nu_t( \eta) ,\qquad z,\xi\in \H.$$

Here $[\overline z,\eta,\overline \eta,\xi]$ is the cross ratio
$\left( (\overline z-\xi)(\overline \eta-\eta) \right) / \left(
(\overline z-\eta)(\overline \eta-\xi) \right) $,
while $\d \nu_t$ is the measure on
the upper half-plane $\H$ defined by $\d 
\nu_t=(\Im\eta)^{t-2}\,\d\overline {\eta}
\,\d \eta$, and $c_t$ is a constant.

For $k,l$ in a weakly dense subalgebra $\widehat {\cal A}_{t}$, that 
will be constructed later
in the paper, the following Hochschild $2$-cocycle is well defined:
$${\cal C}_t(k,l)=\hbox{the derivative at }t,\hbox{ from above, of 
}s\to k*_sl.$$

Clearly
$${\cal C}_t(k,l)=
{ c_t'\over{ c_t}}(k*_tl)+c_t\int_\H k(\overline {z},\eta)l(\overline 
{\eta},\xi)
[\overline {z},\eta,\overline {\eta},\xi]^t\ln[z,\eta,\overline 
{\eta},\xi]\,\d \nu_0(\eta).$$

In what follows we will prove that ${\cal C}_t$ is
  always a completely positive Hochschild $2$-cocycle
(for example in the sense introduced  in [CoCu]).
More precisely,
for all $k_1,k_2,\ldots ,k_N$ in $\widehat {\ac}_{t_0}$,
$l_1,l_2,\ldots ,l_N$ in $\ac_t$, we have that
$$\sum\limits_{i,j}\tau_{\ac_t}(l_i^*{\cal C}_t(k_i^*,k_j)l_j)\leq 0.$$
  This also holds true
for more general, discrete, subgroups of $\PSLtwoR$.

In the case of $\PSLtwoZ$, it turns out that ${\cal C}_t(k,l)$ behaves
like the corresponding cocycle obtained from the generator of a 
quantum dynamical
semigroup: that is, there exists a
(necessary completely diffusive, i.e., completely conditionally 
negative) ${\cal L}$
such that
$${\cal C}_t(k,l)={\cal L}_t(k*_tl)-k*_t{\cal L}_t(l)-{\cal L}_t(k)*_tl.$$

It turns out that ${\cal L}$ is defined on a unital, dense subalgebra $\dc_t$
of $\ac_t$, and that ${\cal L} (k)$ belongs to the algebra of 
unbounded operators
affiliated with $\ac_t$.
  Moreover, by a restricting to a smaller, dense, {\it but not unital}
  subalgebra $\dc^0_t$,  the completely positive part of ${\cal L}$ will take
  values in the predual $L^1(\ac_t)$.

The construction  of ${\cal L}_t$ is done by using automorphic
forms. Let $\Delta$ be the unique (normalized) automorphic form for $\PSLtwoZ$
in order 12. Then $\Delta$ is not vanishing in $\H$, so that the 
following expression
$$\eqalign{\ln \varphi(\overline {z},\xi)&=
\ln\left(\overline {\Delta (z)}\Delta(\xi)\zxi^{12}\right)\cr
\noalign{\medskip}
&=\overline {\ln\Delta (z)}+\ln\Delta (\xi)+12\ln \zxi,\qquad z,
\xi\in \H ,}$$
is well defined, and diagonally $\Gamma$-invariant, for a suitable 
choice of the logarithmic
function.

Let $\Lambda$ be the multiplication operator on $\ac_t$, 
corresponding to pointwise
(Schur) multiplication of a symbol ${k}$ by $\ln \varphi$.
Then $\Lambda$ is  defined on a weakly dense subalgebra $\dc_t$ of $\ac_t$.
If $\{a,b\}$ denotes the Jordan product $\{a,b\}=ab+ba$, then
$${\cal L}(k)=\Lambda(k)-{1\over2}\{T,k\}$$
where $T$ is related to $\Lambda(1)$ in a sense made explicit in \S 5.
Moreover, by adding a suitable constant times the identity operator 
to the linear map
$-\Lambda$, we get a completely positive map, defined on a weakly dense
subalgebra.

   By analogy with Sauvageot's construction [Sau],
  the Hochschild $2$-cocycle ${\cal C}_t$ corresponds
to a construction of a  cotangent bundle associated
with the deformation. Moreover, there is a ``real and imaginary
part'' of ${\cal C}_t$. Heuristically, this
is analogous to the decomposition of $d$, the exterior derivative, on a K\"ahler
manifold,
  into
$\delta$ and
$\overline{\delta}$ (we owe this analogy to A. Connes).

   The construction of the ``real part'' of ${\cal C}_t$ is done as 
follows. One considers the
  ``Dirichlet form''
${\cal E}_t$ associated
to ${\cal C}_t$, which is defined as follows:
$${\cal E}_t(k,l)=\tau_{\ac_t} ({\cal C}_t(k,l)),$$
defined for $k,l$ in a weakly dense, unital subalgebra $\widehat \ac_t$.
  Out of this one constructs the operator $Y_t$ defined by
$$\langle Y_t(k),l\rangle_{L^2(\ac_t)}={\cal E}_t(k,l),\qquad k,l\in 
\widehat \ac_t.$$

   The imaginary part of ${\cal C}_t$ is rather
  defined as 2-cyclic cohomology cocycle. The formula for this cyclic 
[Ra, Ra2] cocycle is:
$$\Psi_t(k,l,m)=\tau_{\ac_{t}}([{\cal C}_t(k,l)-(\nabla 
Y_t)(k,l)]m),\qquad k,l,m\in \widehat\ac_t,$$
with $$(\nabla Y_t)(k,l)=Y_t(k,l)-kY_tl-Y_t(k)l.$$
  This is a construction similar to one used in [CoCu].

  Let $\chi$ be the antisymmetric form defined on ${\cal D}^0_t$, a 
weakly dense subalgebra
of $\ac_t$, by the formula
$$\chi_t(k,l)={1\over2}[\langle\Lambda k,l\rangle-\langle k,\Lambda (l)\rangle].$$

Then there is a nonzero constant $\beta$, depending on $t$, such that
$$\Psi_t(k,l,m) + \beta\tau_{\ac_t}(klm)= 
\chi_t(kl,m)-\chi_t(k,lm)+\chi_t(mk,l),$$
for $k,l,m$ in ${\cal D}^0_t$.

   We will show in the paper that $L^2(\ac_t)$
  can be identified with the Bargmann-type Hilbert space
of diagonally $\Gamma$-invariant functions on $\H\times\H$,
that are square-summable on $F\times \H$, analytic in the second 
variable and antianalytic in
the first variable. Here $F$ is a fundamental domain for $\PSLtwoZ$ 
in $\H$, and on
$F\times \H$ we consider the invariant measure
$$d(z,\eta)^{2t}\,\d\nu_0(z)\,\d\nu_0(w)=\left({ 
(\Im z)^{1/2}(\Im \eta)^{1/2}\over{
\left| \zueta\right| }}\right)^{2t}\,\d\nu_0(z)\,\d\nu_0(w).$$

  With this identification, the ``real part'' of ${\cal C}_t$
is implemented (on $\widehat \ac_t$) by the analytic Toeplitz
operator on $L^2(\ac_t)$
  (compression of multiplication) of symbol $\ln d$.
The ``imaginary part'' of  ${\cal C}_t$ is implemented (on the 
smaller algebra ${\cal D}^0_t$)
  by the Toeplitz operator, on $L^2(\ac_t)$,
of symbol ${1\over{12}}\ln \varphi$.\vskip6pt

The expression that we have obtained for
${\cal C}_t(k,l)={\cal L}_t(k*_tl)-k*_t{\cal L}_t(l)-{\cal L}_t(k)*_tl,$
  ${\cal L}(k)=\Lambda(k)-{1\over2}\{T,k\}$,
   is  in concordance with known results
in quantum dynamics:
Recall that in Christensen and Evans [CE], by improving a result 
due to Lindblad
[Li] and [GKS], it is proved that for every uniformly
norm-continuous semigroup
  $(\Phi_t)_{t\geq
0}$ of completely positive maps on a von Neumann algebra $\ac$,
the generator ${\cal L}={ \d\;\over{ \d t}}\Phi_t$ has 
the following form:
$${\cal L} (x)=\Psi (x)-{1\over2}\{\Psi (1),x\}+\im [H,x],$$
where $\Psi\colon \ac\to \ac$ is a completely positive map and $H$ is a 
bounded selfadjoint
operator.

For a  semigroup of completely positive maps that is only strongly 
uniformly continuous,
the generator has a similar form, although ${\cal L} (x)$, for $x$ in 
$\ac$,  is  defined as a
quadratic form affiliated to the von Neumann algebra $\ac$.

Conversely, given ${\cal L} $, a minimal semigroup may be constructed 
under certain conditions
(see, e.g., [CF, Ho, MS, GS, Dav]), although the semigroup
might not be conservative (i.e., unital) even if ${\cal L}(1)=0$.

If ${\cal L}(x)=\Lambda (x)+(G^*x+xG)$, let $\widehat {\Lambda }_t{x}=
\e^{-tG^*}x \e^{-tG}$. Then in the case of $\ac=B(H)$, the corresponding
semigroup $\Phi_t$, satisfying the master equation
$${ \d\;\over{ \d t}}\langle \Phi _t(x),\xi,n\rangle=
\langle {\cal L} (\Phi_t (x))\xi,\eta \rangle $$
for $\xi, \eta$ in a dense domain, is constructed by the Dyson 
expansion [Ho]
$$\Phi_t(x)=\widehat {\Lambda }_t(x)+\sum\limits_{n\geq 0}
\mathop{\int\cdots\int}\limits_
{0\leq t_1\leq t_2\leq \cdots \leq t_n<t}\widehat {\Lambda }_{t_1}\circ
\Lambda \circ \widehat {\Lambda }_{t_2-t_1}\circ \cdots \circ \Lambda \circ
\widehat {\Lambda }_{t-t}\,\d t_1\,\d t_2\cdots \d t_n,$$
which is proved to be convergent [CF, MS].

It is not clear if a minimal conservative semigroup exists for the quantum
dynamical generator ${\cal L} _t$ constructed in our paper.
The quantum dynamical generators ${\cal L} _t$ constructed in this paper
have the following formal property:

Assume that there exists a family of completely positive maps
$(\Phi_{s,t})_{s\geq t}$, with
$\Phi_{s,t}\colon \ac_t\to \ac_s$ satisfying the following
  variant of the master equation:
$$\left. { \d\;\over{ \d s}}(\Psi_{s,t}(\Phi_{s,t}(X)))\right| _{s=s_0}=
{\cal L} _{s_0}(\Psi_{s_0,t}(\Phi_{s_0,t}(X))).
\eqno{(0.1)}
$$
Then $\Phi_{s,t}$ would satisfy the Chapmann-Kolmogorov condition:
$$\Phi_{s,t}\Phi_{s,v}=\Phi_{s,v};\qquad s\geq t\geq v,\ \Phi_{s,s}=\Id .$$
Moreover,
$$\left. { \d\;\over{ \d s}}(\Phi_{s,t}(X)*_s\Phi_{s,t}(Y))\right| _{s=s_0}$$
would be
$${\cal C}_{s_0}(\Phi_{s_0,t}(X),\Phi_{s_0,t}(Y))+
{\cal L} _{s_0}(\Phi_{s_0,t}(X))*{s_0}\Phi_{s_0,t}(Y)+
\Phi_{s_0,t}(X)*_{s_0}{\cal L} _{s_0}(\Phi_{s_0,t}(Y)),$$
which by the cocycle property would be
$${\cal L}_{s_0}(\Phi_{s_0,t}(X*_{s_0}Y)).$$

Thus ${ \d\;\over{ \d s}}(\Phi_{s,t}(X)*_s\Phi_{s,t}(Y))=
{ \d\;\over{ \d s}}\Phi_{s,t}(X*_tY)$.  If unicity (conservativity) 
holds, it would follow that
$\Psi_{s,t}\Phi_{s,t}(X)$ would be a
  (unital) multiplicative map from $\ac _t$ into
$\ac_s$.

At  present we do not know if this conservativity condition of the
minimal solution and the subsequent considerations hold true.\vskip24pt

{\sm Acknowledgement. }{\it This work was initiated while the author 
was visiting the Erwin
Schroedinger Institute in Wien. This work
  was completed while the \hbox{author} was visiting IHP and IHES,
to which the author is grateful for the excellent conditions and warm 
reception. The author
acknowledges enlightening discussions with
  L.~Beznea, P.~Biane, A.~Connes, P.~Jorgensen, R.~Nest,
J.L.~Sauvageot, and L.~Zsido.}\vskip24pt
\centerline{{\ninebf DEFINITIONS}}
\nobreak\vskip12pt

We recall first some notions associated with Berezin's deformation [Be]
of the upper half-plane that were proved in [Ra]
  (see also [Ra1]), in the $\Gamma$-equivariant context.

We consider the Hilbert space $H_t=H^2(\H,\d \nu_t)$, $t>1$, of 
square-summable analytic
functions on the upper half-plane $\H$, with respect to the measure
$\d \nu_t=(\Im z)^{t-2}\,\d\overline{z}\,\d z$. $\d \nu_0$ is the
$\PSLtwoR$-invariant measure on $\H$. This space occurs as the Hilbert space
for the series
of projective unitary irreducible representations $\pi_t$ of $\PSLtwoR$
on $H_t$, $t>1$ [Sal, Puk].

Recall that $\pi_t(g)$, $g=\left(\matrix{a&b\cr
\noalign{\smallskip}
c&d\cr}\right)$
in $\PSLtwoR$, are defined by means of left translation (using the 
M\" obius action
of $\PSLtwoR$ on $\H$)
by the formula
$$(\pi_t(g)f)(z)=f(g^{-1}z)(cz+d)^{-t},\qquad z\in \H,\ f\in H_t.$$
Here the factor $(cz+d)^{-t}$ for
$g=\left(\matrix{a&b\cr
\noalign{\smallskip}
c&d\cr}\right)$
is defined by using a preselected branch of $\ln (cz+d)$ on $\H$, 
which is always
possible [Sal]. If $t=n$ is an integer $\geq 2$, then $\pi_t$ is 
actually a representation
of $\PSLtwoR$, in the discrete series.

Let $\Gamma$ be a discrete subgroup of finite covolume in $\PSLtwoR$ 
and consider
the von Neumann algebra $\ac_t=\{\pi_t(\Gamma)\}'\subseteq B(H_t)\}$
consisting of all operators that commute with $\pi_t(\Gamma)$.

By generalizing a result of [AS], [Co], [Co1], [GHJ], it was proved 
in [Ra] that
$\{\pi_t(\Gamma)\}''$ (the enveloping von Neumann algebra of the image of
$\Gamma$ through $\pi_t$) is isomorphic to ${\cal L} (\Gamma, {\sigma_t})$,
which is the enveloping von Neumann algebra of the image of the left
regular cocycle representation of $\Gamma $ into $B(l^2(\Gamma))$). Thus
${\cal L} (\Gamma, \sigma_t)$ is a $\twoone$ factor. Here
$\sigma_t$ is the cocycle coming from the projective unitary 
representation $\pi_t$.

Therefore, $H_t$, as a left Hilbert module over
  $\{\pi_t(\Gamma)\}''\simeq {\cal L}(\Gamma,\sigma_t)$,
has Murray--von Neumann dimension (see, e.g., [GHJ])
equal to $\left( ( t-1 ) / \pi \right) \mathop{\hbox{covol}}(\Gamma)$
(this generalizes to projective unitary representations, by the formula 
in [AS, Co, GH]).
The precise formula is
$$\mathop{\hbox{dim}}_{{\cal L}(\Gamma,{\sigma_t})}H_t=
\mathop{\hbox{dim}}_{\{\pi_t (\Gamma)\}''}H_t=
{ t-1\over{ \pi}}\mathop{\hbox{covol}}(\Gamma).$$
Hence the commutant $\ac_t$ is isomorphic to
${\cal L} (\Gamma,\sigma_t)_{\left( ( t-1 ) / \pi \right) 
\mathop{\rm covol}(\Gamma)}$.
We use the convention to denote by $M_t$, for a type $\twoone$  factor $M$,
the isomorphism class of $eMe$, with $e$  an idempotent of trace $t$.
If $t>1$, then one has to replace $M$ by $M\otimes M_N(\C)$ (see [MvN]).

When $\Gamma=\PSLtwoZ$, the class
  of the cocycle $\sigma_t$ vanishes (although not
  in the bounded cohomology, see
[BG]).
Consequently, since  in this case [GHJ]
$${ t-1 \over \pi } 
\mathop{\hbox{covol}}(\Gamma)={ t-1 \over 12 },$$
it follows that when $\Gamma =\PSLtwoZ$ we have
$$\ac_t\simeq L(\PSLtwoZ)_{( t-1 ) / 12 }.$$

We want to analyze the algebras $\ac_t$ by means of Berezin's deformation
of $\H$. Recall that the Hilbert space $H_t$ has reproducing vectors
$e^t_z$, $z\in \H$, that are defined by the condition $\langle 
f,e^t_z\rangle=
f(z)$, for all $f$ in $H$. The precise formula is
$$e^t_z(\xi)=\langle e^t_{\xi},e^t_z\rangle=
{ c_t\over{ \zxi^t}},\qquad \xi\in \H,\ c_t={ t-1\over{ 4\pi}}.$$

Each operator $A$ in $B(H_t)$ then has
a reproducing kernel $\widehat {A}(\overline{z},\xi)$. To obtain the
Berezin symbol,
  one  normalizes
so that the symbol of $A=\Id $ is the identical function $1$.

Thus the Berezin symbol of $A$ is a bivariable function
  on $\H\times\H$, antianalytic in the
first variable, analytic in the second, and given by
$$\widehat {A}(\overline{z},\xi)={ \langle Ae^t_z,e^t_{\xi}\rangle\over{
\langle e^t_z,e^t_{\xi}\rangle}},\qquad \overline{z},\xi\in \H.$$

We have that $\langle Ae^t_z,e^t_{\xi}\rangle$ is a reproducing kernel for
$A\in B(H_t)$, and hence the formula for the symbol $\widehat {AB}$ of 
the composition of
two operators $A,B$ in $B(H_t)$ is computed as
$$\widehat {AB}(\overline{z},\xi)
\langle e^t_z,e^t_{\xi}\rangle\langle ABe^t_z,e^t_{\xi}\rangle
=\langle e^t_z,e^t_{\xi}\rangle \int_\H
\langle Ae^t_z,e^t_{\eta}\rangle \langle Be^t_{\eta},e^t_{\xi}\rangle
\,\d\nu_t(\eta).$$\vskip6pt

\noindent {\sm Definition} 0.1.\kern.3em
{\it
By making explicit the kernels involved in the product, one obtains 
the following formula:
Let $\widehat {A}(\overline{z},\xi)=k(\overline{z},\xi)$,
$\widehat{B}(\overline{z},\xi)=l(\overline{\eta},\xi)$,
and let $(k*_tl)(\overline{z},\xi)$
be the symbol of $AB$ in $H_t$. Then
$$(k*_tl)(\overline{z},\xi)=c_t\int_{\H}(k(\overline{z},\xi))
(l(\overline{\eta},\xi))[\overline{z},\eta,\overline{\eta},\xi]^t\,\d \nu_0(\eta)
\eqno{(0.2)}$$
with $[\overline{z},\eta,\overline{\eta},\xi]=
\left( (\overline{z}-\xi)(\overline{\eta}-\eta) \right) / \left(
(\overline{z}-\eta)(\overline{\eta}-\xi) \right)$.

Here one uses the choice of the branch of $\ln (\overline{z}-\xi)\in 
[-\pi,\pi]$
that appears  in the definition of $e^t_z$ {\rm (}see {\rm [Sal]).} }\vskip6pt

The above definition can be extended, when the integrals are 
convergent, to an (associative)
operation on the space of bivariable kernels, by the formula (0.2).
One problem that remains  open is  to determine when a given bivariable function
represents a bounded operator on $H_t$.

Let $d(\overline{z},\eta)=\left( 
(\Im z)^{1/2}(\Im \eta)^{1/2} \right) / \left( 
\left| \zueta\right| \right)$ 
for $z,\eta$ in $\H$. Then $ d(\overline{z},\eta)^2$
is the hyperbolic cosine of the hyperbolic distance between
$z,\eta$ in $\H$.
The following criterion was proven in [Ra].\vskip6pt
{\sm Criterion} 0.2.\kern.3em
{\it
  Let $h$ be a bivariable function on $\H\times \H$,
antianalytic in the first variable, and analytic in the second variable.
Consider the following norm:
  $\|h\|_t^{\,\widehat {}}$ is the maximum of the two
quantities
$$\matrix{\displaystyle \sup_{z\in 
\H}\int\left| h(\overline{z},\eta)\right| 
( d(z,\eta))^t\,\d\nu_0(\eta),\cr
\noalign{\medskip}
\displaystyle \sup_{\eta\in \H}\int\left| h(\overline{z},\eta)\right| 
( d(z,\eta))^t\,\d\nu_0(z).\cr}$$

Then $\|h\|_t^{\,\widehat {}}$ is a norm on $B(H_t)$, finer than
  the uniform norm, and the vector
space of all elements in $B(H_t)$ whose kernels have finite
$\|\,\cdot \,\|_t^{\,\widehat {}}$ norm is an involutive, weakly dense, 
unital, normal subalgebra of
$B(H_t)$. We denote this algebra by $\widehat {B(H_t)}$.
}\vskip6pt

In [Ra] we proved a much more precise statement about the algebra
$\widehat {B(H_t)}$:\vskip6pt

{\sm Proposition} 0.3 [Ra].\kern.3em
{\it
The algebra of symbols corresponding
to $\widehat {B(H_t)}$ is closed under all the product operations $*_s$,
for $s\geq t$.
In particular $\widehat {B(H_t)}$ embeds continuously into
$\widehat {B(H_t)}$ and its image is closed under the product in
$\widehat {B(H_t)}$.
}\vskip6pt

Since this statement will play an essential role in proving that the 
domains of some linear maps
  in our paper, are algebras, we'll briefly recall the proof of this 
proposition:

Assume that $k,l$ are kernels such that
$\|k\|_t^{\,\widehat {}},\|l\|_t^{\,\widehat {}}<\infty$.
Consider the product of $k,l$ in $\ac_s$. We are estimating
$$\int\left| (k*_sl)(\overline{z},\xi)\right| \left| d(z,\xi)\right| ^t\,\d\nu_0(\xi).$$
This should be uniformly bounded in $z$.

The integrals are bounded by
$$\mathop{\int\kern-6.pt \int}\limits_{\H^2}\left| k(\overline{z},\eta)\right| 
\left| l(\overline{\eta},\xi)\right| \left| [\overline{z},\eta,\overline{\eta},\xi]\right| ^s\left| d(z,\xi)\right| ^t
\,\d\nu_0(\eta)\,\d\nu_0(\zeta).$$

Since obviously $$\left| [\overline{z},\eta,\overline{\eta},\xi]\right| ^s=
\displaystyle \left[{  d(z,\eta) d(\eta,\xi)\over{  d(z,\xi)}}\right]^s,$$
the integral is bounded by
$$\mathop{\int\kern-6.pt 
\int}\limits_{\H^2}\left| k(\overline{z},\eta)\right| \left|  d(\overline{z},\eta)\right| ^t
\left| l(\overline{\eta},\xi)\right| \left|  d(\overline{\eta},\xi)^t\right| \cdot M(z,\eta,\xi)
\,\d\nu_0(\eta,\xi).$$

If we can show that $M(z,\eta,\xi)$ is a bounded function on 
$\H\times \H\times \H$,
then the last integral will be bounded by
$\|M\|_{\infty}\|k\|_t^{\,\widehat {}}
\|l\|_t^{\,\widehat {}}$.

But it is easy to see that
$$M(z,\eta,\xi)=
\left|{  d(\overline{z},\eta)
 d(\overline{\eta},\xi)\over{  d(\overline{z},\xi)}}\right|^{s-t}=
\left| [\overline{z},\eta,\overline{\eta},\xi]\right| ^{s-t}.$$

This is a diagonally $\PSLtwoR$-invariant function on $\H\times \H\times \H$.
Since $ d(z,\eta)$ is an intrinsic notion of the geometry on $\H$
we can replace $\H$ by $\D$, the unit disk. Then the expression of $ d(z',\xi')$
becomes: $\displaystyle { (1-\left| z'\right| ^2)^{1/2} (1-\left| \xi'\right| ^2)^{1/2} \over{ 
\left| 1-\overline{z'}\xi'\right| }}$,
$z',\xi'\in \D$.
We thus consider $M$ as a function of three variables $z',\eta',\xi'\in \D$.
By $\PSLtwoR$-invariance when computing the maximum we may let
$\eta = 0$ and we have
$$\eqalign{\displaystyle M(z',0,\xi') &=
\left|{  d(\overline{z},0) d(0,\xi')\over{  d(z',\xi')}}\right|^{s-t}\cr
\noalign{\medskip}
&=
\left|{ (1-\left| z'\right| ^2)^{1/2} (1-\left| \xi'\right| ^2)^{1/2} \over{  d(z',\xi')}}
\right|^{s-t}\cr
\noalign{\medskip}
&=
\left| (1-z'\xi')\right| ^{s-t}\leq 2}$$
since $t>1$.
This completes the proof of Proposition 0.3.\hfill\qed\vskip6pt

In [Ra] we proved that there is a natural symbol map
$\Psi_{s,t}\colon B(H_t)\to B(H_s)$ defined as follows:\vskip6pt
{\sm Definition} 0.4.\kern.3em
{\it
Let $\Psi_{s,t}\colon B(H_t)\to B(H_s)$
be the map that assigns to every operator
$A$ in $B(H_t)$ of Berezin symbol $\widehat {A}(\overline {z},\eta)$,
$\overline {z},\eta\in \H$,  the  operator
$\Psi_{s,t}(A)$ on $B(H_s)$ whose Berezin symbol {\rm (}as operator on 
$H_s${\rm )} coincides
with the symbol of $A$. Then $\Psi_{s,t}$ is continuous
on $B(H_s)$.
}\vskip6pt

A proof of this will be given in Section~1 and we will in fact prove even more,
  that is, that
$\Psi_{s,t}$ is a completely positive map.

Obviously one has
$$\eqalign{\Psi_{s,t}\Psi_{s,v}&=
\hbox to3em{$\displaystyle\Psi_{s,v}$\hfill}
\hbox{\qquad for }s\geq t\geq v>1,\cr
\noalign{\medskip}
\Psi_{s,s}&=
\hbox to3em{$\Id$\hfill}
\hbox{\qquad for }s>1.}$$

Assume $k,l$ represent two symbols of bounded operators in $B(H_t)$. Then the
product $k*_sl$ makes sense for all $s\geq t$. The following
definition of differentiation of the product structure then appears 
naturally. In this
way we get a canonical Hochschild $2$-cocycle associated with the deformation.\vskip6pt
{\sm Definition-Proposition} 0.5 [Ra].\kern.3em
{\it Fix $1<t_0<t$. Let
$k,l$ be operators in $\widehat {B(H_{t_0})}$. Consider $k*_sl$ for $s\geq t$,
and differentiate pointwise the symbol of this expression at $s=t$. Denote the
corresponding kernel by ${\cal C}_t(k,l)=k*_t'l$. Then ${\cal 
C}_t(k,l)$ corresponds
to a bounded
operator in $B(H_t)$. Moreover, ${\cal C}_t(k,l)$ has the following expression:
$$\eqalign{{\cal C}_t(k,l)&=\displaystyle \left. 
{ \d\;\over{ \d s}}(k*_sl)\right| _{s=t},\cr
\noalign{\medskip}
{\cal C}_t(k,l)(\overline {z},\xi)&=
\displaystyle { c_t'\over{ c_t}}(k*_sl)(\overline {z},\xi)\cr
\noalign{\medskip}
&\displaystyle\qquad +c_t\int_{\H}k(\overline {z},\eta)l
(\overline {\eta},\xi)[\overline{z},\eta,\overline{\eta},\xi]^t
\ln [\overline{z},\eta,\overline{\eta},\xi]\,\d\nu_0(\eta).\cr}$$

Moreover, by differentiation of the associativity property,
  it follows that
${\cal C}_t(k, l)$ defines a Hochschild two-cocycle on the 
weakly dense subalgebra
$\widehat {B(H_{t_0})}$ {\rm (}viewed as a subalgebra of
$B(H_t)$ through the symbol
map\/{\rm ).}
}\vskip6pt

We now specialize this construction  for operators $A\in 
\ac_t=\{\pi_t(\Gamma )\}'$,
that is, operators that commute with the image of $\Gamma $ in $B(H_t)$.
We have the following lemma, which  was proved in [Ra].\vskip6pt

{\sm Lemma} 0.6 [Ra].\kern.3em
{\it Let $\Gamma $ be a discrete subgroup of finite covolume
in $\PSLtwoR$.
Assume $F$ is a fundamental domain of $\Gamma $ in $\H$ {\rm (}of
finite area $\nu_0(F)$ with respect to the $\PSLtwoR$-invariant 
measure $\d\nu_0$ on
$\H${\rm ).}

Let $\ac_t=\{\pi_t(\Gamma )\}'$, which is a
type $\twoone$ factor with trace $\tau $.
Then

{\rm 1)}\kern.3em Any operator $A$ in $\ac_t$ has a diagonally $\Gamma 
$-equivariant kernel
$k=k_A(\overline {z},\xi)$, $z,\xi\in \H$ {\rm (}that is,
$k(\overline {z},\xi)=
k(\overline{\gamma z},\gamma\xi)$, $\gamma\in \Gamma $, $z,\xi\in \H${\rm ).}

{\rm 2)}\kern.3em The trace $\tau_A(k)$ is computed by
$${ 1\over{ \nu_0(F)}}\int_Fk(\overline {z},z)\,\d\nu_0(z).$$

{\rm 3)}\kern.3em More generally, let $P_t$ be the projection from 
$L^2(\H,\d\nu_t)$ onto
$H_t$. Let $f$ be a bounded measurable function on $\H$ that
is $\Gamma $-equivariant and let $M_f$ be the multiplication
  operator on $L^2(\H,\d\nu_t)$ by $f$.
  Let $T_f^t=P_tM_fP_t$ be
the Toeplitz operator on $H_t$ with symbol $M_f$.

Then $T_f^t$ belongs to $\ac_t$ and
$$\tau(T_f^tA)={ 1\over{ \nu_0(F)}}\int_Fk_A(\overline {z},z)f(z)\,\d\nu_0(z).$$

{\rm 4)}\kern.3em 
$L^2(\ac_t)$ is identified with the space of all bivariable functions $k$
on $\H\times\H$ that are analytic in the second variable, antianalytic in the
first variable, and diagonally $\Gamma $-invariant.
  The norm of such an element $k$ is given by the formula
$$\|k\|_{2,t}={ 1\over{ \area(F)}}c_t\mathop{\int\kern-6.pt 
\int}\limits_{F\times\H}
\left| k(\overline {z},\eta)\right| ^2 d(z,\eta)^{2t}\,\d\nu_0(z)\,\d\nu_0(\eta).$$
}\vskip6pt

We also note that the algebras $\widehat {B(H_t)}$, and the map 
$\Psi_{s,t}$, $s\geq t$,
have obvious counterparts for $\ac_t$. Obviously $\Psi_{s,t}$ maps
$\ac_t$ into $\ac_s$ for $s\geq t$.\vskip6pt

{\sm Definition} 0.7 [Ra].\kern.3em
{\it
Let $\ac_t=\widehat {B(H_t)}\cap \ac _t$.
Then $\widehat {\ac}_t$ is a weakly dense involutive, unital 
subalgebra of $\ac_t$.

Moreover, $\widehat {\ac}_t$ is closed under any of the operations $*_s$,
for $s\geq t$. This means that $\Psi_{s,t}(k)\Psi_{s,t}(l)\in 
\Psi_{s,t}(\widehat {\ac}_t)$
  for all $k,l$ in $\widehat {\ac}_t$, $s\geq t$.

More generally, $\ac_s$ is contained in $\widehat {\ac}_t$ if $s<t-2$,
and $\widehat {\ac}_r$ is weakly dense in $\ac_t$ if $r\leq t$
{\rm (}and hence $\widehat {\ac}_r$ is weakly dense in
$\ac_t$ if $r\leq t${\rm )
[Ra, Proposition~4.6].}
}\vskip6pt

We also note that, as a consequence of the previous lemma, we can define for
$1<t_0<t$
$${\cal C}_t(k,l)=\left. 
{ \d\;\over{ \d s}}(k*_sl)\right |_{s=t}\hbox{ for }k,l\hbox{ in }
\widehat {\ac}_{t_0}$$
and we have the expression (0.2) of the kernel.

Another way to define ${\cal C}_t(k,l)$ is to fix vectors $\xi,\eta$ 
in $H_t$ and to
consider the derivative
$${ \d\;\over{ \d s}}\langle (k*_sl)\xi,\eta\rangle _{H_t}=
\langle {\cal C}_t(k,l)\xi,\eta\rangle |_{H_t},\ \xi,\eta\in H_t.$$

For $k,l$ in $\widehat {\ac}_{t_0}$, $t_0<t$, this makes sense because
$k*_sl$ is already the kernel of an operator in $\widehat {\ac}_{t_0}$.

\vfill\eject


\centerline{\S o. {\ninebf Outline of the paper}}\vskip12pt

The paper is organized as follows:

In Section 1, we show, based on the facts proved in [Ra], that the symbol maps
$\Psi_{s,t}\colon \ac_t\to \ac_s$, for $s\geq t$, are completely positive,
  unital and trace preserving.
Consequently the derivative of the multiplication operation (keeping 
the symbols fixed)
is a positive Hochschild $2$-cocycle (see [CoCu]).
In particular the trace of this Hochschild cocycle is a (noncommutative)
Dirichlet form (see [Sau]).\vskip6pt

In Section 2 we analyze positivity properties for families of symbols 
induced by
intertwining operators. As in [GHJ], let $S_{\Delta^{\varepsilon}}$ be the
multiplication operator by $\Delta^{\varepsilon}$, viewed as an operator from
$H_t$ into $H_{t+12\varepsilon}$. Then $S_{\Delta^{\varepsilon}}$ is 
an intertwiner
  between $\pi_t|_{\Gamma }$
and $\pi_{t+12\varepsilon}|_{\Gamma }$, with $\Gamma =\PSLtwoZ$.
Here we use the following branch for $\ln (cz+d)=\ln (j(\gamma,z))$, 
which appears
in the definition of $\pi_t(\gamma)$,
$\gamma=\left(\matrix{a&b\cr
\noalign{\smallskip}
c&d\cr}\right)$ in $\PSLtwoZ$, $\gamma\in \Gamma $.
We define $\ln (j(\gamma,z))=\ln(\Delta(\gamma^{-1}z))-\ln 
\Delta(z)$, which is possible
since there is a canonical choice for $\ln\Delta(z)$.

We use the fact that $S_{\Delta^{\varepsilon}}
S_{\Delta^{\varepsilon}}^*$ is a decreasing family of operators, converging to
the identity as $\varepsilon\to 0$. Let
$$\varphi(\overline {z},\xi)=\overline{\Delta(z)}\Delta(\xi)\zxi^{12}.$$
Then
  $$\ln \varphi(\overline {z},\xi)=\overline{ \ln\Delta(z)}+\ln\Delta(\xi)+
12\ln\zxi,$$
  has the  property that
$$\left[\left(-{ 1\over{ 12}}\ln \varphi(\overline {z}_i,z_j)+{ 
c_t'\over{ c_t}}\right)
[(\overline{z}_i-z_j)/(-2\hbox{i})]^{-t}\right]_{i,j}$$
is a positive matrix for all  $z_1,z_2,\ldots ,z_n$ in $\H$ and for all
$t>1$.\vskip6pt

In Section~3 we use the positivity proven in Section~2 to check that 
the operator
of symbol multiplication by $\displaystyle (\ln 
\varphi)\left(\varphi^{\varepsilon}+C_t\right)$
(for a suitable constant $C_t$, depending only on $t$) is well defined on a
weakly dense subalgebra
of $\ac_t$. This operator gives a completely positive map on this subalgebra.

By a principal-value procedure, valid in a type~$\twoone$ factor, we deduce
  that multiplication
by $(-\ln\varphi+c_{t,\widetilde{\ac}})$ is a completely positive map
$\Lambda$, on a weakly dense unital subalgebra $\widetilde{\ac}$ of $\ac_t$
($c_{t,\widetilde{\ac}}$  is a constant that only depends on $t$ and $\ac$).
Multiplication
by $(\ln \varphi)$ maps $\widetilde{\ac}$ into the operator 
affiliated with $\ac_t$.

In particular $\Lambda (1)$ is affiliated with $\ac_t$. We obtain 
this result by checking
that the kernels $-\left( 
(\varphi^{\varepsilon}-\Id ) / \varepsilon \right) $
are decreasing as $\varepsilon\downarrow \varepsilon_0$, $\varepsilon_0>0$
(up to a small linear perturbation), to $\varphi^{\varepsilon_0}\ln \varphi$,
  plus a suitable constant.

This is not surprising as $\Lambda(1)=\ln \varphi (\overline {z},\xi)$
barely fails
the summability criteria  for $L^1(\ac_t)$.

In Section~4, we analyze the derivatives $X_t$, at $t$, of the 
intertwining maps
$\theta_{s,t}\colon \ac_t\to \ac_s$, $s\geq t$, with 
$\theta_{s,t}(k)=S_{\Delta^{(s-t)/12}}k
S_{\Delta^{(s-t)/12}}^*$. The derivatives ($X_t$) are,
  up to a multiplicative constant, the operators defined in  Section 
3. The operator
$X_t$ is defined on a weakly dense unital subalgebra of $\ac_t$.

We take the derivative of the identity satisfied by $\theta_{s,t}$,
which is
$$\theta_{s,t}(k*_tT^t_
{{\varphi}^{(s-t)/12}}*_t l )=\theta_{s,t}(k)*_s\theta_{s,t}(l ).$$

This gives the identity
$$X_t(k*_tl)+k*_tT^t_{\ln \varphi }*_t\varphi ={\cal 
C}_t(k,l)+X_tk*_tl+k*_tX_tl,$$
which holds on a weakly dense (nonunital) subalgebra.\vskip6pt

Based on an estimate on the growth of the function
$\left| \ln\Delta(z)\Delta^{\varepsilon}(z)\right| $, $z\in \H$, for fixed $\varepsilon >0$,
we prove in Section~5 that the positive, affiliated operators $-\Lambda(1)$
and $-T_{\ln \varphi }^t$ are equal operators.
We prove this by showing that there  is an increasing family 
$A_{\varepsilon }$ in $\ac_t$ and
dense  domains $\dc_0,\dc_1$ (where $\dc_0$ is affiliated to $\ac_t$) such that
$\langle A_{\varepsilon }\xi,\xi\rangle\to \langle -\Lambda(1)\xi,\xi\rangle$
for $\xi$ in $\dc_0$ and
$\langle A_{\varepsilon }\xi,\xi\rangle\to \langle -T_{\ln \varphi 
}\xi,\xi\rangle$
for $\xi$ in $\dc_1$.\vskip6pt

In Section~6 we analyze the cyclic cocycle associated with the 
deformation which
is obtained from the positive Hochschild cocycle by discarding a trivial part.

The precise formula is $$\Psi_t(k,l,m)=
\tau_{\ac_t}([{\cal C}_t(k,l)-Y_t(kl)+(Y_tk)l+k(Y_tl)]m),$$
for $k,l,m$ in a dense subalgebra, and
$$\langle Y_tk,l\rangle=-{1\over2}\tau_{\ac_t}({\cal C}_t(k,l^*)).$$

We reprove a result in [Ra], that the cyclic cohomology cocycle
$$\Psi(k,l,m)-\hbox{cst}\ \tau(klm)$$ is
implemented by $$\chi_t(k,l^*)=\langle X_tk,l\rangle -\langle k,X_tl\rangle $$
for $k,l$ in a dense subalgebra. Since the constant in the above 
formula is nonzero,
this corresponds to nontriviality of $\Psi_t$ on this dense subalgebra.\vskip6pt

In Section 7 we analyze a dual form of the coboundary for ${\cal C}_t(k,l)$,
in which multiplication by $\varphi $ is rather replaced by the 
Toeplitz operator
of multiplication by $\overline {\varphi }$ (compressed to $L^2(\ac_t)$).
It turns out that the roles of $\Lambda (1)$ and $T_{\ln \varphi }^t$ 
are reversed
in the functional equation satisfied by the coboundary.\vskip6pt

In the appendix, giving up the complete
  positivity requirement and the algebra requirement
on the domain of the corresponding maps, we find some more general
   coboundaries for ${\cal C}_t$, which were hinted at in [Ra].

\vfill\eject

\centerline{\S 1. {\ninebf Complete positivity for the Hochschild $2$-cocycle
associated with }}
\centerline{{\ninebf the deformation}}\vskip12pt

In this section we prove the positivity condition on the Hochschild $2$-cocycle
associated with Berezin's deformation.

Denote for $z,\eta$ in $\H$ the expression
$$ d(\overline {z},\eta)=
\displaystyle { (\Im z)^{1/2}(\Im\eta)^{1/2}\over{ \zueta}}$$
and recall that $\left|  d(\overline {z},\eta)\right| ^2$ is the hyperbolic cosine 
of the hyperbolic
distance between $z,\eta\in \H$.

In [Ra] we introduced the following seminorm, defined  for $A\in 
B(H_s)$,  given by the
kernel $k=k_A(\overline {z},\xi)$, $z,\xi\in \H$:
$$
\eqalign{\|A\|_s^{\,\widehat {}} &= \|k\|_s^{\,\widehat 
{}}\cr
&= \max
\left(\sup_{z\in \H}\int \left| k(\overline {z},\eta)\right| \left|  d(z,\eta)\right| ^s
\,\d\nu_0(\eta),
\sup_{z\in \H}\int \left| k(z,\eta)\right|  \left|  d(z,\eta)\right| ^s
\,\d\nu_0(z)\right)
 .}
$$

The subspace of all elements $A$ in $B(H_s)$ (respectively $\ac_s$) such that
$\|A\|_s^{\,\widehat {}}$ is finite is a closed, involutive Banach 
subalgebra of $B(H_s)$
(respectively $\ac_s$) that we denote by $\widehat {B(H_s)}$
(respectively $\widehat{\ac}_s$).

In [Ra] we proved that in fact $\widehat{\ac}_s$ (or $\widehat {B(H_s)}$) is
also closed under any of the products $*_t$, for $t\geq s$, and that 
there is  a universal
constant $c_{s,t}$, depending on $s,t$, such that
$$\|k*_tl\|_s^{\,\widehat {}}\leq c_{s,t}\|k\|_s^{\,\widehat  
{}}\|l\|_s^{\,\widehat {}},\qquad
k,l\in \widehat {\ac}_s.$$

Also $\widehat {\ac}_s$ (or $\widehat {B(H_s)}$) is weakly dense in $\ac_s$
(respectively $B(H_s)$).

Let $\Psi_{s,t}$,  $s\geq t>1$, be the map that associates to any $A$
in $B(H_t)$ (respectively $\ac_t$) the corresponding element in $B(H_s)$
(respectively $\ac_s$) having the same symbol (that is, $\Psi_{s,t}(A)\in \ac_s$
has the same symbol as $A$ in $\ac_t$). Then $\Psi_{s,t}$ maps
$\ac_t$ continuously into $\ac_s$.

In the next proposition we prove that $\Psi_{s,t}$ is a completely 
positive map. This is based
on the following positivity criterion proved in [Ra].\vskip6pt

{\sm Lemma} 1.1 (Positivity criterion).\kern.3em
{\it A kernel $k(\overline {z},\xi)$ defines
  a positive bounded
operator in $B(H_t)$, of norm less than $1$, if and only if for all $N$ in $\N$
and for all $z_1,z_2,\ldots ,z_{N}$ in $\H$ we have that the 
following matrix inequality
holds:
$$0\leq \left[{ k(\overline {z}_i,z_j)\over{ \zizj^t}}\right]_{i,j=1}^N\leq
\left[{ 1\over{ \zizj^t}}\right]_{i,j=1}^N.$$
}\vskip6pt

This criterion obviously holds at the level of matrices of elements in 
$M_p(\C)\otimes \ac_t$.\vskip6pt

{\sm Lemma} 1.2 (Matrix positivity criterion).\kern.3em
{\it If $[k_{p,q}]_{p,q=1}^N$
is a positive matrix of elements in $\ac_t$ then for all $N$ in $\N$, all
$z_1,z_2,\ldots ,z_N$ in $\H$ the following matrix is positive definite:
$$\left[{ k_{p,q}(\overline {z}_i,z_j)\over{ \zizj^t}}\right]_{(i,p),(j,q)
\in \{1,2,\ldots ,P\}\times \{1,2,\ldots
,N\}}.$$

Conversely, if the entries $k_{p,q}$ represent an  element in $\ac_t$ 
and if the above
matrix is positive, then $[k_{p,q}]_{p,q=1}^N$ is a positive matrix in $\ac_t$.
}\vskip6pt

{\it Proof.} Let $[k_{p,q}]_{p,q=1}^N$ be a matrix in $\ac_t$.
Then $k=[k_{p,q}]$ is positive if and only if for all vectors 
$\xi_1,\xi_2,\ldots ,\xi_N$
in $H_t$ we have that
$$\sum\limits_{p,q}\langle k_{p,q}\xi_p,\xi_q\rangle_{H_t}\geq 0.$$

Since $k_{p,q}=P_tk_{p,q}P_t$, where $P_t$ is the projection from
$L^2(\H,\nu_t)$ onto $H_t$, it turns out that this is equivalent with 
the same statement
which now must be valid for all $\xi_1,\xi_2,\ldots ,\xi_N$ in
$L^2(\H,\nu_t)$.

Thus we have that
$$\sum\limits_{p,q}\mathop{\int\kern-6.pt 
\int}\limits_{\H^2}{k_{p,q}(\overline{z},w)\over
\zwz^t} \xi_p(z)\overline{\xi_q(w)}\,\d\nu_t(z)\,\d\nu_t(w)\geq 0$$
for all
$\xi_1,\xi_2,\ldots ,\xi_N$ in $L^2(\D,\nu_t)$.

We let the vectors $\xi_p$ converge  to the Dirac distributions, for 
all $p=1,2,\ldots ,N$,
$\displaystyle \sum\limits_i\lambda_{ip}\delta_{z_i}
(\Im z_i)^{-(t-2)}$, for all $p=1,2,\ldots ,N$. By the above 
inequality we get
$$\sum\limits_{i,j,p,q}{ k_{p,q}(\overline {z}_i,z_j)\over{ \zizj^t}}
\overline{\lambda}_{ip}\lambda_{jq}\geq 0$$
for all choices of $\{\lambda_{ij}\}$ in $\C$.
This corresponds exactly to the fact that the matrix
$$\left[{ k_{p,q}(\overline {z}_i,z_j)\over{ \zizj^t}}\right]_{(p,i),(q,j)}$$
is positive.\hfill\qed\vskip6pt

{\sm Proposition} 1.3.\kern.3em
{\it The map $\Psi_{s,t}\colon \ac_t\to \ac_s$ which 
sends an element
$A$ in $\ac_t$ into the corresponding element in $\ac_s$, having the 
same symbol,
is unital and completely positive.}\vskip6pt

{\it Proof.} This is a consequence of the fact
[ShS, Sal] that the matrix
$$\displaystyle\left[ { 1\over{ \zizj^{\varepsilon }}}\right]_{i,j}$$
is a positive matrix for all $\varepsilon$, all $N$, all
$z_1,z_2,\ldots ,z_n$ in $\H$. Indeed $$\displaystyle { 1\over{ 
\zxi^\varepsilon}}$$
(or $ 1/(1-\overline{z}\xi)^{\varepsilon }$)
is a reproducing kernel for a space of analytic  functions, even if
$\varepsilon <1/2$.\hfill\qed\vskip6pt

We will now follow Lindblad's [Li] argument to deduce that ${\cal C}_t(k,l)$
is a completely positive Hochschild $2$-cocycle.
We recall first the definition of the cocycle ${\cal C}_t$ associated 
with the deformation.\vskip6pt

{\sm Definition} 1.4.\kern.3em
{\it
  Fix $t>s_0>1$. Then the following formula defines a
Hochschild $2$-cocycle on $\widehat {\ac}_{s_0}$.
$$\eqalign{{\cal C}_t(k,l)(\overline {z},\xi)=&\displaystyle
\left.{ \d\;\over{ \d s}}\right|_{{s=t\atop s<t}}(k*_sl)(\overline {z},\xi)\cr
\noalign{\medskip}
=&\displaystyle { c_t'\over{ c_t}}\left(k*_tl\right)(\overline {z},\xi)+
c_t\int_{\H}k(\overline {z},\eta)
l(\overline {\eta},\xi)[\overline{z},\eta,\overline{\eta},\xi]^t\ln
[\overline{z},\eta,\overline{\eta},\xi]\,\d\nu_0(\eta).}$$

Indeed, it was proven in {\rm [Ra]} that the above integral is absolutely 
convergent for $k,l$
in $\widehat {\ac}_{s_0}$, for any $s_0<t$.
}\vskip6pt

The above definition may be thought of also in the following way.
Fix vectors $\xi,\eta$ in $H_t$ and fix $k,l$ in $\ac_{s_0}$.
  Then $k*_tl$, and $k*_sl$ make sense for all
$s_0<s<t$ and they represent bounded operators in $\ac_t$. Thus the 
following derivative
makes sense:
$$\left.{ \d\;\over{ \d s}}\right|_{{s=t\atop s<t}}
\langle k*_sl\xi,\eta\rangle_{H_t},$$
and it turns out to be equal to
$$\langle {\cal C}_t(k,l)\xi,\eta\rangle _{H_t}.$$

In the following lemma  we use the positivity of $\Psi_{s,t}$
to deduce the complete positivity of ${\cal C}_t$.
We recall the following formal formula for ${\cal C}_t$ that was 
proved in [Ra1], [Ra2].\vskip6pt

{\sm Lemma} 1.5 [Ra].\kern.3em
{\it
  Let $1<t_0<t$ and let $k,l,m$ belong to
$\widehat {\ac}_{t_0}$. Then the following holds:
$$\tau_{\ac_t}({\cal C}_t(k,l)*_tm)=
\left.{ \d\;\over{ \d 
s}}\tau_{\ac_s}((k*_sl*_sm)-(k*_tl)*_sm)\right|_{{s=t\atop s>t}}.$$

In a more precise notation, the second term is
$$\left.{ \d\;\over{ \d s}}\tau_{\ac_s}([\Psi_{s,t}(k)*_s\Psi_{s,t}(l)-
\Psi_{s,t}(k*_tl)]*_s\Psi_{s,t}(m))\right|_{{s=t\atop s>t}}.$$
}\vskip6pt

The proof of the lemma is trivial, as long as one uses the absolute convergence
of the integrals, which  follows from the fact that the
kernels belong to an algebra $\widehat {\ac}_{t_0}$, for some $t_0<t$.\vskip6pt

The positivity property that  we are proving for ${\cal C}_t(k,l)$, 
is typical for
coboundaries of the form $D(ab)-D(a)b-aD(b)$, where $D$ is the generator of
dynamical semigroup. It is used by Sauvageot to construct the 
cotangent bimodule
associated with a dynamical semigroup, and many of the properties in [Sau]
can be transferred to ${\cal C}_t$ with the same proof. Such positive 
(or negative)
cocycles appear in the work of Connes and Cuntz (see also [CH]).\vskip6pt

{\sm Proposition} 1.6.\kern.3em
{\it Fix $1<t_0<t$, and for $k,l$ in $\widehat {\ac}_{t_0}$,
define
$$
{\cal C}_t(k,l) :=
   \left.{ \d\;\over{ \d s}}(k*_sl)\right|_{{s>t\atop s=t}}.
$$
Then for all $k_1,k_2,\ldots ,k_N$ in $\widehat {\ac}_{t_0}$,
$l_1,l_2,\ldots ,l_N$ in $\ac_t$, we have that
$$\sum\limits_{i,j}\tau_{\ac_t}(l_i^*{\cal C}_t(k_i^*,k_j)l_j)\geq 0.$$
This is the same as requiring the matrix $({\cal C}_t(k_i^*,k_j))_{i,j}$
to be negative in $M_N(\ac_t)$.
}\vskip6pt

{\it Proof.} For   $s\geq t$, let $f(s)$ be defined by the formula
$$f(s)=\tau 
\left(\sum\limits_{i,j}(k_i^**_sk_j-k_i^**_tk_j)*_t(l_i^**_tl_j)\right).$$
Using the $\Psi_{s,t}$ notation, this is
$$f(s)= \sum\limits_{i,j}\tau\left(\left(\Psi_{s,t}(k_i^*)
\Psi_{s,t}(k_j)-\Psi_{s,t}(k_i^*k_j)\right)
\Psi_{s,t}(l_i^*l_j)\right).$$

By the previous lemma, $f'(t)$ is equal to
$\tau({\cal C}_t(k_i^*,k_j)l_jl_i^*)$. In these terms,  to prove the 
statement we must prove that
$f'(t)\geq 0$. Clearly $f(t)=0$.

By the generalized Cauchy--Schwarz--Stinespring inequality for 
completely positive maps,
and since $\Psi_{s,t}$ is unital, we get that the matrix
$$D_{ij}=[\Psi_{s,t}(k_i^*)\Psi_{s,t}(k_j)-\Psi_{s,t}(k_i^*k_j)]$$
is non-positive. Since $Z_{ij} =\Psi_{s,t}(l_jl_i^*)$ is another positive
matrix in $M_n(\ac_s)$, we obtain that
$$f(s)=\tau_{\ac_s\otimes M_N(\C)}(DZ)$$
is negative.

So $f(s)\leq 0$ for all $s\geq t$, $f(0)=0$. Hence
$\left.{ \d\;\over{ \d s}}f(s)\right|_{{s=t;\;s>t}}$
is negative.\hfill\qed\vskip24pt

\centerline{\ninebf Appendix (to Section 1)}
\nobreak\vskip12pt

We want to emphasize the properties of the trace ${\cal 
E}_t(k,l)=-\tau({\cal C}_t(k,l))$,
$k,l\in \ac_t$. Clearly ${\cal E}_t$ is a positive form on $\ac_t$, and in fact
it is obviously positive definite. Following [Sau], one can prove that
${\cal E}_t$ is a Dirichlet form.
The following expression holds for ${\cal E}_t$.

{\sm Lemma} 1.7.\kern.3em
{\it For $1<t_0<t$, $k,l\in \widehat {\ac}_{t_0}$ we have that
$${\cal E}_t(k,l)=\mathop{\int\kern-6.pt \int}\limits_{F\times 
\H}k(\overline{z},\eta)
\overline{l(\overline {z},\eta)}\left| d(\overline {z},\eta)\right| ^{2t}\ln \left|  d(z,\eta)\right| 
\,\d\nu_0(z,\eta),$$
where $F$ is a fundamental domain for $\Gamma $ in $\H$, and
$$\left|  d(\overline {z},\eta)\right| =\displaystyle \left| { \Im z^{1/2}
\Im \eta^{1/2}\over{ \zueta}}\right| $$
is the hyperbolic cosine of the hyperbolic distance between $z$ and 
$\eta$ in $\H$.
}\vskip6pt

Recall that $L^2(\ac_t)$ is identified [Ra] with the Bargmann-type 
Hilbert space
of functions $k(\overline {z},\eta)$ on $\H\times \H$ that are 
antianalytic in the
first variable, analytic in the second, diagonally $\Gamma $-invariant (that
is, $k(\overline{\gamma}z,\gamma\eta)=k(\overline {z},\eta)$, 
$\gamma\in \Gamma $,
$z,\eta$ in $\H$), and square-summable:
$$\|k\|^2_{L^2(\ac_t)}=c_t
\mathop{\int\kern-6.pt \int}\limits_{F\times \H}\left| k(\overline{z},\eta)\right| ^2
\left|  d(\overline {z},\eta)\right| ^{2t}\,\d\nu_0(z,\eta).$$

Let $\pc_t$ be the projection from the Hilbert space of square-summable
functions $f$ on $\H\times\H$ that are $\Gamma $-invariant and 
square-summable:
$$c_t\mathop{\int\kern-6.pt \int}\limits_{F\times \H}\left| f(\overline{z},\eta)\right| ^2
\left|  d(\overline {z},\eta)\right| ^{2t}\,\d\nu_0(\eta)\,\d\nu_0(z)<\infty.$$

The following proposition is easy to prove, but we won't make any use 
of it in this
paper.\vskip6pt

{\sm Proposition} 1.8.\kern.3em
{\it Let $\varphi $ be a bounded measurable $\Gamma 
$-invariant
function on $\H\times \H$.
Let ${\cal T}_{\varphi }$ be the Toeplitz operator of multiplication 
by $\varphi $
on the Hilbert space $L^2(\ac_t)$, that is, ${\cal T}_{\varphi 
}k=\pc(\varphi k)$,
$k\in L^2(\ac_t)$. Then ${\cal T}_{\varphi }k=\pc_tk\pc_t$, where the 
last composition is in
$\ac_t$, by regarding $k$ as an element affiliated to $\ac_t$.
}\vskip6pt

{\it Remark.} In this setting the positive form ${\cal E}_t$ may be identified
with the quadratic form on $L^2(\ac_t)$ induced by the unbounded 
operator ${\cal T}_{\ln d}$
where $d=\left|  d(\overline {z},\eta)\right| $ is defined as above.

\vfill\eject

\centerline{\S 2. {\ninebf Derivatives of some one-parameter families 
of positive operators}}\vskip12pt

In this section we consider some parametrized families of completely 
positive maps
that are induced by automorphic forms (and fractional powers thereof).
The automorphic forms are used as intertwining operators between the
different representation spaces of $\PSLtwoZ$, consisting of 
analytic functions.

It was proved in [GHJ] that automorphic forms $f$ for $\PSLtwoZ$ of
weight $k$, provide bounded multiplication operators $S_f\colon H_t\to H_{t+k}$.
The boundedness property comes exactly from the fact that one of the conditions
for an automorphic form $f$ of order $k$ is
$$\sup_{z\in \H}\left| f(z)\right| ^2\Im z^k\leq M,$$
which is exactly the condition that the operator of multiplication by $f$
from $H_t$ into $H_{t+k}$ be norm bounded by $M$.

Secondly, the automorphic forms have the (cocycle) $\Gamma 
$-invariance property
as functions on $\H$, that is,
$$f(\gamma^{-1} z)=(cz+d)^{-k}f(z),\qquad z\in \H,\
\gamma=\left(\matrix{a&b\cr
\noalign{\smallskip}
c&d\cr }\right)\in \PSLtwoZ.$$
Since $\pi_t(\gamma)$, $\pi_{t+k}(\gamma)$ act on the corresponding 
Hilbert space
of analytic func-\break tions on $\H$ by multiplication with the automorphic 
factor $\left(cz+d\right)^{-t}$,
respectively
$\left(cz+d\right)^{-t-k}$, this implies exactly that
$$\pi_{t+k}(\gamma)S_f=S_f\pi_{t}(\gamma).$$

Let $f,g$ be automorphic forms of order $k$. Let $F$ be a fundamental 
domain for the group
$\PSLtwoZ$ in $\H$. It was proved in [GHJ] that the trace (in $\ac_t$)
of $S_f^*S_g$ is equal to the Petersson scalar product
$${1\over{\area F}}\langle f,g\rangle =
{1\over{\area F}}\int_Ff(z)\overline {g(z)}(\Im z)^k\,\d\nu_0(z).
\eqno{(2.1)}$$

In the next lemma we will prove that the symbol of $S_fS_g^*$, as an 
operator on
$H_t$ belonging to $\ac_t$ (the commutant of $\PSLtwoZ$) is (up to a
normalization constant) $f(z)g(\xi)\zxi^k$.

In particular this shows that the above formula $(2.1)$ is explained by the trace
formula $\tau_{\ac_t}(k)=\left( 1/\area F\right) \int_F 
k(\overline
z,z)\,\d\nu_0(z)$, applied to the operator $k=S_fS_g^*$.

The role of the factor $\zxi^k$ is to make the function
$\overline {f(z)}g(\xi)\zxi^k$ diagonally $\PSLtwoZ$-invariant.
It is easy to observe that $S_f^*S_g$ is the Toeplitz operator on $H_t$
with symbol $\overline {f(z)}g(z)(\Im z)^k$.
Note that, to form $S_fS_g^*$, we have the restriction $k<t-1$, because $S_g^*$
has to map $H_t$ into a space $H_{t-k}$ that makes sense.

We   observe that the symbol of $S_fkS_g^*$ for an operator $k$
on $H_{t+k}$ of symbol $k=k(\overline {z},\xi)$ is
$${ c_{t-p}\over{ c_t}}\overline {f(z)}g(\xi)\zxi^{p}k(\overline {z},\xi),$$
if $f$, $g$ are automorphic forms of order $p$.

This also explains  the occurrence of operators of multiplication
  with symbols $\Phi(\overline {z},\xi)$ on
the space $L^2(\ac_t)$, in this setting.
In the terminology of the Appendix in the previous section those are
the Toeplitz operators with analytic symbol $\Phi(\overline {z},\xi)$,
a diagonally $\PSLtwoZ$-invariant function.
In the present setting, to get a bounded operator, we map $L^2(\ac_{t+k})$
into $L^2(\ac_t)$, by multiplying by
$\overline {f(z)}g(\xi)\zxi^k$.
In this section we will analyze the derivatives of a family of such 
operators.

Let $\Delta (z)$ be the unique automorphic form for $\PSLtwoZ$ in dimension
12 (this is the first order for which there is a nonzero space of
automorphic forms).

We rescale this form by considering the normalized function
$\Delta_1= \Delta / c$, where the constant $c$ is
chosen so that
$$\sup_{z\in \H}\left| \Delta_1(z)\right| ^2(\Im z)^{12}\leq 1.$$
{\it In the sequel we will omit the subscript $1$ from $\Delta$.}
This gives that the norm  $\|S_{\Delta}\|$, as an operator from $H_t$
into $H_{t+12}$, is bounded by $1$.

As $\Delta $ is a nonzero analytic function
  on the upper half-plane, one can choose an analytic branch for $\ln\Delta$.
Consider the
$\Gamma $-invariant function
$$\varphi(\overline {z},\xi)=\ln\overline {\Delta (z)}+\ln \Delta (\xi)+
12\ln\zxi,$$
which we also write as
  $$\varphi(\overline {z},\xi)=\ln(\overline {\Delta (z)}\Delta (\xi)\cdot
\zxi^{12}).$$

Defining $\pi_t(\gamma )$ for $\gamma $ in $\PSLtwoZ$  involves a 
choice of a branch for
$\ln (cz+d)$, $\gamma=\left(\matrix{a&b\cr
\noalign{\smallskip}
c&d\cr }\right)$. We define $\pi_t(\gamma )$, $\gamma \in 
\PSLtwoZ$,
by using the factor  $(cz+d)^{-t}$ corresponding to the following  choice of
the logarithm for $\ln(cz+d)$:

$$\ln \Delta (\gamma ^{-1}z)-\ln \Delta (z)=\ln (cz+d),$$
$\displaystyle z\in \H,\
\gamma=\left(\matrix{a&b\cr
\noalign{\smallskip}
c&d\cr }\right)$
in $\PSLtwoZ$.

By making this choice for $\pi_t$ restricted to $\Gamma $, we do not 
change the
algebra $\ac_t$, which is the commutant of $\{\pi_t(\Gamma )\}$, but 
we have the following.

  With the above choice for $\ln (cz+d)$ and thus for $\pi_t(\gamma )$,
$\gamma \in \PSLtwoZ$, and for any $\varepsilon >0$,
we have that $S_{\Delta ^\varepsilon} $ is
a bounded operator between $H_t$ and $H_{t+12\varepsilon }$, which 
intertwines $\pi_t$
and $\pi_{t+12\varepsilon }$ for all $t>1$, $\varepsilon >0$.

In the following lemma we make the symbol computation for operators of the form
$S_fS_g^*$. Recall that $H_t$, $t>1$ is the Hilbert space of analytic 
functions on $\H$
that are square-summable under 
$\d\nu_t=(\Im z)^{t-2}\,\d\overline{z}\,\d z$.\vskip6pt

{\sm Lemma} 2.1.\kern.3em
{\it Let $f,g$ be analytic functions on $\H$,
$k$ a strictly positive integer, and $t>1$. Assume
that $M_f=\sup_z\left| f(z)\right| ^2\Im z^k$,
$M_g=\sup_z\left| g(z)\right| ^2\Im z^k$ are finite quantities.

Let $S_f$, $S_g$ be the multiplication operators from $H_t$ into $H_{t+k}$ by
the functions $f,g$. Then $S_f$, $S_g$ are bounded operators of norm at most
$M_f$, $M_g$ respectively.

Moreover, the symbol of $S_fS_g^*\in B(H_t)$ is given by the formula
$${ c_{t-k}\over{ c_t}}\overline{f(z)}f(\xi)
\zxi^k.$$
}\vskip6pt

{\it Proof.} Before starting the proof we'll make the following 
remark that should
explain the role of the constant $c_t$ ($ = 
(t-1)/ 4\pi$)
in this computation.\vskip6pt

{\it Remark.}
The quantity
$c_t$ is a constant that appears due to the 
normalization in the
definition of $H_t$, where we have chosen
$$\|f\|_{H_t}^2=\int_{\H}\left| f(z)\right| ^2(\Im z)^{t-2}\,\d\overline{z}\,\d z.$$
Consequently  the reproducing vectors $e_z^t$
(defined by $\langle f,e_z^t\rangle =f(z)$, $f\in H_t$, $z\in \H$) are given
by the following formula [Ba, Mi]:
$$e_z^t(\xi)=\langle e_z^t,e_{\xi}^t\rangle ={ c_t\over{ 
\zxi^t}},\qquad z,\xi\in \H.$$

Consequently the normalized symbol of an operator $A$ in $B(H_t)$ is 
given by the
formula $k_A(\overline {z},\xi)=
\langle Ae_z^t,e_{\xi}^t\rangle / \langle 
e_z^t,e_{\xi}^t\rangle $,
$z,\xi$ in $\H$.

In the product formula we have that the symbol $k_{AB}(\overline {z},\xi)$
of the product of two operators $A,B$ on $H_t$ with symbols $k_A$, $k_B$ is
given by
$$\langle e_z^t,e_{\xi}^t\rangle k_{AB}(\overline {z},\eta)=
\langle ABe_z^te_{\xi}^t\rangle =
\int_{\H}\langle Ae_z^t,e_{\eta}^t\rangle \langle 
Be_{\eta}^t,e_{\xi}^t\rangle
\,\d\nu_t(\eta).$$

Thus
$$\eqalign{&k_{AB}(\overline {z},\xi)\cr
&\qquad
=\displaystyle \langle e_z^t,e_{\xi}^t\rangle \int_{\H}k_A(\overline {z},\eta)
\langle e_z^t,e_{\eta}^t\rangle k_B(\overline {\eta},\xi)
\langle e_{\eta}^t,e_{\xi}^t\rangle \,\d\nu_t(\eta)\cr
\noalign{\medskip}
&\qquad
=\displaystyle {\zxi^t\over c_t}\int _{\H}k_A(\overline {z},\eta)
{ c_t\over{\zueta^t}}
k_B(\overline {\eta},\xi){ c_t\over{ \etxi^t}}\,\d\nu_t(\eta)\cr
\noalign{\medskip}
&\qquad
=\displaystyle c_t\int _{\H}k_A(\overline {z},\eta)
k_B(\overline {\eta},\xi)[\overline{z},\eta,\overline{\eta},\xi]
\,\d\nu_0(\eta).}$$

This accounts for the constant $c_t$ that occurs in front of the 
product formula
(otherwise if we proceed as in [Ba] and include the constant $c_t$ in 
the measure
$\d\nu_t$, the constant will still show up in the product formula). \vskip6pt

In the proof
of the lemma we use the following observation.\vskip6pt

{\sm Observation} 2.2.\kern.3em
{\it Let $f,H_t,S_f$ be as in the statement of 
Lemma~{\rm 2.1.} Let $e^t_z$,
$e_z^{t+k}$ be the evaluation vectors at $z$, in the spaces $H_t$ 
and $H_{t+k}$.
Then
$$S_f^*e_z^{t+k}=f(z)e_z^{t},\qquad z\in \H.$$
}\vskip6pt

{\it Proof.} Indeed, since we will prove the boundedness of $S_f$, we 
can check this by
evaluating on a vector $g$ in $H_t$.
We have
$$\langle S_f^*e_z^{t+k},g\rangle _{H_t}=\langle 
e_z^{t+k},(S_f)g\rangle _{H_{t+k}}=
\langle e_z^{t+k},fg\rangle =\langle \overline{fg,e_z^{t+k}}\rangle =
\overline{ fg(z)}.$$

On the other hand:
$$\langle \overline{f(z)},e^t_z,g\rangle_{H_t}=\overline{f(z)}
\langle e^t_z,g\rangle_{H_t}=
\overline{f(z)}\langle \overline{g,e^t_z}\rangle_{H_t}=
\overline{f(z)}\overline{g(z)}.$$

This shows the equality of the two vectors.\hfill\qed\vskip6pt

We can now go on with the 
proof of Lemma~2.1.
It is obvious that $S_f,S_g$ are unbounded operators of norms $M_f,M_g$.
Indeed, for $S_f$ we have that
$$\eqalign{\|S_fg\|^2_{H_{t+k}}=&\displaystyle 
\int_{\H}\left| (S_fg)(z)\right| \,\d\nu_{t+k}(z)\cr
\noalign{\medskip}
=&\displaystyle \int_{\H}\left| (fg)(z)\right| ^2\,\d\nu_{t+k}(z)\cr
\noalign{\medskip}
=&\displaystyle \int_{\H}\left| f(z)\right| ^2\left| g(z)\right| ^2
(\Im z)^k(\Im z)^{t-2}\,\d\overline{z}\,\d z\cr
\noalign{\medskip}
=&\displaystyle \int_{\H}\left| g(z)\right| ^2(\left| f(z)\right| ^2(\Im z)^k)\,\d\nu_{t}(z)\cr
\noalign{\medskip}
\leq &\displaystyle M_f\int_{\H}\left| g(z)\right| ^2\,\d\nu_{t}(z).}$$

Hence $\|S_f\|\leq M_f$.

To prove the second assertion, observe that the symbol $k(\overline {z},\xi)$
of $S_fS_g^*$, as an operator on $H_t$, is given by the following formula:
$$\eqalign{k(\overline {z},\xi)=
&\displaystyle { \langle S_fS_g^*e^t_z,e^t_{\xi}\rangle _{H_t}\over{
\langle e^t_z,e^t_{\xi}\rangle _{H_t}}}\cr
\noalign{\medskip}
=&\displaystyle { \langle S_g^*e^t_z,S_f^*e^t_{\xi}\rangle \over{
\langle e^t_z,e^t_{\xi}\rangle _{H_t}}}\cr
\noalign{\medskip}
=&\displaystyle { \overline{g(z)}\langle e^{t-k}_z,
\overline{f(\xi)}e^{t-k}_{\xi}\rangle _{H_{t-k}}\over{
\langle e^t_z,e^t_{\xi}\rangle _{H_t}}}\cr
\noalign{\medskip}
=&\displaystyle \overline{g(z)}f(\xi){ \langle e^{t-k}_z,
e^{t-k}_{\xi}\rangle _{H_{t-k}}\over{
\langle e^t_z,e^t_{\xi}\rangle }}\cr
\noalign{\medskip}
=&\displaystyle \overline{g(z)}f(\xi){ \displaystyle { c_{t-k}\over{ 
\zxi^{t-k}}}
\over{ \displaystyle { c_t\over{ \zxi^t}}}}\cr
\noalign{\medskip}
=&\displaystyle { c_{t-k}\over{ c_t}}\overline{g(z)}f(\xi)\zxi^k,\qquad
z,\xi\hbox{ in }\H.}$$

This also works also for $k$ not an integer (as $\ln \zxi$ is chosen
once for all).\hfill\qed\vskip6pt

Let us finally note that the same arguments might be used to prove the
following more general statement.\vskip6pt

{\it Remark.} Let $f,g$ be analytic functions as in the statement of 
the lemma, and
let $k$ be an operator in $\ac_t$. Then $S_fkS_g^*$, which belongs to 
$\ac_{t+k}$
(if we think of $S_f,S_g$ as bounded operators mapping $H_t$ into 
$H_{t+k}$), has
the following symbol:
$$\displaystyle { c_{t}\over{ c_{t+k}}} f(\xi)\overline{g(z)}
\zxi^kk(\overline {z},\xi).$$
To show this, we have to evaluate
$$\eqalign{\displaystyle { \langle 
S_fkS_g^*e^{t+k}_z,e^{t+k}_{\xi}\rangle _{H_{t+k}}\over{
\langle e^{t+k}_z,e^{t+k}_{\xi}\rangle _{H_t}}}&=
{ f(\xi)\overline{g(z)}\langle 
ke^t_z,e^t_{\xi}\rangle _{H_t}\over{
\langle e^{t+k}_z,e^{t+k}_{\xi}\rangle _{H_{t+k}}}}\cr
\noalign{\medskip}
&=\displaystyle f(\xi)\overline{g(z)}{ \langle 
ke^{t}_z,e^{t}_{\xi}\rangle _{H_{t}}\over{
\langle e^t_z,e^t_{\xi}\rangle _{H_t}}}\cdot
{ \langle e^{t}_z,e^{t}_{\xi}\rangle _{H_{t}}\over{
\langle e^{t+k}_z,e^{t+k}_{\xi}\rangle _{H_{t+k}}}}\cr 
\noalign{\medskip}
&=\displaystyle {c_{t}\over{ c_{t+k}}} f(\xi)\overline{g(z)}
\zxi^kk(\overline {z},\xi).}$$\vskip6pt

In the next lemma we will deduce a positivity condition for kernels 
of operators
that occur as generators of parametrized families $S_{f^\varepsilon} 
S^*_{g^\varepsilon }$,
where $f,g$ are supposed to have a logarithm on $\H$.\vskip6pt

{\sm Lemma} 2.3.\kern.3em
{\it
  Assume that $f$ is a function as in Lemma~{\rm 2.1.} $f$ is analytic
on $\H$ and we assume $M_f=\sup_{z\in 
\H}\left| f(z)\right| ^2(\Im z)^k$ is
less than $1$.

Assume that $f$ is nonzero on $\H$, and choose a branch  for
$\ln f$ and hence for $f^{\varepsilon }$, $\varepsilon $ being 
strictly positive.

Let $\varphi (\overline {z},\eta)$ be the function
$\overline{\ln f(z)}+\ln f(\xi)+k\ln \zxi$ and use
this as a choice for
$\ln [\overline{f(z)}f(\xi)\zxi^k]=\varphi (\overline {z},\xi)$.

Then for all $\varepsilon > 0$ the  kernel
$$k_{\varphi }(\overline {z},\eta)=k_{\varphi ,t,\varepsilon 
}(\overline {z},\xi)=
\varphi ^{\varepsilon }(\overline {z},\xi)\left[
{ c_{t-k\varepsilon} \over{ c_t}}\ln \varphi -k{ c_t'\over{ c_t}}\right]$$
is nonpositive in the sense of $\ac_t$, that is,
$\displaystyle { k_{\varphi }(\overline {z}_i,z_j)\over{ \zizj^t}}$
is a nonpositive matrix
for all choices of $N\in \N$, $z_1,z_2,\ldots ,z_N\in \N$.
}\vskip6pt

{}{\it Proof.} By the choice we just made it is clear that the norm of
{}the operator
$S_{f^\varepsilon} $ is always less than $1$. We will also denote by
$S^t_{f^\varepsilon} $ the corresponding operators, which act as a contraction
from $H_t$ into $H_{t+k\varepsilon }$.

Consider the following operator-valued functions, with values in $H_t$:
$$f(\varepsilon )=
S_{f^\varepsilon }^{t-k\varepsilon }\left(S_{f^\varepsilon 
}^{t-k\varepsilon }\right)^*.$$

Obviously the symbol of $f(\varepsilon )$ is
  $\left( c_{t-k\varepsilon } / c_t\right) \varphi^\varepsilon 
(\overline z,\xi)$,
and moreover, $f(0)=1,\ f(\varepsilon )$ is a decreasing map because for
$0\leq \varepsilon \leq \varepsilon '$, we have that
$$
\displaystyle f(\varepsilon')=S_{f^{\varepsilon '}}^{t-k\varepsilon '}
\left(S_{f^{\varepsilon '}}^{t-k\varepsilon '}\right)^*
=S_{f^\varepsilon }^{t-k\varepsilon }
\left[S_{f^{\varepsilon '-\varepsilon }}^{t-k\varepsilon '}
\left(S_{f^{\varepsilon '-\varepsilon }}^{t-k\varepsilon
'}\right)^*\right]\left(S_{f^\varepsilon }^{t-k\varepsilon }\right)^*
$$

But the operator in the middle has norm less than $1$, and hence we get that
$$f(\varepsilon ')\leq S_{f^\varepsilon }^{t-k\varepsilon }
\left(S_{f^\varepsilon }^{t-k\varepsilon }\right)^*=f(\varepsilon ).$$

Fix $N$, and  $z_1,z_2,\ldots, z_N$ in $\H$.
Then (since the corresponding operators form a decreasing family)
$$g(\varepsilon )=
\displaystyle \left[{f(\varepsilon )(\overline{z_i},z_j)\over 
\zizj^t}\right]_{i,j}$$
is a decreasing family of matrices, and $g(0)=\Id $. Hence 
$g'(\varepsilon )$ must
be a negative (nonpositive) matrix.
  Note that $f(\varepsilon )(\overline{z_i},z_j)=
\left( c_{t-k\varepsilon } / c_t\right) \varphi^\varepsilon 
(\overline{z_i},z_j)$.

But $g'(\varepsilon )$ has exactly the formula stated above, that is,
$$g'(\varepsilon )=\displaystyle {\varphi^\varepsilon (\overline{z_i},z_j)
\left[\ln \varphi(\overline{z_i},z_j)
\displaystyle {c_t-k\varepsilon \over c_t}-k{c'_t\over c_t}\right]\over
{\zizj^t}}$$ This completes the proof.\hfill\qed\vskip6pt

By collecting the terms together, and since
$c'_t / c_t = 1 / (t-1)$,
we obtain, for all $\varepsilon \geq 0$, the following\vskip6pt

{\sm Lemma} 2.4.\kern.3em
{\it
  With the notations from the previous lemma, for all $\varepsilon 
\geq 0$, the kernel
$$k_\varphi=k_{\varphi,\varepsilon ,t}=
\varphi^\varepsilon \left[\ln \varphi-\displaystyle {1\over
{t-1-k\varepsilon }}\right]$$  is nonpositive in $\ac_t$.
  Precisely this means that for all choices of $N$ in $\N$ and
$z_1,z_2,\ldots,z_N$ in $\H$ we have
that $$\displaystyle {k_\varphi (\overline{z_i},z_j)\over {\zizj^t}}$$ is
a nonpositive matrix.
}\vskip6pt

The following observation will be used later in the proofs.\vskip6pt

{\sm Observation} 2.5.\kern.3em
{\it For any $s>1$,
the   identity
$$\tau_{\ac_s}(S_g^*S_g)=
\displaystyle{
c_s\over{
c_{s+k}}}\tau_{\ac_t}(S_g S_g^*),$$
  holds true for any automorphic form $g$ of order $k$.
}\vskip6pt

{\it Proof.}
We have that $(S_g^*S_g)$ is the Toeplitz operator (on $H_s$) with
symbol $$\left| g(z)\right| ^2(\Im z)^k.$$
  Hence the trace
$\tau_{\ac_s}( S_g^* S_g)$ is $\displaystyle {1\over {\area(F)}}
\int_{F}\left| g\right| ^2(\Im z)^k \,\d \nu_0(z).
$

On the other hand the symbol of $(  S_g S_{g^*})$
(which  is viewed here
as an operator on $H_{s+k}$) is equal to
$$(z,\xi)\rightarrow
\displaystyle{
\left\langle S_{g^*} e_z^{s+k}, S_{g^*} e_{\xi}^{s+k}\right\rangle
\over{
\left\langle e_z^{S+k},e_{\xi}^{e+k}\right\rangle }}=
\displaystyle{
c_s\over{
c_{s+k}}}
\overline{g(z)}g(\xi)\zxi^k,
$$
and hence the trace of the latter symbol
is
$$\displaystyle{
c_s\over{
c_{s+k}}}
\displaystyle {1\over {\area (F)}}
\int_{F}\left| g(z)\right| ^2(\Im z)^k \,\d \nu_0(z).
\eqno{\qed}$$

\vfill\eject

\centerline{\S 3. {\ninebf Properties of the (unbounded) 
multiplication maps by}}
\centerline{\ninebf  $\ln [\overline{\Delta (z)}\Delta 
(\xi)\zxi^{12}]$ on different
spaces of kernels}\vskip12pt
Let $\varphi (\overline {z},\xi)=
\ln (\overline{\Delta (z)}\Delta (\xi)\zxi^{12})$.
In this section we want to exploit the negativity properties of the kernels
$$\varphi ^{\varepsilon/12 }
\left(\displaystyle{1\over 12} \ln \varphi -{ 1\over{ t-1-\varepsilon 
}}\right).$$
By $\odot$ we denote the operation
of pointwise multiplication of symbols. It is the analogue of Schur
  multiplication on matrices
or on a group algebra. When no confusion is possible we will omit the symbol
$\odot$ and just replace it by $\,\cdot \,$.

We want to draw conclusions on the properties of the  multiplication 
maps,  defined on a
suitable dense subspace of $L^2(\ac_t)$, by the formula
$$\Lambda _{\varepsilon }(k) = k\odot
  \displaystyle \left(\varphi ^{\varepsilon/12 }
\left[
  {1\over 12}\ln \varphi -{ 1\over{ t-1-\varepsilon }} 
\right] \right).$$

For functions $k(\overline {z},\eta)$ on $\H\times \H$ that are 
positive, but do
not necessary represent a positive operator, we will introduce the following
definition.\vskip6pt

{\sm Definition} 3.1.\kern.3em
{\it
  A function $k(\overline {z},\eta)$ on $\H\times \H$
that is analytic for $\eta$ and antianalytic for $z$ will be called 
{\rm positive for
$\ac_t$} if the following
  matrix
$$\displaystyle  {\left[k(\overline {z}_i,z_j)\over{ \zizj^t} 
\right]^n_{i,j=1}}$$
is positive for all choices of $N\in \N$ and
$z_1,z_2,\ldots ,z_N$ in $\H$.

The space of such kernels will be denoted by ${\cal S}_t$.
}\vskip6pt

The following remark is a trivial consequence of
  the fact that the Schur product of two positive matrices is positive,
and a consequence of the description
  for positivity of kernels of operators in $\ac_t$ given in Section 1.\vskip6pt

{\sm Proposition} 3.2.\kern.3em
{\it For all numbers $r,s>1$,
the vector space $(\ac_r)_+ \odot {\cal S}_s$ is contained in
${\cal S}_{r+s}$ and $(\ac_r)_+\subseteq {\cal S}_r$.}\vskip6pt

{\it Proof.} Just observe that in fact ${\cal S}_s\odot{\cal S}_r$
  is contained in ${\cal S}_{s+r}$.\hfill\qed\vskip6pt

The problem that we address in this section comes from the fact that 
the operator
$$-\Lambda _{\varepsilon ,r,s}(k)=
k\odot\varphi ^{\varepsilon/12 }
\left[-{1\over 12}\ln \varphi+ { 1\over{ r-1-\varepsilon }}\right]$$
maps $k\in(\ac_s)_+$ into ${\cal S}_{r+s}$. Also
$\Lambda _{\varepsilon ,r,s}(k)(\overline{z},z)$ is integrable on $F$,
  so it is tempting to infer
that
$\Lambda _{\varepsilon ,r,s}(k)$ belongs to $L^1(\ac_{r+s})$. In fact, 
we conjecture that
  a kernel $k(\overline {z},\eta)$ in ${\cal S}_t$ that is also 
diagonally integrable
on $F$ corresponds to an element in $L^1(\ac _t)_t$. Since we are unable to
  prove the conjecture directly, we will use monotonicity properties for the
derivatives of $\varphi^{\varepsilon }$.

If no constants were involved, we would
simply say that $\varphi^{\varepsilon }(-\ln \varphi )$ is the increasing
limit of the derivatives, since the second derivative would be negative.
This doesn't hold exactly, but the constants involved are small enough and
have a negligible effect on the previous line of reasoning. This is 
done in the
following lemma.\vskip6pt

{\sm Lemma} 3.3.\kern.3em {\it
Let $k$ be a positive kernel in $\ac_s$. Fix $v>1$ and let
$\varepsilon >0$ be small enough. Consider the following elements in 
$\ac_{s+\varepsilon }$
defined by the kernels
$$\lambda _{\varepsilon ,v,s}(k)(\overline {z},\xi)=
\left[{ v-1-\varepsilon \over{ v-1}}\varphi ^{\varepsilon/12 }
(\overline {z},\xi)\right]k(\overline {z},\xi).$$

Note that up to a multiplicative constant $\lambda _{\varepsilon 
,v,s}(k)$ is the
kernel of $S_{\Delta_{}^{\!\varepsilon/12}}
kS_{\Delta_{}^{\!\varepsilon/12}}^* $ in $\ac_{s+\varepsilon }$.
Let $\widetilde{\lambda }_{\varepsilon ,v,s}(k)$ be the image {\rm (}through
$\Psi_{v+2s,v+\varepsilon }${\rm )} of this kernel in $\ac_{v+2s}$.

Then $\widetilde{\lambda }_{\varepsilon ,v,s}(k)$ is a decreasing family of
positive kernels representing elements in $\ac_{v+2s}$, and there 
exists a negative element
$M(k)=M_{\varepsilon ,v,s}(k)$ in $-L^1(\ac_{v+2s})_+$ such that 
$M(k)$ is the derivative
with respect to $\varepsilon $:
$$\displaystyle M(k)={ \d\;\over{ \d\varepsilon }}\widetilde{\lambda 
}_{\varepsilon ,v,s}(k).
\eqno{(3.1)}$$

The derivative is computed in the  strong convergence topology,
  on a dense domain $\dc\subseteq
H_{v+2s}$ affiliated with $\ac_{v+2s}$.

The symbol of $M_{\varepsilon ,v,s}(k)$ as an operator in $H_{v+2s}$ is equal
to
$$\Lambda _{\varepsilon ,v,s}(k)(\overline {z},\xi)=
{ v-1-\varepsilon \over{ v-1}}k(\overline {z},\xi)\varphi ^{\varepsilon/12 }
\left[{1\over 12}\ln \varphi-{ 1\over{ v-1-\varepsilon }}\right].$$
}\vskip6pt

{\it Proof.}  For simplicity of the proof we will use the notation
$\varphi_1=\varphi^{1/12}$. We prove first that the family
$\lambda _{\varepsilon ,v,s}(k)$ is a decreasing family in
$\ac_{v+s+\varepsilon }$ and hence in $\ac_{v+2s}$.

Indeed, $\left( ( v-1-\varepsilon ) / ( v-1 ) \right) \varphi_1 ^{\varepsilon }
(\overline {z},\xi)$ is a decreasing family of operators in $\ac_v$,
and hence by  Proposition~3.2 it follows that
$$\displaystyle { v-1-\varepsilon \over{ v-1}}\varphi_1 ^{\varepsilon }
(\overline {z},\xi)k(\overline {z},\xi)$$ is a decreasing family
in ${\cal S}_{v+s}$, and hence in ${\cal S}_{v+s+\varepsilon }$. 
Since we  know that
these operators are already bounded in $\ac_{v+\varepsilon }$, the 
first part of
the statement follows immediately.

Denote by $G(\varepsilon )=G(\varepsilon )(\overline {z},\xi)$ the 
kernel represented by
$$\displaystyle { v-1-\varepsilon \over{ v-1}}\varphi_1 ^{\varepsilon }
(\overline {z},\xi)k(\overline {z},\xi),$$
which represents therefore a (decreasing) family in $\ac_{v+s+\varepsilon }$
and hence in $\ac_{v+2s}$. Fix $\varepsilon _0>0$ and let
$$g_\varepsilon (\overline {z},\xi)={ G(\varepsilon )(\overline {z},\xi)-
G(\varepsilon_0 )(\overline {z},\xi)\over{ \varepsilon -\varepsilon _0}}.$$

Then $g_{\varepsilon }$ is a negative (nonpositive) element in
$\ac_{v+2s}$. We want to find a formula for $g_{\varepsilon 
'}-g_{\varepsilon }$.
Obviously when $\varepsilon $ converges to $\varepsilon _0$, the kernel
$g_{\varepsilon }$ converges (at least pointwise) to the kernel
$\Lambda_{\varepsilon _0,v,s}(k)(\overline {z},\xi)$.

It is elementary calculus to find
 for
$H_{\varepsilon ',\varepsilon }(\overline {z},\xi)=g_{\varepsilon 
'}(\overline {z},\xi)-
g_{\varepsilon }(\overline {z},\xi)$
 the following pointwise expression:
$$H_{\varepsilon ,\varepsilon '}=(\varepsilon -\varepsilon ')
\int_0^1\left(\int_0^1 tG''(\varepsilon (t,s))\,\d s\right)\,\d t
\eqno{(3.2)}$$
where $\varepsilon (t,s)=s[(1-t)\varepsilon _0+t\varepsilon 
]+(1-s)[(1-t)\varepsilon _0+
t\varepsilon ']$ belongs to the interval determined by
$\varepsilon ,\varepsilon ',\varepsilon _0$.

This formula holds at the level of kernels (that is, by evaluating 
both sides on
any given points $z,\xi\in \H$).

On the other hand, because
$$G(\varepsilon )={ 1\over{ v-1}}
[(v-1-\varepsilon )\varphi_1 ^{\varepsilon }]k,$$
we may compute immediately that
$$\eqalign{G'(\varepsilon )=&\displaystyle
{ 1\over{ v-1}}[(v-1-\varepsilon )
\varphi_1 ^{\varepsilon }\ln \varphi_1 -\varphi_1 ^{\varepsilon }]k,\cr
\noalign{\medskip}
G''(\varepsilon )=&\displaystyle
{ 1\over{ v-1}}[(v-1-\varepsilon )\ln ^2\varphi_1
-2\ln \varphi_1 ]\varphi_1^\varepsilon \cdot k.}$$

Furthermore we have the following expression for $G''(\varepsilon )$:
$$\eqalign{G''(\varepsilon )=&\displaystyle
{ v-1-\varepsilon \over{ v-1}}\varphi_1 ^{\varepsilon }\cdot k
\left[\ln ^2\varphi_1 -{ 2\over{ v-1-\varepsilon }}\ln \varphi_1 \right]\cr
\noalign{\medskip}
=&\displaystyle { v-1-\varepsilon \over{ v-1}}\left\{k\left[
\varphi_1 ^{\varepsilon /2}\left(-\ln \varphi_1 +
{ 1\over{ v-1-\varepsilon }}\right)\right]^2-
{ k\varphi_1^{\varepsilon }\over{ (v-1-\varepsilon )^2}}\right\}.}$$

Thus we obtain further that
$$G''(\varepsilon )=
{ v-1-\varepsilon \over{ v-1}}k\left[
\varphi_1 ^{\varepsilon /2}\left(-\ln \varphi_1 +
{ 1\over{ v-1-\varepsilon }}\right)\right]^2-
{ k\varphi_1^{\varepsilon }\over{ (v-1-\varepsilon )(v-1)}}.$$

But because of the previous Proposition 3.2 and
Lemma 2.4, we have that
$$\left[\varphi_1 ^{\varepsilon /2}\left(-\ln \varphi_1 +
{ 1\over{ v-1-\varepsilon }}\right)\right]^2=
\left[\varphi_1 ^{\varepsilon /2}\left(-\ln \varphi_1 +
{ 1\over{ v-\varepsilon /2-1-\varepsilon /2}}\right)\right]^2$$
represents the square of an element
$$\varphi_1 ^{\varepsilon /2}\left(-\ln \varphi_1 +{ 1\over{ 
v-1-\varepsilon }}\right)$$
in ${\cal S}_{v-\varepsilon /2}$.
The square of the above element consequently belongs to
${\cal S}_{2v-\varepsilon }$.

Hence $${\cal R}(\varepsilon )=k\displaystyle
\left[\varphi_1 ^{\varepsilon /2}\left(-\ln \varphi_1 +
{ 1\over{ v-1-\varepsilon }}\right)\right]^2,$$
as a kernel, belongs to
${\cal S}_{s+2v-\varepsilon }\subseteq {\cal S}_{s+2v}$.

In conclusion, we have just verified that
$$\displaystyle G''(\varepsilon )={\cal R}(\varepsilon )-{ k\varphi_1 
^{\varepsilon }\over{
(v-1)(v-1-\varepsilon )}},$$
where ${\cal R}(\varepsilon )$ belongs to ${\cal S}_{2v+s}$.

Moreover, it is obvious that
$$Q(\varepsilon )={ k\varphi_1 ^{\varepsilon }\over{
(v-1)(v-1-\varepsilon )}}$$
is a bounded element in $\ac_{v+2s}$, and that $Q(\varepsilon )$ is 
consequently
bounded by a constant $C$, independent of all the variables $v,s,\varepsilon $:
$$Q(\varepsilon )\leq C\cdot\Id \hbox{\qquad in }\ac_{2v+s},$$
  and hence
$$Q(\varepsilon )\leq C\cdot\Id \hbox{\qquad in }{\cal S}_{2v+s}.$$

We put this into the integral formula for
$$H_{\varepsilon ,\varepsilon '}=g(\varepsilon )-g(\varepsilon ')=
{ G(\varepsilon )-G(\varepsilon _0)\over{ \varepsilon -\varepsilon _0}}-
{ G(\varepsilon') -G(\varepsilon _0)\over{ \varepsilon '-\varepsilon _0}}$$
to obtain that $H_{\varepsilon ,\varepsilon '}$ is of the form
$(\varepsilon -\varepsilon ')\left[\rc-\qc\right]$ where $\rc$ 
belongs to ${\cal S}_{2v+s}$
and $\qc$ belongs to
$(\ac_{v+2s})_+ \hbox{ and } 0\leq \qc\leq C\cdot\Id .$

But then $\rc$ belongs to $\ac_{v+2s}\cap {\cal S}_{v+2s}$, and hence
$\rc \in (\ac_{v+2s})_+$.

Thus, in $\ac_{v+2s}$, we have (assuming 
$\varepsilon -\varepsilon '\geq 0$) that
$$H_{\varepsilon ,\varepsilon '}
\geq -(\varepsilon -\varepsilon ')\qc\geq -(\varepsilon -\varepsilon ')C,$$
so
$$H_{\varepsilon ,\varepsilon '}\geq -(\varepsilon -\varepsilon ')C,$$
and therefore
$$g_{\varepsilon }-g_{\varepsilon '}\geq C(\varepsilon '-\varepsilon ).$$
Hence for $\varepsilon \geq \varepsilon '\geq \varepsilon _0$ we have that
$g_{\varepsilon }+C_{\varepsilon }\geq g_{\varepsilon '}+C_{\varepsilon '}$ in
$\ac_{v+2s}$, for a fixed positive constant $C$.

Now recall that
$$g(\varepsilon )={ G(\varepsilon )-G(\varepsilon _0)\over{ 
\varepsilon -\varepsilon _0}}$$
and that $G(\varepsilon )$ itself was a decreasing family in $\ac_{v+2s}$,
so that $g(\varepsilon )$ are negative elements in $\ac_{v+2s}$.

Denote for simplicity $h(\varepsilon )=-g(\varepsilon )$.
Then what we just obtained is the following:

The operators $h(\varepsilon )$ are positive elements
in $(\ac_{v+2s})_+$. Moreover,
$h(\varepsilon )-C{\varepsilon }\leq h(\varepsilon ')-C{\varepsilon '}$ if
$\varepsilon \geq \varepsilon '$, i.e., $h(\varepsilon )-C\varepsilon $
is a decreasing family. By adding a big constant $k$ to $h(\varepsilon )$
we have that $K+h(\varepsilon )-C{\varepsilon }$ is a decreasing family
of positive elements in $(\ac_{s+2v})_+$.

Thus as $\varepsilon $ decreases to $\varepsilon _0$ we have that
$K+h(\varepsilon )-C{\varepsilon }$ is an increasing family
of positive operators in $\ac_{v+s}$.

Moreover, the trace of $h(\varepsilon )$ is equal to
$-\tau (g(\varepsilon ))$, which is
$$-\int_{F}{ \left[
{ v-1-\varepsilon \over{ v-1}}\varphi_1 ^{\varepsilon }(\overline {z},z)-
\varphi_1 ^{\varepsilon _0}(\overline {z},z){ v-1-\varepsilon 
_0\over{ v-1}}\right]
k(\overline {z},z)
\over{ \varepsilon -\varepsilon _0}}\,\d\nu.$$
This integral converges (in $L^1(F\,\d\nu_0)$) to
$$\eqalign{\displaystyle &\int_{F}
\left.{ \d\;\over{ \d\varepsilon }}{ v-1-\varepsilon \over{ v-1}}
\varphi_1 ^{\varepsilon }(\overline {z},z)\right|_{\varepsilon =\varepsilon _0}
k(\overline {z},z)\,\d\nu_0(z)\cr
\noalign{\medskip}
&\qquad =-{ v-1-\varepsilon _0\over{ v-1}}\int_{F}\varphi_1 ^{\varepsilon_0 
}(\overline {z},z)
\left[\ln \varphi_1 -{ 1\over{ v-1-\varepsilon_0 }}\right]k(\overline 
{z},z)\,\d\nu_0(z)},$$
which is finite (the convergence is dominated here for
  example by $C\varphi_1 ^{\varepsilon '}$,
for some $\varepsilon '\leq \varepsilon_0 $).

Thus $K+h(\varepsilon )-C(\varepsilon )$ are an increasing family in $\ac_t$
(as $\varepsilon$ decreases to $ \varepsilon_0 $)
  and the supremum of the traces (in $L^1(\ac_t)$)
is finite. By Lesbegue's Dominated
  Convergence Theorem in $L^1(\ac_t)$, the limit of
$K+h(\varepsilon )-C(\varepsilon )$ exists in $L^1(\ac_t)$
and convergence is in the strong operator topology on a dense domain 
affiliated with $\ac_t$.\hfill\qed\vskip6pt

{\sm Corollary} 3.4.\kern.3em
{\it Let  $\Lambda_{\varepsilon ,v,s}(k)$ be 
the map, defined  in the
previous lemma, that associates to any positive $k$ in $\ac_{s}$ a 
positive element
in $\ac_{v+2s}$, whose kernel is given by the formula:
$$\Lambda_{\varepsilon ,v,s}(k)(\overline {z},\xi)={ v-1-\varepsilon 
\over{ v-1}}
k(\overline {z},\xi)\varphi ^{\varepsilon/12 }\left[-{1\over 12}\ln \varphi +
{ 1\over{ v-1-\varepsilon }}\right].$$
  Then $-\Lambda_{\varepsilon ,v,s}$ is
  a completely positive map
  from $\ac_{s}$ into $L^1(\ac_{v+2s})$.}\vskip6pt

{\it Proof.} From the previous lemma we know that 
$\Lambda_{\varepsilon ,v,s}(k)$
is well defined and belongs to $L^1(\ac_{v+2s})$. On the other hand
$\Lambda_{\varepsilon ,v,s}(k)$ is obtained by multiplication with a 
positive kernel
in ${\cal S}_v$, and hence (as in the proof of the complete positivity for
$\Psi_{s,t}$) we obtain that $\left[\Lambda_{\varepsilon 
,v,s}\left(k_{pq}\right)\right]_{p,q }$
is a positive in $M_N(\C)\otimes L^1(\ac_{v+2s})$, if $[k_{p,q}]$ is a
positive matrix in $M_N(\C)\otimes\ac_{s}$.\hfill\qed\vskip6pt

{\sm Corollary} 3.5.\kern.3em
{\it Let $\varepsilon _0>0$ and $t>3+\varepsilon _0$.
Let $\Lambda_{\varepsilon }$ be defined, on the space of all symbols $k$
representing operators in 
$\bigcup_{1<s<t-2-\varepsilon _0}\ac_{s}$,
by the formula $$\Lambda_{\varepsilon _0}(k)=\displaystyle
\left.{ \d\;\over{ \d\varepsilon }}
(k\odot\varphi ^{\varepsilon })\right|_{\varepsilon =\varepsilon _0}.$$
Note that the pointwise derivative of kernels is
$(\varphi ^{\varepsilon _0}\ln \varphi )\odot k$.

Then $\Lambda_{\varepsilon _0}(k)$ belongs to $L^1(\ac_{t})$, and moreover, the
derivative is valid in the sense of the strong operator topology on a 
dense domain
affiliated to $\ac_{t}$.

Fixing $1<s<t-2-\varepsilon _0$, there exists a sufficiently large constant $C$
{\rm (}depending on $s,t,\varepsilon _0${\rm )} such that
$-\left[\Lambda_{\varepsilon _0}+Ck\odot \varphi ^{\varepsilon }\right]$
{\rm (}and hence $-\left[\Lambda_{\varepsilon _0}+C\cdot\Id \right]${\rm )}
becomes a completely positive operator from $\ac_{s}$ into $\ac_{t}$.
}\vskip6pt

{\it Proof.} Because of the condition $s<t-2-\varepsilon _0$, we can always
find a constant $C$, by the previous lemma, such that the previous lemma
applies to
the operator
$\Lambda_{\varepsilon _0}+Ck\odot\varphi ^{\varepsilon }$.\hfill\qed\vskip6pt

{\sm Corollary} 3.6.\kern.3em
{\it Fix $t>3$. For every $1<s<t-2$ and for every $k$
in $\ac_{s}$ there exists an {\rm (}eventually unbounded\/{\rm )} operator
$\Lambda(k)$ {\rm (}of symbol multiplication by $\ln \varphi ${\rm )}
that is affiliated
with $\ac_{t}$, and there exists a dense domain $\dc$ in $H_t$ that 
is affiliated
with $\ac_{t}$, such that the derivative
$$\left.
{ \d\;\over{ \d\varepsilon }}\langle k\odot \varphi ^{\varepsilon 
}\xi,\eta\rangle _{H_t}
\right|_{\varepsilon =0}$$
exists for all $\xi,\eta$ in $\dc$ and is equal to
$$\langle \Lambda (k)\xi,\eta\rangle.$$

Moreover, there exists a constant $C$, depending only on $s,t$,
such that for any positive matrix $[k_{p,q}]_{p,q=1}^p$ in
$M_N(\ac_{s})_+$, the operator matrix
$$-\left[(\Lambda+C\cdot\Id )(k_{p,q})\right]_{p,q=1}^N$$
represents a positive operator, affiliated with $\ac_{t}$.
}\vskip6pt

{\it Remark.} The operator $k\odot \varphi ^{\varepsilon }$ appearing 
in the previous
statement is bounded. Indeed, modulo a multiplicative constant,
$k\odot \varphi ^{\varepsilon }$ is the symbol of $S_{\Delta^\varepsilon} k
S^*_{\Delta^\varepsilon} $. If $k\in \ac_{s}$, then
$S_{\Delta^\varepsilon} k
S^*_{\Delta^\varepsilon} $
belongs to $\ac_{s+12\varepsilon }$, and since $s<t$, by choosing 
$\varepsilon $
small enough, we can assume that $S_{\Delta^\varepsilon} k
S^*_{\Delta^\varepsilon} $
represents a bounded operator on $H_t$, and hence that the expression
$\langle k\odot \varphi ^{\varepsilon }\xi,\eta\rangle_{H_t}$ makes sense
for all $\xi,\eta$ in $H_t$.\vskip6pt

Before going to the proof of the statement of Corollary 3.6, we prove 
the following lemma
(which will be used in the proof of Corollary 3.6)
concerning the operator $k\odot \varphi ^{\varepsilon }$ and the 
range of the operator
$\Lambda_{\varepsilon,v,s }(k)$ defined in Lemma~3.3.\vskip6pt

{\sm Lemma} 3.7.\kern.3em
{\it With the notations from Lemma~{\rm 3.3,} let
$k$ be an operator in $\ac_{s}$, $v>1$, $\varepsilon >0$. Let
$M_{\varepsilon ,v,s}(k)$ be the derivative
{\rm (}with respect to $\varepsilon ${\rm ),} which
belongs to $L^1(\ac_{v+2s})$, of the decreasing family
$$\widetilde {\lambda }_{\varepsilon ,v,s}(k)={ v-1-\varepsilon \over{ v-1}}
\varphi ^{\varepsilon }(\overline {z},\xi)k(\overline {z},\xi).$$
Then the range and init space of the unbounded
  operator $M_{\varepsilon ,v,s}(k)$ are contained
{\rm (}and  dense\/{\rm )} in the closure of the range
of $S_{\Delta^\varepsilon} 
\subseteq H_{2v+s}$ {\rm (}more
precisely in closure of the range of $S^{2v+s-\varepsilon 
}_{\Delta^\varepsilon }${\rm ).}
}\vskip6pt

{\it Proof.} Indeed, by what we have just proved, $M_{\varepsilon ,v,s}(k)$
is the strong operator topology limit
(on a dense domain affiliated with the von Neumann algebra),
  as $\varepsilon'$  decreases to
$\varepsilon $, of the operators
$$\displaystyle { G_{\varepsilon '}(\overline {z},\xi)-
G_{\varepsilon }(\overline {z},\xi)\over{ \varepsilon '-\varepsilon
}}.$$
  Recall that $G_{\varepsilon '}(\overline {z},\xi)$
  was the symbol (modulo a multiplicative constant)
of the operator $S_{\Delta^{\varepsilon '}}k
S^*_{\Delta^{\varepsilon '}}$.

Then by applying
$( G_{\varepsilon '}-G_{\varepsilon '} ) / ( 
\varepsilon '-\varepsilon )$
to any vector $\xi$ in $H_t$, the outcome is already a vector in the 
closure of the range of
$S_{\Delta^\varepsilon} $.
This property is preserved in the limit. By selfadjointness the same 
is valid for the init
space.\hfill\qed\vskip6pt

We proceed now to the proof of Corollary~3.6.\vskip6pt

{\it Proof of Corollary}~3.6. We start by constructing first the 
domain $\dc$.
For $\varepsilon_0>0$ let  $\dc_{\varepsilon _0}\subseteq H_t$  be the range
of $(S^t_{\Delta^{\varepsilon _0}})^*$, considered as
  an operator from $H_{t+12\varepsilon _0}$
into $H_t$. $\dc$ will be the increasing union
(with respect to $\varepsilon_0$) of $\dc_{\varepsilon _0}$.

Let $B_{\varepsilon _0}$ be a right inverse, as an unbounded operator
  for the operator
$S_{\Delta^{\varepsilon _0}}$. Thus $B_{\varepsilon _0}$
  acts from a domain dense in the closure of
range $S^t_{\Delta^{\varepsilon _0}}$ into $H_t$. It is clear that 
$B_{\varepsilon _0}$
is an intertwiner affiliated with the von Neumann algebras $\ac_t$ and
$\ac_{t+12\varepsilon _0}$ (by  von Neumann's theory
of unbounded operators,
affiliated to a $\twoone$ factor [MvN]).

Thus, denoting by $P_{\varepsilon _0}$ the projection onto the 
closure of the range of
$S^t_{\Delta^{\varepsilon _0}}$ in $H_{t+12\varepsilon _0}$, the 
following properties
hold true:
$$(S^t_{\Delta^{\varepsilon _0}})B_{\varepsilon _0}=P_{\varepsilon _0.}$$
By taking the adjoint, we obtain
$$B^*_{\Delta_{\varepsilon _0}}\left(S^t_{\Delta^{\varepsilon _0}}\right)^*=
P_{\varepsilon _0}.$$

All compositions make sense in the algebra of unbounded operators 
affiliated with
$\ac_{t}$, and $\ac_{t+12\varepsilon _0}$.
On $H_{t+12\varepsilon _0}$, we let $M_{\varepsilon _0}(k)$ be the 
$L^1$ operator
given by Corollary 3.4, whose symbol is $$k\odot
\varphi ^{\varepsilon_0 }\ln \varphi ,$$
for $k$ in $\ac_{s}$.

We define $\Lambda_{\varepsilon _0}(k)$ by the following composition:
$$\Lambda_{\varepsilon _0}(k)={ c_{t+12\varepsilon _0}\over{ c_t}}
B_{\varepsilon _0}M_{\varepsilon _0}(k)B^*_{\varepsilon _0}.$$
We want to prove that $\Lambda_{\varepsilon _0}$ does not depend on 
$\varepsilon_0$.
Obviously (by [MvN]), the operator $\Lambda_{\varepsilon _0}(k)$
is affiliated with $\ac_{t}$.

Moreover, for $\xi,\eta$ in $\dc_{\varepsilon _0}$, which are thus
of the form
$$\xi=S^*_{\Delta^{\varepsilon _0}}\xi_1,\qquad
\eta=S^*_{\Delta^{\varepsilon _0}}\eta_1,$$
  for some $\xi_1$, $\eta_1$
in $H_{t+12\varepsilon _0},$
we have that
$$\langle \Lambda_{\varepsilon _0}(k)\xi,\eta\rangle=
\langle \Lambda_{\varepsilon _0}(k) S^*_{\Delta^{\varepsilon _0}}\xi_1,
S^*_{\Delta^{\varepsilon _0}}\eta_1\rangle _{H_t}.$$
This is equal to
$$\eqalign{ \displaystyle &{c_{t+12\varepsilon _0}\over{ c_t}}\langle
B_{\varepsilon _0}M_{\varepsilon _0}(k)B^*_{\varepsilon _0}
S^*_{\Delta^{\varepsilon _0}}\xi_1,S^*_{\Delta^{\varepsilon 
_0}}\eta_1\rangle _{H_t}\cr
\noalign{\medskip}
&\qquad =\displaystyle {c_{t+12\varepsilon _0}\over{ c_t}}\langle
P_{\varepsilon _0}M_{\varepsilon _0}(k)P_{\varepsilon _0}
\xi_1,\eta_1\rangle _{H_{t+12\varepsilon _0}}.}$$
Because of Lemma~3.7, we know that this is further equal to
$$\displaystyle {c_{t+12\varepsilon _0}\over{ c_t}}\langle
M_{\varepsilon _0}(k)\xi_1,\eta_1\rangle .$$

We use the above chain of equalities  to deduce that the definition of
$\Lambda_{\varepsilon _0}(k)$ is  independent of the choice of
$\varepsilon _0$.

  Indeed, assume  we  use another
$\varepsilon _0'$, which we  assume to be bigger than $\varepsilon _0$. Assume
$\xi=S^*_{\Delta^{\varepsilon _0'}}\xi_2$. This is further equal to
$S_{\Delta^{\varepsilon _0}}^*S^*_{\Delta^{(\varepsilon 
_0'-\varepsilon _0)}}\xi_2$.

Then, by redoing the previous computations we arrive to the term
$$\displaystyle {c_{t+12\varepsilon _0'}\over{ c_t}}\langle
M_{\varepsilon _0'}(k)\xi_2,\eta_2\rangle .$$

But on the other hand in this situation
$$\eqalign{\langle M_{\varepsilon _0}(k)\xi_1,\eta_1\rangle =&
\displaystyle {c_{t+12\varepsilon _0}\over{ c_t}}
\langle M_{\varepsilon _0}(k)S^*_{\Delta^{(\varepsilon 
_0'-\varepsilon _0)}}\xi_2,
S^*_{\Delta^{(\varepsilon _0'-\varepsilon _0)}}\eta_2\rangle \cr
\noalign{\medskip}
\displaystyle =&\displaystyle {c_{t+12\varepsilon _0}\over{ c_t}}
\langle S_{\Delta^{(\varepsilon _0'-\varepsilon _0)}}M_{\varepsilon _0}(k)
S^*_{\Delta^{(\varepsilon _0'-\varepsilon _0)}}\xi_2,
\eta_2\rangle .}$$

To show independence of the choice of $\varepsilon_0$,
  we need consequently to prove that
$$\displaystyle {c_{t+12\varepsilon _0}\over{ c_t}}
\langle S_{\Delta^{(\varepsilon _0'-\varepsilon _0)}}M_{\varepsilon _0}(k)
S^*_{\Delta^{(\varepsilon _0'-\varepsilon _0)}}\xi_2,
\eta_2\rangle $$
is equal to
$\displaystyle {c_{t+12\varepsilon _0}\over{ c_t}}
\langle M_{\varepsilon _0'}(k)\xi_2,\eta_2\rangle $.

Now all the operators are in $L^1$. Moreover, the symbol of
$$S_{\Delta^{(\varepsilon _0'-\varepsilon _0)}}M_{\varepsilon _0}(k)
S^*_{\Delta^{(\varepsilon _0'-\varepsilon _0)}}$$
is $$\displaystyle {c_{t+12\varepsilon _0}\over{ c_{t+12\varepsilon _0'}}}
\varphi ^{\varepsilon _0'-\varepsilon _0}$$ times the symbol of
$M_{\varepsilon _0}(k)$.

  But the symbol of  $M_{\varepsilon _0}(k)$ is
$\varphi ^{\varepsilon _0}\ln \varphi $ divided by
$c_{t+12\varepsilon _0} / c_t$.

This shows independence of the choice of $\varepsilon _0$ (some care 
has to be taken
when choosing $\xi_1,\eta_1$, $\xi_2,\eta_2$ given $\xi,\eta$).
We always choose them in the init space of $S^*_{\Delta^{\varepsilon _0}}$,
respectively $S^*_{\Delta^{\varepsilon' _0}}$.
By the von Neumann theorem we will be able to choose a common
  intersection domain for these operators.

Consequently, to check that the derivative
of $\langle h\odot \varphi ^{\varepsilon }\xi,\eta\rangle _{H_t}$
at $\varepsilon =0$ is equal to the  operator $\Lambda(k)$ introduced 
in the statement of
Corollary 3.6, we only
have to check this for vectors $\xi,\eta$, that we assume to be of the form
$$\xi=S^*_{\Delta^{\varepsilon _0}}\xi_1,\qquad 
\eta=S^*_{\Delta^{\varepsilon _0}}\eta_1.$$

Then, modulo a multiplicative constant,
$\langle k\odot\varphi ^{\varepsilon }\xi,\eta\rangle _{H_t}$ becomes
$$\langle k\odot\varphi ^{\varepsilon +\varepsilon _0}
\xi_1,\eta_1\rangle _{H_{t+\varepsilon }}.$$

By a change of variables, the derivative at $0$ of
$\langle k\odot\varphi ^{\varepsilon }\xi,\eta\rangle _{H_t}$
becomes the derivative at $\varepsilon _0$ of the later expression:
$\langle k\odot\varphi ^{\varepsilon }\xi_1,\eta_1\rangle 
_{H_{t+\varepsilon }}$.
Up to a multiplicative constant,  this derivative exists and it is equal to
$\langle M(k)\xi_1,\eta_1\rangle $, which is by definition
$\langle \Lambda(k)\xi_1,\eta_1\rangle $.

Finally, observe that for any constant $C$,
$\langle (\Lambda(k)+C)\xi_1,\eta_1\rangle $ is equal to
$$B_{\varepsilon _0}(M_{\varepsilon _0}(k)+C'S_{\Delta^{\varepsilon _0}}
  kS^*_{\Delta^{\varepsilon _0}})
B^*_{\varepsilon _0}$$
for a constant $C'$ obtained from $C$ by multiplication by a normalization
factor depending on $t$ and $\varepsilon _0$.

Consequently, if $[k_{p,q}]_{p,q=1}^N$ is a positive matrix in $\ac_{s}$, 
then by using
the complete positivity result of Lemma~3.3, we infer that the matrix
$$-\left[M_{\varepsilon _0}\left(k_{p,q}\right)+C'\varphi 
^{\varepsilon _0}\odot
k_{p,q}\right]_{p,q=1}^N$$
represents a positive operator in
$M_p(\C)\otimes(\ac_{t+\varepsilon _0})_+$.

Since $\varphi ^{\varepsilon _0}\odot k_{p,q}$ is 
$S_{\Delta^{\varepsilon _0}}k_{p,q}
S^*_{\Delta^{\varepsilon _0}}$, we get that
$-\left[\Lambda (k_{p,q})\right]_{p,q=1}^N$ is a positive matrix of 
operators affiliated
to $M_p(\C)\otimes\ac_{t}$.\hfill\qed\vskip6pt

{\it Remark.}  If want to deal with less general operators 
(paying the price of
not including the identity operator in the domain of $\Lambda $), then
we can take operators of the form
$S_{\Delta^{\varepsilon _0}}kS_{(\Delta^{\varepsilon _0})^*}$ that
belong to $\ac_{s}$, $s<t-2$, $s-\varepsilon _0>1$, and then $\Lambda (k)$
will be in $L^1(\ac_{t})$, for such a kernel $k$, directly from Lemma~3.3.

\vfill\eject

\centerline{\S 4. {\ninebf Construction of an (unbounded) coboundary for}}
\centerline{\ninebf the Hochschild cocycle in Berezin's deformation}\vskip12pt
In this section we  analyze the Hochschild $2$-cocycle
$$\displaystyle \left.{\cal C}_{t}(k,l) =
 { \d\;\over{ \d s}}(k*_sl)\right|_{{s=t\atop s>t}}$$
that  arises in the Berezin deformation. We prove that the operator introduced
in the previous section (\S 3) may be used to construct an operator ${\cal L}$
(defined on a dense subalgebra of $\ac_t$), taking values in the 
algebra of unbounded
operators affiliated with $\ac_{t}$. ${\cal L}$ will be defined on a 
dense subalgebra
of $\ac_{t}$.

The equation satisfied by ${\cal L} $ is
$${\cal C}_t(A,B)={\cal L} _t(A*_tB)-A*_t{\cal L} _t(B) -{\cal L} _t(A)*_tB$$
and this will be fulfilled in the form sense (that is, by taking the 
scalar product
with some vectors $\xi,\eta$ in both sides).

The fact that ${\cal L} $ takes its  values in the unbounded
operators affiliated with $\ac_{t}$ presents
some inconvenience, but we recall that  in
the setting of type~$\twoone$ factors,  by von Neumann theory [MvN],
the algebra of unbounded (affiliated) operators is a well behaved algebra (with
respect to composition, sum and the adjoint operations).

In fact we will prove that ${\cal L} $ comes with two summands
$${\cal L} (k)=\Lambda (k)-{1\over2}\{T,k\}$$
where $-T$ is positive affiliated with $\ac_{t}$ and $-\Lambda $
  a completely positive (unbounded)
map. In the next section we prove that  $T$
is  $\Lambda (1)$.

For technical reasons (to have an algebra domain for  
${\cal L} $), we  require that $k\in \widehat {\ac}_{s}$, $s<t-2$, since
we know (by [Ra]) that the space of operators in $\ac_{t}$ 
represented by such kernels
is closed under taking the $*_t$ multiplication (the multiplication 
in $\ac_{t}$).

The operator $\Lambda $ will be (up to an additive multiple of the identity)
multiplication of the symbol by $\ln \varphi $. This operation is 
made more precise
in Corollary 3.6.

If $k$ is already of the form
$S_{\Delta^{\varepsilon_0}}kS_{\Delta^{\varepsilon_0}}^*$,
for some $k$ in $\ac_{s-\varepsilon _0}$, $s-\varepsilon _0>1$, $s<t-2$, then
$\Lambda (k)$ is an operator in $L^1(\ac_{t})$. In order to have the
identity $\Id$ in the domain) we allow $\Lambda $ to take its values in 
the operators
affiliated with $\ac_{t}$.

Consequently $\Lambda (1)$ is just the positive operator, affiliated
with $\ac_{t}$, which corresponds to the symbol
$\ln \varphi =\ln (\overline {\Delta (z)}\Delta (\xi)\zxi^{12})$
plus a suitable multiple of the identity.

To deduce the expression for ${\cal C}_t(k,l)$ one could argue 
formally as follows:
$$\eqalign{{\cal C}_t(k,l)(\overline {z},\xi)&=
{ c_t'\over{ c_t}}(k*_tl)(\overline {z},\xi)\cr
\noalign{\medskip}
&\qquad+\displaystyle c_t\int_{\H}k(\overline {z},\eta)l(\overline {z},\xi)
[\overline{z},\eta, \overline{\eta},\xi]^t\ln
[\overline{z},\eta, \overline{\eta},\xi]\,\d \nu_0(\eta).}
\eqno{(4.1)}$$
At this point to get  a $\Gamma$-invariant expression, we should decompose
$\ln [\overline{z},\eta, \overline{\eta},\xi]=
\ln \left[\left( (\overline z-\xi)(\overline {\eta}-\eta) \right) / \left( 
(\overline z-\eta)(\overline {\eta}-\xi) \right) \right]$
as a sum of $\Gamma$-invariant functions. The easier way to do that 
would be to write
$$\ln \varphi [\overline{z},\eta, \overline{\eta},\xi]={ 1\over{ 12}}
[\ln \varphi (\overline {z},\xi)+\ln \varphi (\overline {\eta},\eta)-
\ln \varphi (\overline {\eta},\xi)-\ln \varphi (\overline {z},\eta)].$$

If we use this expression back in $(4.1)$ we would get four terms which are
described as follows.

The term corresponding to $\ln \varphi (\overline {z},\xi)$ will come in front
of the integral and give
$${ 1\over{ 12}}\ln\varphi (\overline {z},\xi)(k*_tl)(\overline {z},\xi).$$

The term corresponding to $\ln \varphi (\overline{z},\eta)$ would multiply
$k(\overline{z},\eta)$ and would correspond formally to
${ 1\over{ 12}}
[(\ln \varphi) k]*_tl$.

The term corresponding to $\ln \varphi (\overline {\eta},\eta)$ would give the
following integral:
$$c_t\int_{\H}k(\overline {z},\eta)(\ln\varphi (\overline {\eta},\eta))
l(\overline {\eta},\xi)[\overline{z},\eta, \overline{\eta},\xi]^t
\,\d \nu_0(\eta).$$

This is formally ${ 1\over{ 12}}k*_tT^t_{\ln\varphi }*_tl $.
If $\ln \varphi $ were a bounded function and $T^t_{\ln\varphi }$ the Toeplitz
operator with this symbol, this  expression would make perfect sense.

Putting this together we would get
$$\eqalign{k*_t'l={\cal C}_t(k,l)&= { c_t'\over{ c_t}}k*_tl+
{ 1\over{ 12}}\ln\varphi (k*_tl)-\left[ \left( 
{ 1\over{ 12}}\ln \varphi \right) k\right] *_tl\cr
\noalign{\medskip}
&\qquad -\displaystyle k*_t\left[ \left( 
{ 1\over{ 12}}\ln \varphi \right) l\right] +
k*_tT^t_{( 1 / 12)\ln \varphi }*_tl.}$$


This would give that ${\cal C}_t(k,l)$ is implemented by the operator
$${\cal L} (k)=\left({ 1\over{ 12}}\ln \varphi -{ c_t'\over{ c_t}}\right)k
-{1\over2}\left\{T_{( 1 / 12)\ln\varphi },k\right\},$$
where by $\{a,b\}$ we denote the Jordan product $ab+ba$.

This means that ${\cal C}_t(k,l)$
is  implemented by the operator ${\cal L} (k)$, which  resembles the
  canonical form  of a generator of a  dynamical semigroup:
  a positive map ($-\ln \varphi $ is a positive kernel, when adding a 
constant) minus a Jordan
product.


To justify such a formula and the convergence of the integrals 
involved seems to
be a difficult task, so we will follow a different but more 
rigorous approach, which
consists in defining the operator $(-\ln \varphi )k$, as
  in the previous section, as a strong operator topology
  derivative.

To that end we introduce a family of completely positive maps
that canonically connect the fibers of the deformation. These maps 
arise from automorphic
forms, viewed (as in [GHJ]) as intertwining operators.

In the next lemma we give a precise meaning for the operator 
$T^t_{\ln\varphi }$,
which is the (unbounded) Toeplitz operator acting on $H_t$ with 
symbol $\ln\varphi $.\vskip6pt

{\sm Lemma} 4.1.\kern.3em
{\it We define $T=T^t_{\ln\varphi }$,  as a quadratic form, by
$$\langle T^t_{\ln\varphi }\xi,\xi\rangle _{H_t}=\displaystyle
\int_{\H}(\ln \varphi )\left| \xi\right| ^2\,\d\nu_t,$$ on the domain
$$\dc=\left\{\xi\in H_t\biggm| \displaystyle \int_{\H} (\ln \varphi )
\left| \xi\right| ^2\,\d \nu_t<\infty\right\}.$$
Clearly, $\dc$ is dense in $H_t$, as it contains
$\dc_0=\bigcup_{\varepsilon >0}\Range S_{\Delta^{\varepsilon 
}} $,
where $S_{\Delta^{\varepsilon }} $ is viewed as the operator of multiplication
by $\Delta^{\varepsilon }$ from $H_{t-\varepsilon }$ into $H_t$.

Moreover, $T^t_{\ln \varphi }$ is the restriction to
$H_t$ of the multiplication
operator by $\ln \varphi $ on $L^2(\H,\nu_t)$.
For $\xi,\eta$ in $\dc_0$ we have that
$$\langle T^t_{\ln \varphi }\xi,\eta\rangle _{H_t}=
\left. { \d\;\over{ \d\varepsilon }}\langle T^t_{ \varphi^\varepsilon }\xi,\eta
\rangle _{H_t}\right|_{\varepsilon =0}.$$
}\vskip6pt

{\it Proof.} All that stated above is obvious: the last statement is 
justified because,
if $S_{\Delta^{\varepsilon _0}}\colon H_{t-12\varepsilon_0}\to H_t$,
  then $S^*_{\Delta^{\varepsilon _0}}
T^t_{\ln \varphi }S_{\Delta^{\varepsilon _0}}$ is obviously equal to
$T^{t-\varepsilon_0 }_{\ln \varphi \varphi ^{\varepsilon _0}}$.\hfill\qed\vskip6pt

In the next lemma we explain the role of automorphic forms as 
comparison operators
between different algebras $\ac_{t}$ (they are
  a sort of tool for making a differentiable field
out of the algebras $\ac_{t}$).\vskip6pt

{\sm Definition} 4.2.\kern.3em
{\it
For $s\geq t$, let $\theta_{s,t}\colon \ac_{t}\to \ac_{s}$
be the completely positive map associating to $k$ in $\ac_{t}$
  the bounded operator in $\ac_{s}$
defined as
  $$\theta_{s,t}(k)=\left(S_{\Delta^{(s-t) / 12}}\right)k
\left(S_{\Delta^{(s-t) / 12}}\right)^*$$
Clearly the symbol of $\theta_{s,t}(k)$ is
$${ c_t\over{ c_s}}k(\overline {z},\xi)(\varphi(\overline 
{z},\xi))^{(s-t) / 12}.$$
Also we have $\theta _{s,t}(\theta _{t,v}(k))=\theta _{s,v}(k)$ for all
$s\geq t\geq v$.
}\vskip6pt

The following property is a trivial consequence of the definition of
$\theta _{s,t}$. It expresses the fact that $\theta _{s,t}$
has an almost multiplicative structure,
  as follows.\vskip6pt

{\sm Lemma} 4.3.\kern.3em
{\it For $s\geq t$ the following holds for all $k,l$ in
$\ac_{t}$:
$$\theta _{s,t}\left(k*_tT^t_{\varphi ^{(s-t) / 12}}*_tl\right)=
\theta _{s,t}(k)*_s\theta _{s,t}(l).$$
}\vskip6pt

{\it Proof.} This is obvious since $\theta _{s,t}(k)\theta _{s,t}(l)$ (with
product in $\ac_{s}$) is equal to $$S_{\Delta^{(s-t) / 12}}k
S^*_{\Delta^{(s-t) / 12}}S_{\Delta^{(s-t) / 12}}l
S^*_{\Delta^{(s-t) / 12}}.$$
But  an  obvious  formula
shows
$S^*_{\Delta^{(s-t) / 12}}S_{\Delta^{(s-t) / 12}}$
is equal to $T^t_{\varphi ^{(s-t) / 12}}$ (see, e.g., 
[Ra]).\hfill\qed\vskip6pt

We intend  next  to differentiate the above formula, in $s$,
while keeping $t$ fixed. In order to do this we will need to differentiate
$\theta _{s,t}(k)$. One problem that arises is the fact that
{\it a priori}
$\theta _{s,t}(k)$ belongs to $\ac_{s}$ rather than $\ac_{t}$. But if $k$
belongs to some $\ac_{t_0}$, with $t_0<t$, and $s$ is sufficiently 
close to $t$,
then $\theta _{s,t}(k)$ will be (up to a multiplicative constant) represented
by the symbol of $\theta _{s+t-t_0,t_0}(k)$. Since $s$ was small,
this defines (via $\Psi_{t, t-t_0+s}$)
a bounded operator in $\ac_{t}$.
Thus for such $k$ it makes sense to define
$\langle \theta _{s,t}(k)\xi,\eta\rangle _{H_t}$ for all vectors $\xi,\eta$
in $H_t$.

We differentiate this expression with respect to $s$. The existence of the 
derivative, in the strong
operator topology, was already
shown in the  previous section. We reformulate Corollary~3.6, in the 
new setting.\vskip6pt

{\sm Lemma} 4.4.\kern.3em
{\it Let $t>3$ and let $k$ belong to $\ac_{s_0}$, where
$1<s_0<t-2$. Then there exists a dense domain $\dc_0$
{\rm (}eventually depending on
$k${\rm )} that is affiliated with
$\ac_{t}$ such that the following expression:
$$\langle X_t(k)\xi,\eta\rangle _{H_t}=\left.{ \d\;\over{ \d s}}\right|_{s=t}
\langle \theta _{s,t}(k)\xi,\eta\rangle _{H_t}$$
defines a linear operator $X_t$ on $\dc$, which is affiliated  with $\ac_{t}$
{\rm (}and 
hence closable
{\rm [Ha]).}

Moreover, for a sufficiently large constant
$C$ {\rm (}depending on $s,t${\rm ),} 
$-X_t+C\cdot\Id $
becomes a completely positive map with values in the operators affiliated
to $\ac_{t}$.

Consider the {\rm (}non-unital\/{\rm )} subalgebra
$\widetilde {\ac}_{s_0}\subseteq 
\ac_{s_0}$,
which is also weakly dense, consisting of all operators $k$ in $\ac_{s_0}$
that are of the form
$S_{\Delta^{\varepsilon _0}}k_1S^*_{\Delta^{\varepsilon _0}}$ {\rm (}where 
$S_{\Delta^{\varepsilon _0}}$
maps $H_{s_0-12\varepsilon _0}$ into $H_s${\rm ),} where $k_1$
belongs to $\ac_{s_0-12\varepsilon _0}$ and $s_0-12\varepsilon _0$ is 
assumed bigger
than $1$.

  Then $X_t$ also maps $\widetilde {\ac}_{s}$ into
$L^1(\ac_{t})$. For such a $k$ the limit in the definition of $X_t$ 
is in the strong
operator topology on a dense, affiliated domain.
}\vskip6pt

Before going into the proof we make the following remark (which is 
not required
for the proof).\vskip6pt

{\it Remark.}  Since the kernel of the operator $\theta _{s,t}(k)$ (in
$\ac_{s}$) is equal to
$${ c_t\over{ c_s}}\cdot k(\overline {z},\xi)
[\varphi (\overline {z},\xi)]^{(s-t) / 12},$$
it follows that $X_t{(k)}$ is associated (in a sense that doesn't have to be
made precise for the proof) to
$$\left(-{ c_t\over{ c'_t}}+{ 1\over{ 12}}\ln\varphi \right)k,$$
which appeared in the formula in the introduction.\vskip6pt

{\it Proof of Lemma} 4.4. Because of the form of the symbol we  may use
  Corollary~3.6.\hfill\qed\vskip6pt

The main result of our paper shows that, by accepting an unbounded
coboundary, the Hochschild $2$-cocycle appearing in Berezin's deformation
is trivial, and the coboundary (which is automatically dissipative)
  has a form very similar to the canonical
expression of a generator of a quantum dynamical semigroup.

First we deduce a direct consequence out of the formula in Lemma~4.3.\vskip6pt

{\sm Proposition} 4.5.\kern.3em
{\it Fix a number $t>3$. Consider the algebra
$\widetilde {\ac}_t\subseteq \ac_{t}$ consisting of all $k\in \ac_{s}$
for some $s<t-2$ that are of the form
$S_{\Delta^{\varepsilon _0}}k_1S^*_{\Delta^{\varepsilon _0}}$, for 
some $\varepsilon _0$
{\rm (}such that $s-\varepsilon _0>1${\rm )} and
$k_1\in \ac_{s_0-\varepsilon }$.

Let $X_t$ be the operator defined in the previous lemma. Then
$$\left.{ \d\;\over{ \d s}}\right|_{s=t}\theta _{s,t}\left(
k*_tT^t_{\varphi ^{(s-t) / 12}}*_tl\right)=
X_t(k*_tl)+k*_tT^t_{(1 / 12)\ln \varphi }*_tl,
\eqno{(4.2)}$$
$$\left.{ \d\;\over{ \d s}}\right|_{s=t}[\theta _{s,t}
(k)*_s\theta _{s,t}(l)]=X_t(k)*_tl+{\cal C}_t(k,l)+k*_tX_t(l),
\eqno{(4.3)}$$
for all $k,l\in \widetilde{\ac}_t$. Consequently the two terms on the 
right-hand side
of $(4.2)$ and $(4.3)$ are equal, that is,
$$X_t(k*_tl)+k*_tT^t_{(1 / 12)\ln \varphi }*_tl=
X_t(k)*_tl+{\cal C}_t(k,l)+k*_tX_t(l).$$
}\vskip6pt

Before proceeding to the the proof of Proposition 4.5, we note that
  $\widetilde{\ac}_t$ is indeed
an algebra (see also the end of this section).
Assume that
$k,l$ are given, but that\break they correspond to two  different choices of
$s$, say $s,s'$, with $s'<s$. Because\break $\Psi_{s,s'}$ maps $\ac_{s'}$ 
into $\ac_{s}$,
we can assume $s=s'$. Then when $k=S_{\Delta^{\varepsilon 
_0}}k_1S^*_{\Delta^{\varepsilon _0}}$,\break
$l=S_{\Delta^{\varepsilon' _0}}'l_1S^*_{\Delta^{\varepsilon'_0}}$,
and say $\varepsilon _0'>\varepsilon _0$.
Then we replace the expression of $l$ as\break
$S_{\Delta^{\!\varepsilon _0}_{}}\left[S_{\Delta^{\!\varepsilon' 
_0-\varepsilon _0}}l_1\left(
S_{\Delta^{\!\varepsilon' _0-\varepsilon 
_0}}\right)^*\right]S^*_{\Delta^{\!\varepsilon _0}_{}}$,
and choose the new $l_1$ to be $S_{\Delta^{\!\varepsilon' 
_0-\varepsilon _0}}l_1\left(
S_{\Delta^{\!\varepsilon' _0-\varepsilon _0}}\right)^*\!.$\vskip6pt

{\it Proof of Proposition} 4.5. We will give separate proofs for each 
of the equalities
(4.2), (4.3). Of course these are the product formula for 
derivatives, but the complicated nature
of the operator functions obliges us to work on the nonunital algebra
$\widetilde{\ac}_t$. This might be
just a technical condition that perhaps could be dropped.\vskip6pt

{\it Proof of equality} (4.2).
$${ \d\;\over{ \d s}}\theta _{s,t}\left(k*_tT^t_{\varphi ^{(s-t) / 12
}}*_tl\right)=
X_t(k*_tl)+k*_tT^t_{(1 / 12)\ln \varphi ^*}*_t l.$$

We start with the left-hand side:
Denote $P_s=k*_tT^t_{\varphi ^{(s-t) / 12}}*_tl$, for fixed 
$h,l\in \widetilde{\ac}_t$.

We have to evaluate (against
$\langle\,\cdot \, , \xi\rangle \eta$, where $\xi,\eta$
belong to a suitable dense domain $\dc$ affiliated with $\ac_{t}$)
the expression:
$${ \theta _{s,t}(P_s)-P_t\over{ s-t}}=\theta _{s,t}\left({ 
P_s-P_t\over{ s-t}}\right)+
{ \theta _{s,t}(P_t)-P_t\over{ s-t}}.$$

The second term converges when $s\searrow t$, since $P_t=k*_tl$, to
$\left.{ \d\;\over{ \d s}}\right|_{s=t}\theta _{s,t}(P_t)$,
which is by definition $X_t(k*_tl)$. Here we rely on the fact that 
$\widetilde {\ac}_t$
is an algebra, so that $k*_tl$ belongs to the domain of $X_t$.

For the limit of the expression
$\theta _{s,t}\left( ( P_s-P ) / ( s-t ) \right)$,
we note that $$\displaystyle { P_s-P_t\over{ s-t}}=k*_tT^l_{h_{s,t}}*_tl,$$
  with
$$h_{s,t}={ \varphi ^{(s-t) / 12}-\Id \over{ s-t}}.$$

Assume now that $k=S_{\Delta^{\varepsilon 
_0}}k_1S^*_{\Delta^{\varepsilon _0}}$,
$l=S_{\Delta^{\varepsilon _0}}l_1S^*_{\Delta^{\varepsilon _0}}$.

Then $\theta _{s,t}(k*_tT^t_{h_{s,t}}*_tl)$ is equal to
$$S_{\Delta^{(s-t) / 12}}\left(\left( S_{\Delta^{\varepsilon _0}}
k_1S_{\Delta^{\varepsilon_0}}^*\right)
*_tT^t_{h_{s,t}}*_t\left(S_{\Delta^{\varepsilon _0}}l_1
S^*_{\Delta^{\varepsilon _0}}\right)\right)
S^*_{\Delta^{(s-t) / 12}}.$$

This is easily seen to be equal to
$$\eqalign{&S_{\Delta^{( (s-t) / 12 )+\varepsilon _0}}
\left[k_1*_{t-\varepsilon _0}
T^t_{\varphi ^{\varepsilon _0}h_{s,t}}*_{t-\varepsilon }l_1\right]
S^*_{\Delta^{( (s-t) / 12 )+\varepsilon _0}}\cr
\noalign{\medskip}
&\qquad =\theta _{s,t-\varepsilon _0}\left(k_1*_{t-\varepsilon _0}
T^t_{\varphi ^{\varepsilon _0}h_{s,t}}
*_{t-\varepsilon _0}l_1\right).}$$

Denote $\widetilde {P}_s=k_1*_{t-\varepsilon _0}
T^t_{\varphi ^{\varepsilon _0}h_{s,t}}
*_{t-\varepsilon _0}l_1$.

As $s$ decreases to $t$, we have (as $\varphi ^{\varepsilon _0}\ln \varphi $
is bounded) that $\widetilde {P}_s$ converges in the uniform operator 
topology, to
$k_1*_{t-\varepsilon _0}T^t_{( 1 / 12)\varphi^{ \varepsilon 
_0}\ln\varphi }
*_{t-\varepsilon _0}l_1$. This is because
$\varphi ^{\varepsilon _0}\left( ( \varphi ^{(s-t) / 12
}-\Id ) / ( s-t ) \right) $
converges uniformly to ${ 1\over{ 12}}\varphi 
^{\varepsilon _0}\ln\varphi $,
since $\varphi $ is a bounded function.

Also, if $s$ is sufficiently close to $t$, $\theta _{s,t-\varepsilon 
}(k_1)$ defines
for every $k_1\in \ac_{t-\varepsilon_0}$ a bounded operator on 
$\ac_{t}$. Indeed,
$\theta _{s,t-\varepsilon _0}(k_1)$ has symbol (up to a 
multiplicative constant) equal to
$\varphi ^{( s-t+\varepsilon _0 ) / 12}k_1$. This is well defined as
the kernel of an operator in $\ac_{t}$, since $k_1\in
\ac_{t_0-\varepsilon _0}$.


Thus $\theta _{s,t_0-\varepsilon }$ can be thought of as a completely 
positive map
from $\ac_{t_0-\varepsilon }$ into $\ac_{t}$. Moreover,
$\theta _{s,t-\varepsilon _0}(1)$, which is
$S_{\Delta^{( s-t_0+\varepsilon ) / 12}}
S^*_{\Delta^{( s-t_0+\varepsilon ) / 12}}$,
is less than a  constant $C$ (not depending on $s$) times the identity.

Hence the linear maps $\theta _{s,t_0-\varepsilon }$, acting from 
$\ac_{t_0-\varepsilon }$
into $\ac_{t}$, are uniformly bounded.
Consequently, when evaluating
$$\left| \left\langle \left( 
\theta _{s,t_0-\varepsilon }(\widetilde {P}_s)-
\theta _{t,t_0-\varepsilon }(\widetilde {P}_t)
\right)\xi,\eta\right\rangle \right| ,$$
we can majorize by
$$\left| \langle \theta _{s,t_0-\varepsilon }(\widetilde {P}_s-\widetilde 
{P}_t)\xi,\eta\rangle \right| +
\left| \langle \left(\theta _{t,t_0-\varepsilon }-\theta _{s,t_0-\varepsilon }\right)
(\widetilde {P}_t)\xi,\eta\rangle \right| .$$

The first term goes to zero by uniform continuity of $\theta 
_{s,t_0+\varepsilon }$
with respect to $s$ (and since $\|\widetilde {P}_s-\widetilde {P}_t\|\to 0$).
The second goes to zero because of the pointwise strong-operator-topology 
continuity of the map
$s\to \theta _{s,t_0-\varepsilon }$.

Thus $\theta _{s,t}\left( ( P_s-P_t ) / ( s-t ) \right)$ 
converges to
$\theta _{t,t_0-\varepsilon }(\widetilde {P}_t)$,
which was
$$S_{\Delta^{\varepsilon _0}}\left(k_1*_{t-\varepsilon _0}
T^t_{( 1 / 12)\varphi ^{\varepsilon _0}\ln \varphi }*l_1\right)
S^*_{\Delta^{\varepsilon _0}},$$
which is equal to $k*_tT^t_{\ln \varphi }*_tl$.\vskip6pt

{\it Proof of equality} (4.3).
$$\left.{ \d\;\over{ \d s}}\right|_{s=t}\theta _{s,t}(k)*_s\theta _{s,t}(l)=
X_tk*_tl+{\cal C}_t(k,l)+k*_tX_t.$$

We  verify this equality by evaluating  
it on $\langle \,\cdot \, ,\xi\rangle \eta$,
$\xi,\eta\in \dc$, where $\dc$ is a dense domain affiliated to $\ac_t$.

We write the expression
$${ \theta _{s,t}(k)*_s\theta _{s,t}(l)-k*_tl\over{ s-t}}$$
as
$${ \theta _{s,t}(k)*_s\theta _{s,t}(l)-
\theta _{s,t}(k)*_t\theta _{s,t}(l)\over{ s-t}}+
{ \theta _{s,t}(k)*_t\theta _{s,t}(l)-k*_tl\over{ s-t}}.$$

We will analyze first the first summand and prove that
$${ \theta _{s,t}(k)*_s\theta _{s,t}(l)-
\theta _{s,t}(k)*_t\theta _{s,t}(l)\over{ s-t}}
\eqno{(4.4)}$$
converges to ${\cal C}_t(k*_tl)$.

We use the symbols of $k,l$, and then
find that
the symbols of $\theta _{s,t}(k)$,
$\theta _{s,t}(l)$ are $\varphi ^{( s-t ) / 12}k$,
$\varphi ^{( s-t ) / 12}l$, up to multiplicative constants that 
we ignore here
(because the argument has a qualitative nature).

Then the symbol of the expression in (4.4)
is
$$\int_{\H}\left(k\varphi ^{( s-t ) / 12}\right)(\overline {z},\eta)
\left(l\varphi ^{( s-t ) / 12}\right)(\overline {\eta},\xi)
{ [\overline{z},\eta,\overline{\eta},\xi]^s-
[\overline{z},\eta,\overline{\eta},\xi]^t\over{ s-t}}\,\d \nu_0(\eta) .
\eqno{(4.5)}$$

By  the mean-value theorem, with $\alpha _s(v)=sv+(1-v)t$, the 
expression becomes
$$\int^{1}_0\int_{\H}\left(k\varphi ^{( s-t ) / 12}
\right)(\overline {z},\eta)
\left(l\varphi ^{( s-t ) / 12}\right)(\overline {\eta},\xi)
[\overline{z},\eta,\overline{\eta},\xi]^{\alpha_s (v)}
\ln [\overline{z},\eta,\overline{\eta},\xi]
\,\d \nu_0(\eta) \,\d v.$$

Similarly (by ignoring the numerical factors due to
the constants $c_s$) we have that
$${\cal C}_t(\theta _{s,t}(k),\theta _{s,t}(l))$$
contains the integral
$$\int_{\H}\left(k\varphi ^{( s-t ) / 12}\right)(\overline {z},\eta)
\left(l\varphi ^{( s-t ) / 12}\right)(\overline {\eta},\xi)
[\overline{z},\eta,\overline{\eta},\xi]^{t}
\ln [\overline{z},\eta,\overline{\eta},\xi]
\,\d \nu_0(\eta) .
\eqno{(4.6)}$$

Taking the difference, we obtain the following integral:
$$\eqalign{\displaystyle
&\int_{\H}\left(k\varphi ^{( s-t ) / 12}\right)(\overline {z},\eta)
\left(l\varphi ^{( s-t ) / 12}\right)(\overline {\eta},\xi) \cr
\noalign{\medskip}
\displaystyle &\qquad
\cdot \left({ [\overline{z},\eta,\overline{\eta},\xi]^s-
[\overline{z},\eta,\overline{\eta},\xi]^t\over{ s-t}}-
[\overline{z},\eta,\overline{\eta},\xi]^{t}
\ln [\overline{z},\eta,\overline{\eta},\xi]\right)
\,\d \nu_0(\eta) .}
\eqno{(4.7)}$$

By Taylor expansion, this is $(s-t)$ times a term involving  the integral
$$\int_{\H}\left(k\varphi ^{( s-t ) / 12}\right)(\overline {z},\eta)
\left(l\varphi ^{( s-t ) / 12}\right)(\overline {\eta},\xi)
[\overline{z},\eta,\overline{\eta},\xi]^{s'}
\left(\ln [\overline{z},\eta,\overline{\eta},\xi]\right)^2
\,\d \nu_0(\eta),
\eqno{(4.8)}$$
where $s'$ is the interval determined by $s$ and $t$.


We have to prove that the integral of the absolute value of the
integrand in the above
integral is bounded by a constant independent of the choices
of $s'$ (and $s$).

We write $\left| [\overline{z},\eta, \overline{\eta},\xi]\right| =
d(\overline{z} , \eta)  d(\overline{\eta},\xi) /
d(\overline{z},\xi)$. Then
$$\left| \ln\left| [\overline{z},\eta, \overline{\eta},\xi]\right| \right| \leq
\left| \ln\left|  d(\overline{z},\eta)\right| \right| +\left| \ln\left|   d(\overline{\eta},\xi)\right| \right| +
      \left| \ln\left|   d(\overline{z},\xi)\right| \right| .$$
Also we note that the logarithm in
$\ln [\overline{z},\eta, \overline{\eta},\xi]$ has
bounded imaginary part, as the branches in
$$\ln \zueta,\quad \ln \etxi,\quad
\ln \zxi$$ have imaginary part in the fixed segment $[0,2\pi]$.

Thus the term that we have to evaluate will involve terms of the form
$$\int_{\H}\left| k'(\overline{z},\eta)\right|  \left| l' (\overline{\eta},\xi)\right| 
\left|  d(\overline{z},\eta)\right| ^{s'} \left|  d(\overline{\eta},\xi)\right| ^s \,\d \nu_0(\eta),$$
where $k'{(\overline{z},\eta)}$ could be
$k\varphi^ {( s-t ) / 12} (\overline{z},\eta)$,
eventually multiplied
by a power ($1$ or $2$) of $\ln d (\overline{z},\eta).$
A similar assumption holds for $l$.

By  the Cauchy--Schwarz inequality, this expression is bounded by products of:
$$\left(
\int_{\H}\left| k'(\overline{z},\eta)\right| ^2
\left|  d(\overline{z},\eta)\right| ^{2s} \, d\nu_0(\eta)\right)^{1/2}
\left(\int_{\H}\left| l'(\overline{z},\xi)\right| ^2
\left|  d(\overline{\eta },\xi)\right| ^{2s} \,\d \nu_0(\eta)\right)^{1/2}.\eqno{(4.9)}$$

But such expressions are finite, because we know that $k,l$ are in
$\ac_{t_0}$ for some fixed $t_0<t$, and hence in $L^2(\ac_{t_0})$, and
consequently the integral
$$\int_{\H}\left| k(\overline{z},\eta)\right| ^2\left|  d(\overline{z},\eta)\right| ^{2t_0}
\,\d \nu_0(\eta)\eqno{(4.10)}$$
is finite.

Moreover, $\varphi(\overline{z},\xi)=
\overline{\Delta(z)}\Delta(\xi)\zxi^{12}$.
This  term is bounded in absolute value by $\left|  d(\overline{z},\xi)\right| ^{-12}$,
so that $\left| \varphi(\overline{z},\xi)\right| ^{( s-t ) / 12}$ is bounded by
$\left|  d(\overline{z},\xi)\right| ^{s-t}$.
Also $\left| \ln \left|  d(\overline{z},\eta)\right| \right| \left|  d(\overline{z},\eta)\right| ^ {\varepsilon}$
is bounded for any choice of $\varepsilon$.

Thus, by choosing $\varepsilon$ small enough, the finiteness of the
integral in (4.10) implies the finiteness of the integral in (4.9).

Thus the integral in (4.7) tends to zero.
This  completes the proof that
$$\displaystyle {\theta_{s,t}(k) x_s \theta_{s,t}(l)-\theta_{s,t}(k) x_t
\theta_{s,t}(l)\over {s-t}}$$
converges to ${\cal C}_t(k,l)$.

The remaining term to be analyzed is
$$\left.\displaystyle {\d\;\over{\d s }}\right|_{s=t}
\displaystyle {\theta_{s,t}(k) *_t \theta_{s,t}(l)-k *_tl
\over {s-t}}.$$
We have to show that the limit is
$(X_t(k)) *_tl+k *_t (X_t(l))$ (evaluated on vectors
$\xi,\eta$ in a dense domain affiliated to ${\ac_t}$).

Fix the vector $\xi,\eta$. Then we have to analyze the following sum:
$$\left\langle \displaystyle{
\theta_{s,t} (l)-l\over {s-t}} \xi,
\theta_{s,t}(k^{*})\eta\right\rangle +
\left\langle \displaystyle{
\theta_{s,t} (k)-k\over {s-t}} l\xi,\eta\right\rangle .$$
The second term obviously converges to
$\langle [X_t(k)*_t]l \xi,\eta\rangle $,
and
the first term is also convergent to
$\langle X_t(l) \xi k^*,\eta\rangle $, because
$\left( (
\theta_{s,t} (l)-l ) / (s-t) \right) \xi$, for $\xi$ in a dense
domain $\dc$, converges in norm to $X_t(l)\xi$. Indeed in Corollary~3.6
  we proved that
$(
\theta_{s,t} (l)-l ) / (s-t)$ converges strongly to
$X_t$, on a dense domain topology, because the convergence (for
$l=S_{\Delta^{\varepsilon_0}} l_1(S^*_{\Delta^{\varepsilon_0}}) $,
  $l_1\in \ac_{t_0-\varepsilon_0}$)
comes by proving that the partial fractions
$-(
\theta_{s,t} (l)-l ) / (s-t)$
increase (modulo $(s-t)$ times a constant) to $-X_tl$.

This completes the proof.\hfill\qed\vskip6pt


We are now able to formulate our main result. We recall
first the context of this result. The algebras $\ac_t$ are
the von Neumann algebras (type $\twoone$ factors) associated
with Berezin's deformation of $\H/\PSLtwoZ$.
These algebras can be realized as subalgebras of $B(H_t)$
where $H_t$ is the Hilbert space $H^2(\H ,\Im z^{t-2}
\,\d\overline{z}\,\d z)$.

As such, every operator $A$ in $\ac_t$ (or $B(H_t)$)
is given by a reproducing kernel: $k_A$, which is a bivariable function
on $\H$, analytic in the second variable, antianalytic in the first and
$\PSLtwoZ$-invariant. The symbols are normalized so that
the symbol of the identity is the constant function $1$.

By using these symbols (that represent the deformation)
we can define the $*_t$ product of two symbols $k,l$ by
letting $k*_tl$  be the product symbol in the algebra
$\ac_t$.

The Hochschild $2$-cocycle associated with the deformation
is defined by
$$\left.{\cal C}_t(k,l)=
\displaystyle{\d\;\over {\d s}} (k*_sl)\right|_{s=t}.$$
The Hochschild cocycle condition is obtained by
differentiation of the associativity identity.

The cocycle ${\cal C}_t$ is well defined on a weakly dense,
unital subalgebra $\widehat {\ac}_t$ of $\ac_t$.
A sufficient condition that an element in $\ac_t$, represented by
a symbol $k$, belongs to $\widehat\ac_t$, is that the quantity
$\| k \|^{\,\widehat {}}_t$, defined as the maximum
of $$\displaystyle\sup_{z\in H}
\int_{\H} \left| k(\overline{z},\eta)\right|  \left| d(\overline{z},\eta)\right| ^t \,\d\nu_0(\eta)$$
and
$$\displaystyle\sup_{\eta \in H}
\int_{\H} \left| k(\overline{z},\eta)\right|  \left| d(\overline{z},\eta)\right| ^t \,\d\nu_0(z),$$ be
finite.

The algebra  $\widehat{\ac}_t$ is
the analogue of the Jolissaint algebra [Jo] for discrete groups.

We proved in Section 2 that the applications $\Psi_{s,t}$ which
map the operator $A\ \hbox{in}\ \ac_t$ into the corresponding
operator in $\ac_s$ having the same symbol are completely positive.

This property proves that ${\cal C}_t$ is completely negative, that is,
for all $l_1,l_2,\ldots ,l_N$ in $\ac_t$,
for all $k_1,k_2,\ldots ,k_N$ in $\widehat\ac_t$, we
have that
$$\sum l_i^*c(k_i^*, k_j)l_j\leq0.$$

This property could be used to construct, as in [Sau],
the cotangent bundle. In fact, here
${\cal C}_t$, or rather $-{\cal C}_t$, plays the role of $\nabla L$, where
$L$ should be a generator of a quantum dynamical semigroup
${\Phi_t}$ (thus $L=
\left.\displaystyle {\d\;\over{\d s }} \Phi_s\right|_{s=0}$)
and we  have
$\nabla L(a,b)=L(a,b)-aL(b)-L(a)b$.

It is well known that $\nabla L$ is completely negative
[Li]. In our case, the role of the quantum dynamical semigroup
is played by the completely positive maps
$\Psi_{s,t}$ that have the property $\Psi_{s,t}\Psi_{t,v}=
\Psi_{s,v}$, $s\geq t\geq v$.
The generator $L$ doesn't make sense here, since $\Psi_{s,t}$
takes its values in different algebras, depending on $s$.

  Instead we use
the derivative of the multiplication operation,
which formally is $$\left.\displaystyle {\d\;\over{\d s }}
\Psi_{s,t}^{-1}(\Psi_{s,t}(k)*_s \Psi_{s,t}(l))\right|_{s=t},$$
as a substitute for $\nabla L$.

All the above is valid for
the general
case of
Berezin's deformation of $H/\Gamma$,
where $\Gamma$ is any discrete subgroup of
$\PSLtwoR$ of finite covolume.

When specializing to $\Gamma=\PSLtwoZ$, we  construct also
the diffusive operator
  ${\cal L}$, which plays the role of
the generator of a dynamical semigroup.

In the next theorem we formulate our main result. We construct 
explicitly an operator
${\cal L}$ such that $${\cal L} (ab)- {\cal L} (a)b-a {\cal L} 
(b)={\cal C}_t(a,b).$$
We will show that ${\cal L}$ is well defined on a weakly dense
(non-unital) subalgebra $\dc_t^0$ and the above relation holds for $a,b\in
\dc_t^0$ (which is obtained by considering suitable subalgebras
of $\widehat{\ac_s}, s<t-2)$.
Moreover, ${\cal L}$ has an expression that is very similar
to the Lindblad [Li, CE, GSK, HP] form of the generator $L$
of a uniformly continuous semigroup. Recall that this
expression is, in the uniformly continuous case,
$$L{(x)}=\Phi(x)-{1\over2} \{\Phi(1), x\}+i[H,x],$$
  where
$\Phi$ is completely positive and $H$ is selfadjoint.

In our case (which is certainly not [Dav]
corresponding to the uniformly continuous case) the generator
${\cal L}(x)$ is defined rather as an unbounded operator
(which is the approach taken in  [Dav, CF, Ho,
GS, MS]).

We prove that there exists a weakly dense, {\it unital} algebra
$\dc_t$ containing $\dc_t^0$, and a linear map $\Lambda$ from
$\dc_t$ into the operators affiliated with $\ac_t$, and
a positive operator that is also affiliated to $\ac_t$,
such that
$${\cal L}(x)=\Lambda(x)-{1\over2}\{T,x\}.$$
Also $\Lambda$ maps $\dc_t^0$ into $L^1(\ac_t)$

Moreover, $\Lambda$ has properties that are very similar
to a completely positive map. We prove that there exists
an increasing filtration $(\bc_{rt})_{1<r<t-2}$ of $\dc_t$, consisting of
weakly dense subalgebras, such that, for a constant
$C_{rt}^0$ depending on $r$, $-[\Lambda +C_{rt}^0\cdot\Id]$ is a completely
positive map on $\bc_{rt}$.

This means that when  restricted to $\bc_{rt}$, ${\cal L}$
has the form
${\cal L}(x)=\Lambda'(x)-{1\over2}\{T',x\}$,
where $-\Lambda'=-[\Lambda+C_r^0\cdot\Id ]$ is a completely positive map
and $T=T+C^0_{rt}\cdot\Id $.\vskip6pt

{\sm Theorem} 4.6.\kern.3em
{\it
Let $\ac_t$, $t>1$, with product operation $*_t$ be the
  von Neumann algebra {\rm (}a type $\twoone$ factor\/{\rm )}
associated with 
Berezin's
deformation of $\H/\PSLtwoZ$.

Let ${\cal C}_t$ be the Hochschild $2$-cocycle associated with the deformation
$$\left.{\cal C}_t(k,l)=\displaystyle {\d\;\over{\d s}} k*_sl\right|_{s=t},$$
which is defined on the weakly dense subalgebra $\widehat\ac_t$.

Then there exists a weakly dense {\rm (}non-unital\/{\rm )} subalgebra
$\dc_t^0$ in $\widehat\ac_t\subseteq\ac_t$ and ${\cal L}_t$, a linear 
operator on
$\dc_t^0$, with values in the algebra of operators affiliated with
$\ac_t$, such that
$${\cal C}_t(k,l)={\cal L}_t(kl)-k{\cal L}_t(l)-
{\cal L}_t(k)l,\qquad k,l\in \dc_t^0.$$
Note that $-{\cal L}_t$ is automatically completely dissipative.

Moreover, ${\cal L}_t$ has the following expression.
There exists a weakly dense, {\rm unital}
subalgebra $\dc_t$, such that $\dc_t^0\subseteq \dc_t\subseteq
\widehat{\ac}_t$, there exists $\Lambda_t$ defined on $\dc_t$ with
values in the operators affiliated to $\ac_t$, and
there exists $T$, a positive unbounded operator affiliated with $\ac_t$ such
that
$${\cal L}_t(k)=\Lambda_t(k)-{1\over2}\{T,k\},\qquad k\in \dc_t^0.$$

Moreover, $\Lambda_t$ has the following complete positivity properties:

{\rm 1)} $\Lambda_t$ maps $\dc_t^0$ into $L^1(\ac_t)$;

{\rm 2)} There exists an increasing filtration of weakly dense, unital
subalgebras $(\bc_{s,t})_{1<s<t-2}$ of $\dc_t$, with
$\bigcup_{s} \bc_{s,t}=\dc_t$, and there exist constants $C_{s,t}$ such that
$-[\lambda_t+C_{s,t}\cdot\Id ]$ is completely positive on $\bc_{s,t}$.}\vskip6pt

{\it Remark.} At the level of symbols the operator
$\Lambda_t$ has a very easy expression, namely $\Lambda_t{(k)}$
is the pointwise multiplication (the analogue of Schur multiplication)
of $k$ with the $\Gamma$-equivariant symbol
$$\ln(\overline{\Delta{(z)}}\Delta{(\xi)}\zxi^{12}).$$
We identify, as in Section 1, $L^2(\ac_t)$ with a Hilbert space
of $\Gamma$-bivariable functions analytic in the first variable
and antianalytic in the second.

Then $\Lambda$ corresponds to the (unbounded) analytic Toeplitz operator
with symbol $\ln\left(\overline{\Delta{(z)}}\Delta{(\xi)}\zxi^{12}\right)$.\vskip6pt

{\it Proof of Theorem} 4.6. This was almost proved in
  Lemma~3.3 and Proposition~4.5, but we have to
identify the ingredients. Here the algebra $\dc_t^0$ is the
union (with respect to $s,\varepsilon_0$)
$$\bigcup_{1<s-\varepsilon_0<s<t-2}
S_{\Delta^{\varepsilon_0}}\ac_{s-\varepsilon_0}S^*_{\Delta^{\varepsilon_0}}.$$
It is obvious that $\dc_t^0$ is an algebra (under the product on
$\ac_t$---the algebra $\dc_t$ is the union
$\bigcup_{1<s<t-2}\ac_t$, viewed as an algebra
of $\ac_t$). The algebra $\bc_{s,t}$ is the union
(with respect to $\varepsilon$)
of
$\bigcup_{1<s-\varepsilon_0}
S_{\Delta^{\varepsilon_0}}\ac_{s-\varepsilon_0}S^*_{\Delta^{\varepsilon_0}}.$

The operator $T$ is the Toeplitz operator with symbol
${1\over{12}}\ln \varphi$, while $\Lambda_t$ is
$X_t$, where $X_t$ was defined in Lemmas 3.3 and 4.4.
In Proposition 4.5 we also proved that
$$
{\cal C}_t(a,b)=X_t(a *_t b)-X_t(a) *_t b-a *_t X_t b+a *_t T^t
_{(\ln \varphi )/12}{*_t} b, \hbox{\qquad for\ all\ }
a,b\in\dc_t^0.$$

Clearly the term $a *_t T^t
_{(\ln \varphi )/12}{*_t} b$ is a cohomologically trivial term,
  and hence ${\cal C}_t(a,b)$
is implemented by ${\cal L}_t(a)=X_t(a)-
{1\over2}\{a,T^t_{(\ln \varphi )/12}\}$.
Hence ${\cal C}_t$ is implemented by
${\cal L}_t=\Lambda_t(a)-
{1\over2}\{a,T^t_{(\ln \varphi )/12}\}$.
All the other properties for $\Lambda_t$ were proven in
Section 3.

One also needs to show that the vector spaces $\dc_t=
\bigcup_{s<t-2} \widehat{\ac}_s$
and
$$\dc_t^0=\bigcup_{1<s-\varepsilon_0<s<t-2}
S_{\Delta^{\varepsilon_0}}
\widehat{\ac}_{s-12\varepsilon_0}
S^*_{\Delta^{\varepsilon_0}}$$ are indeed algebras (in $\ac_t$).
$\dc_t$ is obviously an algebra, since we proved (in [Ra])
that $\widehat{\ac}_s$ is closed under $*_v$ for all
$v\geq s$. Of course, if we take the product of different
$\widehat{\ac}_{s_1}$ and $\widehat{\ac}_{s_2}$
we may embed them in
$\widehat{\ac}_{\max(s_1,s_2)}$.

To prove that $\dc_t^0$ is an algebra (in $\ac_t$) we will
need to show first that we are reduced to proving that
$ S_{\Delta^{\varepsilon_0}} \widehat{\ac}_{s-12
\varepsilon_0} S^*_{\Delta^{\varepsilon_0}}$, for
fixed $s$ and $\varepsilon_0$, is closed under the product $*_t$
in $\ac_t$.

Indeed if we do a product for different $s$, we may
simply take the maximum of  $s',s$. If we do a product corresponding to
different $\varepsilon_0$'s, say $\varepsilon_0$
and $\varepsilon_1$, then we choose $\varepsilon_1$ to be
the largest.

Then observe that for $k\in \ac_{ s-12\varepsilon_1}$
$$S_{\Delta^{\varepsilon_1}} k S^*_{\Delta^{\varepsilon_1}}=
S_{\Delta^{\varepsilon_0}}(S_{\Delta^{\varepsilon_1-\varepsilon_0}} k
S^*_{\Delta^{\varepsilon_1-\varepsilon_0}})S^*_{\Delta^{\varepsilon_0}}.$$

Now $S_{\Delta^{\varepsilon_1-\varepsilon_0}} k
S^*_{\Delta^{\varepsilon_1-\varepsilon_0}}$
has symbol equal to (modulo a multiplicative constant)
$\varphi^{\varepsilon_1-\varepsilon_0}k$. Since
$\left| \varphi\right| \leq d^{-12}$, it follows that
$\left| \varphi\right| ^\varepsilon\leq d^{-12\varepsilon}$ and hence that
$\varphi^{\varepsilon_1-\varepsilon_0}k$ belongs to
$\widehat {\ac}_{s-12\varepsilon_1+12(\varepsilon_1-\varepsilon_0)}$,
which is $\widehat{\ac}_{s-12\varepsilon_0}$.

Thus $S_{\Delta^{\varepsilon_1}}\widehat{\ac}_{s-12\varepsilon_1}
S^*_{\Delta^{\varepsilon_1}}$ is contained in
$S_{\Delta^{\varepsilon_0}}\widehat{\ac}_{s-12\varepsilon_0}S^*_{\Delta^{\varepsilon_0}}.$

Now we are reduced to showing that the product of two elements
$S_{\Delta^{\varepsilon_0}}k_1 S^*_{\Delta^{\varepsilon_0}}$
and
$S_{\Delta^{\varepsilon_0}} l_1 S^*_{\Delta^{\varepsilon_0}}$,
$k_1,l_1 \in \widehat{\ac}_{s-12\varepsilon_0} $,
is again an element in
$S_{\Delta^{\varepsilon_0}}\widehat{\ac}_{s-12\varepsilon_0}
S^*_{\Delta^{\varepsilon_0}}.$

But $$(S_{\Delta^{\varepsilon_0}} k_1
S^*_{\Delta^{\varepsilon_0}} )*_t
(S_{\Delta^{\varepsilon_0}} l_1
S^*_{\Delta^{\varepsilon_0}})$$
coincides with
$$S_{\Delta^{\varepsilon_0}}[k_1*_{t-12\varepsilon_0}
S^*_{\Delta^{\varepsilon_0}}
S_{\Delta^{\varepsilon_0}}
*_{t-12\varepsilon_0}l_1]
S_{\Delta^{\varepsilon_0}}.$$
Because $\widehat{\ac}_{s-12\varepsilon_0}$ is closed under
the product $*_{t-12\varepsilon_0}$, it is sufficient to show that
$$T_{\varphi^{\varepsilon_0}}^{t-\varepsilon_0}=
S^*_{\Delta^{\varepsilon_0}}
S_{\Delta^{\varepsilon_0}}$$ belongs to $\widehat{\ac}_{s-12\varepsilon_0}$.
But this is a general fact contained in the following lemma.\vskip6pt

{\sm Lemma} 4.7.\kern.3em
{\it Assume $f$ is a bounded, measurable,
$\Gamma$-equivariant function on $\H$.
Let $T_f^t$ be the Toeplitz operator on $H_t$ with symbol $f$.
Then $T_f^t$ belongs to $\widehat{\ac}_t$.
Moreover,
$\|T_f^t \|^{\,\widehat {}}_t\leq C\|f\|_{\infty}$,
where $C$ is a constant depending on $t$.}\vskip6pt

{\it Proof.} Note the symbol of $T_f^t$ is given by the formula
[Ra]
$$
s_f(\overline{z},\xi)=\int_{\H}f(a)
[z,a,\overline{a},\xi]^t \,\d \nu_0(a),\qquad z,\xi\in \H .
$$
We have to check that the quantity
$$\sup_{z\in \H}\int_{\H} \left| S_f(z,\xi)\right| 
\left|  d(z,\xi)\right| ^t\,\d \nu_0(\xi)\leq \|f\|_\infty
$$
(and a similar one) is finite.

But the above integral is bounded by
$$\eqalign
{\displaystyle &\mathop{\int\kern-6.pt\int}_{\H^2}\left| f(a)\right| \left| [\overline 
z,a,\overline
a,\xi]\right| ^t\left| d(z,\xi)\right| ^t\,\d\nu_0(a,\xi)\cr
\noalign{\medskip}
&\qquad =\displaystyle \mathop{\int\kern-6.pt\int}_{\H^2}\left| f(a)\right| 
( d(\overline z,a))^t
( d(\overline a,\xi))^t\,\d\nu_0(a,\xi)\cr
\noalign{\medskip}
&\qquad =\displaystyle \int_{\H} f(a)( d(\overline z,a))^t
\left(\int_{\H} d(\overline a,\xi)^t\,\d\nu_0(\xi)\right)\,\d\nu_0(a).}$$
But the inner integral is a constant $K_t$, depending just on $t$ and 
not on $z$.
Thus we get
$$K_t\int_{\H} f(a)( d(\overline z,a))^t\,\d\nu_0(a)\leq K_t^2\|f\|_\infty.
\eqno{\qed}
$$
\hfill\qed

\vfill\eject



\centerline{\S 5. {\ninebf Comparison of $T^t_{\ln\varphi}$ and
$\Lambda (1)(\overline z,\xi)={}$}{\ninerm ``$\ln\varphi (\overline 
z,\xi)-\left( c'_t/c_t\right) $''}}\vskip12pt


In this section we compare $\Lambda(1)$, which was constructed in Section 3,
  with $T^t_{\ln\varphi}$.

  We recall that $\Lambda(1)$
is  (up to an additive constant depending on the deformation parameter $t$)
$$(S_{\Delta^{\varepsilon_0}})^{-1}\left(\left.\displaystyle 
{\d\;\over{\d\varepsilon}}
S_{\Delta^{\varepsilon}}
S^*_{\Delta^{\varepsilon}}\right|_{\varepsilon=
\varepsilon_0\atop{\varepsilon>\varepsilon_0}}\right)
\left(\left(S_{\Delta^{\varepsilon_0}}\right)^*\right)^{-1}$$
where $S_{\Delta^{\varepsilon}}$ is acting on $H_{t+12\varepsilon_0}$,
  while $S_{\Delta^{\varepsilon_0}}$ acts  from $H_t$ into
$H_{t+12\varepsilon_0}$. The inverse $(S_{\Delta^{\varepsilon_0}})^{-1}$ is  an
   unbounded operator with domain dense in the closure of the range of
$S_{\Delta^{\varepsilon_0}}$.
We have explained in Section 4 that
$\Lambda(1)$ corresponds, in a non-specified
way, to the kernel $\ln \varphi (\overline z,\xi)-
\left( c'_t/c_t\right) $.

Both $\Lambda(1)$ and $T_{\ln\varphi }^t$ are  positive and 
affiliated with $\ac_t$.
Also recall, from Section 3, that the above definition for $\Lambda(1)$ 
translates into
the fact that for
  $\wc=\bigcup_{\varepsilon_0}\Range 
\left(S^t_{\Delta^{\varepsilon_0}}\right)^*$
we have that (up to a constant)
$$\Lambda(1)={\d\;\over{\d\varepsilon}}
\left\langle 
S^t_{\Delta^{\varepsilon}}\left(S^t_{\Delta^{\varepsilon}}\right)^*
w,w\right\rangle_{H_t}=
\lim\limits_{\varepsilon\searrow 0}
\left\langle{S_{\Delta^{\varepsilon}} S^*_{\Delta^{\varepsilon}}-
\Id \over{\varepsilon}}w,w\right\rangle$$

Our main result proves that there exists a (possibly different) domain where
$T_{\ln \varphi }^t$ is given by the same formula.

The main result is as follows:\vskip6pt

{\sm Proposition} 5.1.\kern.3em {\it There exists a densely defined ${\cal 
S}_0\subseteq H_t$, which
is a core for $T_{\ln \varphi }^t$
{\rm (}though not affiliated with $\ac_t${\rm )}
such that the following holds true:

Let $G_\varepsilon$ be the bounded operator in $\ac_t$ given by
$\left( 1/\varepsilon \right) \left(S_{\Delta^{\varepsilon}}
S^*_{\Delta^{\varepsilon}}-\Id \right)$.
Clearly $G_\varepsilon$ has kernel $\widehat {G}_\varepsilon
  (\overline z,\xi)=
\left( 1/\varepsilon \right) \left( \left( c_t / (c_t+\varepsilon ) \right) 
\varphi (\overline z,\xi)-\Id \right)$,  and the kernels converge pointwise
{\rm (}as $\varepsilon$ tends to $0${\rm )} to $\ln\varphi
(\overline z,\xi)-\left( c'_t/c_t\right) $.

Then, for all $v_1, v_2$ in ${\cal S}_0$, we have that
$$\langle T_{\ln\varphi }v_1,v_2\rangle=
\lim\limits_{\varepsilon\searrow 0}\langle G_\varepsilon v_1,v_2\rangle.$$
}\vskip6pt

{\it Remark.} By comparison, the same holds true for $\Lambda(1)$:
the only difference is that this happens on a different
domain ${\cal W}$ (in place of ${\cal S}_0$) which is affiliated to
$\ac_t$.\vskip6pt

This will be proved in several steps, divided into
the following lemmas.\vskip6pt

{\sm Lemma} 5.2.\kern.3em
  {\it Let $${\cal S}=\displaystyle
\left\{\sum\limits_{i=1}^N{\lambda_i\over{(z-\overline
a_i)^{\alpha_i}}}\e^{\im\varepsilon_i z} \biggm|
\Re \alpha_i>3,\ \varepsilon_i>0,\ \lambda_i\in \C, 
N\in \N\right\}.$$ Then
${\cal S}$ is contained in all
$H_t$, and dense in all $H_t$, $t>1$.
}\vskip6pt

{\it Proof.}
Actually $\Re \alpha_i>1$
would  be sufficient for the convergence, but for later 
considerations we take $3$
instead  of $1$.
  It is sufficient to consider a single term (so $N=1$). We omit
all the indices for $\alpha,a,\varepsilon$ and let $\lambda=1$. We 
prove first that
$f(z)=\left( 1 / (z-\overline a)^\alpha \right) \e^{\im\varepsilon z}$
belongs to any $H_t$. Indeed we have
$$\displaystyle \int_{\H}\left|{1\over{(z-\overline a)^{ \alpha}}}
\e^{\im\varepsilon z}\right|^2\,\d\nu_t(z)
=\int_{\H} {1\over{\left| z-\overline a\right| ^{\Re
\alpha}}}
\e^{-(\Im z)\varepsilon}
(\Im z)^{t-2}\,\d\overline z\,\d z,$$
which is obviously convergent as $\Re \alpha\geq 2$.\hfill\qed\vskip6pt

In the next lemma, we enlarge that  space ${\cal S}$ to exhaust the 
range of all
$S_{\Delta^{\varepsilon}}$.\vskip6pt

{\sm Lemma} 5.3.\kern.3em
{\it Let ${\cal S}_{0,t}=
\bigcup_{\varepsilon>0}\Delta^\varepsilon{\cal
S}$, $t-\varepsilon>1$.
  Then ${\cal S}_{0,t}$
is dense in all $H_t$, $t>1$.}\vskip6pt

{\it Proof.} We need only look at ${\cal S}\subseteq H_{t-\varepsilon}$ and
  apply the operator $S_{\Delta^{\varepsilon}}$.\hfill\qed\vskip6pt

Next we need a bound on $\Im (\ln\Delta(z))$.
Recall that we are using a choice for $\ln\Delta(z)$ which comes from that
fact that $\Delta(z)$ is nonzero in $\H$.\vskip6pt

{\sm Lemma} 5.4.\kern.3em
{\it Let $\ln\Delta(z)$ be the principal branch of 
the logarithm of the
function $\Delta$. Then
$\left| \Im\ln(\Delta(z))\right| $ is bounded  by a constant times
$\left(\Re z+\left( 
1 / (\Im z)^2\right) \right)$,  as $\Im z\downarrow 0$.}\vskip6pt

{\it Proof.} We let $q=\e^{2\pi \im z}$ and use the following expansion for
$\ln\Delta(z)$:
$$\ln\Delta(z)={\pi \im z\over{12}}+\sum\limits_{n\ggeq 1}\ln(1-q^n).$$

When $r=\left| q\right| =\left| \e^{2\pi \im z}\right| =\e^{-\pi y}$ tends to $1$ we have, with 
$q=r\e^{\im\theta}$,
$z=x+\im y$, that
$$\eqalign{\displaystyle \Im\ln\Delta (z)&={\pi x\over{12}}+
\sum\limits_{n\ggeq 1}\arg((1-r^n\cos n\theta)+\im r^n\sin n\theta)\cr
\noalign{\medskip}
&={\pi x\over{12}}+
\sum\limits_{n\ggeq 1}\tan^{-1}
\left[{r^n\sin(n\theta)\over{1-r^n\cos n\theta}}\right].}$$

As $r\to 1$ this is dominated by
$${\pi x\over{12}}+\sum\limits_{n\geq 1}{r^n\sin 
n\theta\over{1-r^n\cos n\theta}},$$
which in turn is dominated by
$${\pi  x\over{12}}+\sum\limits_{n\ggeq 1}{r^n\over{1-r^n}}.$$
This turns out to be
$$\eqalign{\displaystyle {\pi x\over{12}}&+
(r+r^2+r^3+\cdots)\cr
\noalign{\medskip}
&\qquad +(r^2+r^4+r^6+r^8+\cdots)\cr
\noalign{\medskip}
&\qquad \qquad +(r^3+r^6+r^9+\cdots)\cr
\noalign{\medskip}
&\qquad \qquad \qquad +
(r^4+r^8+\cdots)\cr
\noalign{\medskip}
&\qquad \qquad \qquad \qquad +(r^5+r^{10}+\cdots)\cr
\noalign{\medskip}
&\qquad \qquad \qquad \qquad \qquad +\cdots}$$
and this is dominated by
$${\pi x\over{12}}+{c\over{(1-r)^2}},$$
for some constant $c$.

Letting $r=\e^{-2\pi y}$, and using that
$\lim\limits_{y\to 0}\left( (1-\e^{-2\pi y} ) / y\right) $ is finite,
it follows that
$$\left| \Im \ln\Delta(x)\right| \leq c\left(x+{1\over{y^2}}\right)=
c\left(\Re z+{1\over{(\Im z)^2}}\right) .\eqno{\qed}$$\vskip6pt

{\sm Corollary} 5.5.\kern.3em
{\it For any $\varepsilon>0$, there exists a 
constant $c_\varepsilon$ such
that
$$\left| \Delta ^\varepsilon(z)\ln\Delta (z)\right| 
\leq c_\varepsilon\left(1+{1\over{\Im z}}\right)\left(1+\Re z+
{1\over{(\Im z)^2}}\right) .$$
}\vskip6pt

{\it Proof.}
We write
$$\left| \Delta ^\varepsilon(z)\ln(\Delta (z))\right| \leq\left| \Delta 
\right| ^\varepsilon\ln\left| \Delta (z)\right| +
\left| \Delta (z)\right| ^\varepsilon\left| \Im(\ln\Delta
(z))\right| .$$
We note that
$\left| \Delta (z)\right| ^2\Im z^{12}$ is a bounded function and hence
$$\displaystyle \left| \Delta (z)\right| \leq{c_1\over{(\Im z)^6}}.$$
Also, since
$\left| x^\varepsilon\ln x\right| \leq\const\left(\left[x^{\varepsilon_1},
x^{\varepsilon_2}\right]\right)$ for
$x>0$, where $\varepsilon_1>\varepsilon>\varepsilon_2$, we have that
$$\eqalign{\displaystyle 
\left| \Delta(z)\right| ^\varepsilon\ln\left| \Delta(z)\right| &\leq\const
\left(\left| \Delta(z)\right| ^{\varepsilon_1},\left| \Delta(z)\right| ^{\varepsilon_2}\right)\cr
\noalign{\medskip}
\displaystyle &\leq c\max\left({1\over{(\Im 
z)^{6\varepsilon_1}}},
{1\over{(\Im z)^{6\varepsilon_2}}}\right)\leq c\left(1+{1\over{(\Im 
z)}}\right).}$$
Similarly,
$$\left| \Delta(z)\right| ^\varepsilon\left| \Im(\ln(\Delta (z)))
\right| \leq{c\over{(\Im 
z)^{6\varepsilon}}}
\left[x+{1\over{(\Im z)^2}}\right] .$$
Putting the two  inequalities together, we get
$$\left| \Delta ^\varepsilon(z)\ln(\Delta (z))\right| 
\leq c\left(1+{1\over{(\Im z)}}+
{1\over{(\Im z)^{6\varepsilon}}}\left[x+{1\over{(\Im z)^2}}\right]\right),$$
which is thus smaller than
$$\eqalign{\displaystyle &c\left(\left(1+{1\over{(\Im z)}}\right)+
x\left(1+{1\over{(\Im z)}}\right)+
\left(1+{1\over{(\Im z)}}\right){1\over{(\Im z)^2}}\right)\cr
\noalign{\medskip}
\displaystyle &\qquad \qquad \qquad \qquad \qquad \qquad 
=c\left(1+{1\over{(\Im z)}}\right)\left(1+\Re 
z+{1\over{(\Im
z)^2}}\right) .}\eqno{\qed}$$\vskip6pt

{\sm Corollary} 5.6.\kern.3em
{\it Because $\left| \Delta(z)\right| $ has the order
of growth of $\left| \e^{2\pi\im z}\right| =\e^{-2\pi y}$, $y=\Im z$,
  it follows,  by first splitting
$\Delta ^\varepsilon(z)=\Delta ^{\varepsilon_1}(z)\Delta ^{\varepsilon_2}(z)$,
  that the growth of
$\left| \Delta ^\varepsilon(z)\ln\Delta (z)\right| $ will
come from $1 / \Im z$ as $\Im z\to 0$.
  Thus the above estimate can be improved
to
  $$\left| \Delta ^\varepsilon(z)\ln\Delta (z)\right| \leq
c{\Re z\over{(\Im 
z)^3}}\left(\e^{-\varepsilon_1\Im z}\right).$$
}\vskip6pt

In the next lemma we establish the integral formula for $\langle 
T^t_{\ln \varphi} v,v\rangle$.\vskip6pt

{\sm Lemma} 5.7.\kern.3em
{\it Fix $t\geq 10$. For $v$ in ${\cal S}_{0,t}$,  the integral
$$\mathop{\int\kern-6.pt\int}\limits_{\H^2}
{\ln\varphi({\overline z},\xi)\over{(\overline 
z-\xi)^t}}v(z)\overline{v(\xi)}\,\d
\nu_t(z,\xi)$$
is absolutely convergent and equal to
$$ \int_\H \ln \varphi(\overline z, z)\left| v(z)\right| ^2 \,\d \nu_t (z)=
\langle T^t_{\ln \varphi} v,v\rangle.$$}\vskip6pt

{\it Proof.} We will make use of the fact that $v\in {\cal S}_{0,t}$, so that
$$v(z)=\Delta^\varepsilon(z) v_1(z)$$
for some $\varepsilon >0$ and for some $v_1 \in {\cal S}$
(which is contained in $H_{t-\varepsilon}$).

We start by establishing the absolute convergence of the integral.
The integral of the absolute value of the integrand is
$$\mathop{\int\kern-6.pt\int}\limits_{\H^2}
{\left| \ln\varphi(\overline z,\xi)\right| 
\left| \Delta ^\varepsilon(z)\right| \over{\left| \overline z-\xi\right| ^t}}\left| v_1(z)\right| \left| \Delta 
^\varepsilon(\xi)\right| 
\left| v_1(\xi)\right| \,\d
\nu_t(z,\xi).$$

We expand this into three terms, by using the expression
$$\ln \varphi(\overline{z},\xi)=
\ln\overline{\Delta(z)}+\ln \Delta(\xi)+12\ln (\overline z-\xi),\qquad 
{\rm for\ } z,\xi\in \H.$$
We will analyse each term separately. Since the situations are similar we will
  do only the computation
for the term involving $\left| \ln\Delta (z)\right| $.
The corresponding integral is
$$\mathop{\int\kern-6.pt\int}\limits_{\H^2}
{\left| \ln\Delta (z)\right| \left| \Delta ^\varepsilon(z)\right| 
\over{\left| \overline z-\xi\right| ^t}}
\left| v_1(z)\right| \left| \Delta ^\varepsilon(\xi)\right| 
\left| v_1(\xi)\right| \,\d
\nu_t(z,\xi) .\eqno{(5.1)}$$

Because $(\Im z)^{t/2}(\Im\xi)^{t/2}
/ \left| \overline z-\xi \right| ^t$
is bounded  above $1$, the previous integral is   in turn bounded by the integral
$$\mathop{\int\kern-6.pt\int}\limits_{\H^2}
\left| \ln \Delta (z)\Delta 
^\varepsilon(z)\right| \left| v_1(z)\right| 
\left| \Delta ^\varepsilon(\xi)\right| 
\left| v_2(\xi)\right| \left( \Im z\right) ^{t/2-1}
\left( \Im \xi\right) ^{t/2-1}
\,\d z\,\d\overline z\,\d\xi\,\d\overline \xi .$$

We  use the estimate from Corollary 5.6 to obtain that this integral 
is further bounded
(up to a multiplicative constant $c$) by
$$c\mathop{\int\kern-6.pt\int}\limits_{\H^2}{\Re z\over{(\Im 
z)^3}}\left| v_1(z)\right| \left| v_2(\xi)\right| 
\e^{-\varepsilon_1\Im z}\e^{-\varepsilon(\Im \xi)}(\Im z)^{t/2-2}(\Im 
\xi)^{t/2}
\,\d z\,\d\overline z\,\d\overline \xi\,\d\xi.
$$
This comes to
$$c\mathop{\int\kern-6.pt\int}\limits_{\H^2}(\Re z)\left| v_1(z)\right| 
\left| v_2(\xi)\right| 
\e^{-\varepsilon_1\Im z}\e^{-\varepsilon(\Im \xi)}(\Im 
z)^{t/2-5}(\Im \xi)^{t/2}
\,\d z\,\d\overline z\,\d\overline \xi\,\d\xi.$$
As long as $t/2-5\geq 0$, the term $\e^{-\varepsilon_1(\Im z)}(\Im 
z)^{t/2-5}$ will be
bounded by some $\e^{-\varepsilon'_1(\Im z)}$.

Thus if $t\geq 0$, and with the price of replacing 
$\varepsilon,\varepsilon_1$ with some smaller
ones, in order to kill growth of $(\Im z)^{t/2-5}$ and 
$(\Im\xi)^{t/2-2}$, we get
a multiple of
$$\mathop{\int\kern-6.pt \int}\limits_{\H^2}(\Re z)\left| v_1(z)\right| 
\left| v_2(\xi)\right| 
\e^{-\varepsilon_1\Im z}\e^{-\varepsilon_2\Im \xi}
\,\d z\,\d\overline{z}\,\d\xi \,\d\overline{\xi}.$$

But for $z=x+\im y$, $|v_1(z)|$ involves powers of $
1 / x^3$, which makes
the integral absolutely convergent.  Hence the integral in (5.1) is 
absolutely convergent.

In the next lemma we will prove that for $v\in {\cal S}_0$, the integral
$$c_t\mathop{\int\kern-6.pt \int}\limits_{\H^2}
{ \ln \zxi\over{ \zxi^t}}v_1(z)\overline{v_2(\xi)}
\,\d\nu_t(\overline {z},\xi)$$
is absolutely convergent and equal to
$$\int v_1(z)\overline{v_2(\xi)}\ln (\overline{z}-z)\,\d\nu_t(z)+
{ c_t'\over{ c_t}}.$$

We now complete the proof of Lemma~5.7:
$$\displaystyle \int_{\H}\ln\varphi (\overline {z},z)
\left| v(z)\right| ^2\,\d\nu_t(z)
=\int[\overline{\ln \Delta (z)}+\ln \Delta (z)+12\ln (\overline {z}-z)]
\left| v(z)\right| ^2\,\d\nu_t(z).$$
We analyze each term separately.

We have
$$\eqalign{\displaystyle &\int \ln \Delta (z)v(z)\overline{v(z)}\,\d\nu_t(z)=
\int \ln \Delta (z)\overline{v(z)}c_t
\left(\int{ v(\xi)\over{ 
(z-\overline{\xi})^t}}\,\d\nu_t(\xi)\right)\,\d\nu_t(z)\cr
\noalign{\medskip}
\displaystyle &\qquad =c_t\mathop{\int\kern-6.pt \int}\limits_{\H^2}
{ \ln (\Delta (z))\overline{v(z)}v(\xi)\over{ (z-\overline 
{\xi})^t}}\,\d\nu_t(z,\xi)=
c_t\mathop{\int\kern-6.pt \int}\limits_{\H^2}
{ \ln \Delta (\xi)\overline{v(\xi)}v(z)\over{ \zxi^t}}
\,\d\nu_t(z,\xi).}$$

Similarly,
$$\int \ln \overline{\Delta (z)}v(z)\overline{v(z)}\,\d\nu_t(z)=
c_t\mathop{\int\kern-6.pt \int}\limits_{\H^2}
{\overline{ \ln \Delta (z)}\overline{v(\xi)}v(z)\over{ \zxi^t}}
\,\d\nu_t(z,\xi).$$

We know that the integrals are absolutely convergent and that we may 
integrate  in any
order.
Finally using the next lemma, we will have that
$$\int_\H\ln (\Im z)\left| v(z)\right| ^2\,\d \nu_t(z)
=
c_t\mathop{\int\kern-6.pt \int}\limits_{\H^2}
{ \ln \zxi\over{ \zxi^t}}v(z)\overline{v(\xi)}
\,\d\nu_t(z,\xi)-{ c_t'\over{ c_t}}\langle v,v\rangle _{H_t}.$$

Putting this together we get that
$$\eqalign{&\langle T_{(1/12)\ln\varphi }v,v\rangle\cr
\noalign{\medskip}
&\qquad =c_t\mathop{\int\kern-6.pt \int}\limits_{\H^2}
{ [1/12\ln (\overline{\Delta(z)}\Delta(\xi)\zxi^{12})]\over{
\zxi^t}}v(z)\overline{v(\xi)}
\,\d\nu_t(z,\xi)-{ c_t'\over{ c_t}}\langle v,v\rangle _{H_t}.}
$$
This completes the proof of Lemma~5.7.\hfill\qed\vskip6pt

The following  lemma was used above.\vskip6pt

{\sm Lemma} 5.8.\kern.3em
{\it For $v$ in ${\cal S}_{0,t}$, we have that 
$$\int_{\H}\ln (\Im z)\left| v(z)\right| ^2\,\d \nu_t=
c_t\mathop{\int\kern-6.pt \int}\limits_{\H^2}
{ \ln \zxi\over{ \zxi^t}}v(z)\overline{v(\xi)}
\,\d\nu_t(z,\xi)-{ c_t'\over{ c_t}}\langle v,v\rangle _{H_t}.
$$}\vskip6pt

{\it Proof.} Start with the identity
$$v(\xi)=c_s\int_{\H}{ v(z)\over{ \zxi^s}}\,\d \nu_t(z)=
\langle v,e^t_{\xi}\rangle _{H_t}.$$

We differentiate this with respect to $s$, at
$s=t$ (which is allowed because of the 
fast decay
of the functions in ${\cal S}_{0,t}$).

This gives us
$$0={ c_t'\over{ c_t}}v(\xi)+c_t\int_{\H}{ v(z)
[\ln (\Im z)-\ln \zxi]\over{ \zxi^s}}
\,\d \nu_t(s).$$

Now we integrate over $\H$, with respect to the measure 
$\overline{v(\xi)}\cdot\d \nu_t(\xi)$.

We get
$$\eqalign{0&=\displaystyle { c_t'\over{ c_t}}\|v\|_{H_t}^2+
c_t\mathop{\int\kern-6.pt
\int}\limits_{\H^2}{ v(z)\overline{v(\xi)}\ln \Im z\over{
\zxi^t}}\,\d \nu_t(z,\xi) \cr
\noalign{\medskip}
&\qquad- \displaystyle
\mathop{\int\kern-6.pt \int}\limits_{\H^2}{ \ln \zxi 
v(z)\overline{v(\xi)}\over{
\zxi^t}}\,\d \nu_t(z,\xi).}$$

The second integral is
$$\int_{\H}\left| v(z)\right| ^2\ln (\Im  z)\,\d \nu_t(z).$$

So we get the required identity.

This completes the proof of Lemma~5.7 and also the proof of
Lemma~5.8.\hfill\qed\vskip6pt

We now prove that  the reproducing kernel
$$\displaystyle {1\over{12}}\ln \left(\overline{\Delta (z)}\Delta(\xi)\zxi^{12}\right)-
{ c_t'\over{ c_t}}$$
is the derivative of
$\left( c_t / ( c_{t+12(s-t)} ) \right) S_{\Delta^{( s-t ) / 12}}
S^*_{\Delta^{( s-t ) / 12}}$
on the space ${\cal S}_{0,t}$.\vskip6pt

{\sm Lemma} 5.9.\kern.3em
{\it For $v_1,v_2\in {\cal S}_{0,t}$, we have that
$$\mathop{\int\kern-6.pt \int}\limits_{\H^2}{ \ln \varphi (
\overline z,\xi) v_1(z)\overline{v_2(\xi)}\over{
\zxi^t}}\,\d \nu_t(z,\xi)$$
is the limit, when $\varepsilon \searrow 0$, of 
$$\mathop{\int\kern-6.pt \int}\limits_{\H^2}{ 1/\varepsilon (\varphi
(\overline {z},\xi)^{\varepsilon }-\Id )\over{
\zxi^t}}v_1(z)\overline{v_2(\xi)}\,\d \nu_t(z,\xi).$$
}\vskip6pt

{\it Proof.} The convergence of the integrals involved in the limits 
was proved in
Lemma~5.7. To check the value of the limit we will evaluate the difference.
This is
$$\mathop{\int\kern-6.pt \int}\limits_{\H^2}{ [1/\varepsilon (\varphi
(\overline {z},\xi)^{\varepsilon }-\Id )-\ln \varphi (\overline {z},\xi)]
\over{ \zxi^t}}v_1(z)\overline{v_2(\xi)}\,\d
\nu_t(z,\xi).$$

We use the Taylor formula to express
$$\eqalign{\displaystyle
&\mathop{\int\kern-6.pt \int}\limits_{\H^2}[1/\varepsilon (\varphi
(\overline {z},\xi)^{\varepsilon }-1)-\ln \varphi (\overline {z},\xi)]
v_1(z)\overline{v_2(\xi)}\,\d \nu_t(z,\xi)\cr
\noalign{\medskip}
\displaystyle &\qquad = 
\mathop{\int\kern-6.pt \int}\limits_{\H^2}\varepsilon  \int_0^1
{\varphi ^{\varepsilon r}(\overline {z},\xi)\ln ^2\varphi (\overline 
{z},\xi)\over{
\zxi^t}}v_1(z)\overline{v_2(\xi)}\,\d r\,\d \nu_t(\overline{z},\xi).}$$

The same type of argument as in Lemma~5.7, because
  of the rapid decay of the vectors $v_1,v_2$
in ${\cal S}_{0,t}$, proves that the integral
$$\int_0^1\mathop{\int\kern-6.pt \int}\limits_{\H^2}{ \varphi 
^{\varepsilon r}
(\overline {z},\xi)\ln ^2\varphi (\overline {z},\xi)\over{
\zxi^t}}v_1(z)\overline{v_2(\xi)}\,\d \nu_t(z,\xi)\,\d r$$
is absolutely convergent with a bound independent of $\varepsilon $.
This completes the proof of Lemma~5.9.\hfill\qed\vskip6pt

To complete the proof of  Proposition 5.1, it remains to check the 
fact that the operators
$G_{\varepsilon }=
\left( S_{\Delta^\varepsilon} S_{\Delta^\varepsilon}^* -\Id \right) / 
\varepsilon $
are decreasing, after making a correction of the form
$-G_{\varepsilon }+\varepsilon K$, for a constant $K$.
This is done in the following lemma.\vskip6pt

{\sm Lemma} 5.10.\kern.3em
{\it Consider the bounded operators
$G_{\varepsilon }=
\left( S_{\Delta^\varepsilon} S_{\Delta^\varepsilon}^* -\Id \right) / 
\varepsilon $,
  which are represented by the
kernels
$${1\over \varepsilon}
\left\lbrack{ c_{t-12\varepsilon }\over{ c_t}} \varphi (\overline 
z,\xi)^\varepsilon-1\right\rbrack.
$$
Then, there exists a constant $K$ such that $-G_\varepsilon +
K\varepsilon$ is {\rm (}as $\varepsilon$ decreses to $0${\rm )} an increasing
  family of positive operators in
$\ac_{2t+1}$. }\vskip6pt

{\it Proof.} Note that the kernel of
$S_{\Delta^{\varepsilon }}S^*_{\Delta^{\varepsilon }}$ is
$$\displaystyle { c_{t-12\varepsilon }\over{ 
c_t}}\left[ \overline{\Delta(z)}\Delta (\xi)
\zxi^{12}\right] ^{\varepsilon }$$
Hence the derivative is
$$\displaystyle -12{ c_{t}'\over{ c_t}}+\ln
\left[ \overline{\Delta (z)}\Delta (\xi)
\zxi^{12}\right] .$$

Clearly $S_{\Delta^{\varepsilon }}S^*_{\Delta^{\varepsilon }}$ is
a decreasing family. We will proceed as in Lemma~3.3.

Let $s_{\varepsilon }=S(\varepsilon)$ be the kernel of
$S_{\Delta^{\varepsilon}} S_{\Delta^{\varepsilon}}^*$ (as an operator 
on $H_t$).
Then
$$s_{\varepsilon }(\overline {z},\xi)={ c_{t-12\varepsilon }\over{ c_t}}
(\varphi (z,\xi))^{\varepsilon }.$$
Let $G_{\varepsilon }=
\left( s_{\varepsilon }-\Id \right) / \varepsilon $.

The first derivative of $s_{\varepsilon }$ (with respect to $\varepsilon $) is
$$-12{ c_{t}'\over{ c_t}}\varphi ^{\varepsilon }+
{ c_{t-12\varepsilon }\over{ c_t}}\varphi ^{\varepsilon }\ln \varphi .$$

The second derivative is
$$24{ c_{t}'\over{ c_t}}\varphi ^{\varepsilon }\ln\varphi -
{ c_{t-12\varepsilon }\over{ c_t}}\varphi ^{\varepsilon }(\ln \varphi )^2.$$

This is equal to (as $c_t'=1$)
$${ c_{t-12\varepsilon }\over{ c_t}}\varphi ^{\varepsilon }
\left[(\ln \varphi )^2-{ 24\ln \varphi \over{ c_{t-12\varepsilon }}}\right]=
{ c_{t-12\varepsilon }\over{ c_t}}\varphi ^{\varepsilon }
\left[(\ln \varphi )+{ 12\ln \varphi \over{ c_{t-12\varepsilon }}}\right]^2-
{ c_{t-12\varepsilon }\over{ c_t}}\varphi ^{\varepsilon }
{ 144\over{ (c_{t-12\varepsilon })^2}}.$$

This is further equal to
$${ c_{t-12\varepsilon }\over{ c_t}}\varphi ^{\varepsilon }
\left[\ln \varphi +{ 12\ln \varphi \over{ c_{t-12\varepsilon }}}\right]^2-
{ 144\over{ c_t(c_{t-12\varepsilon })}}\varphi ^{\varepsilon }.$$

For every $r>1$, $S^r_{\Delta^\varepsilon} 
\left(S^r_{\Delta^\varepsilon} \right)^*=
f(\varepsilon )$ is a decreasing family in $\ac_{r}$.
By evaluating the kernel, which is
$$f(\varepsilon )(\overline {z},\xi)={ c_{r-12\varepsilon }\over{ c_r}}
[\varphi (z,\xi)]^\varepsilon ,$$
we get that ${ \d\;\over{ \d\varepsilon }}
f(\varepsilon )(\overline {z},\xi)$
is a positive kernel for $\ac_{r}$. Since
$$\displaystyle c_{r-12\varepsilon }={ r-12\varepsilon -1\over{ \pi}},$$
we obtain that
$${ c_{r-12\varepsilon }\over{ c_r}}\varphi ^{\varepsilon }\ln \varphi -
12{ c'_r\over{ c_r}}\varphi ^{\varepsilon }$$
represents a negative kernel for $\ac_{r}$.

We recall, from Section 3, that a kernel $k=k(\overline{z},\xi)$ is
positive for $\ac_r$ (even if $k$ does not necessary represent an
operator in $\ac_r$) if
$\left[
\displaystyle{
k(\overline{z_i}, z_j)
\over{
(\overline{z_i}-z_j)^r}}
\right]_{i,j=1}^N$ is a positive
matrix for all choices of $z_1,z_2,\ldots ,z_N $ in $H$, and
for all $N$ in $\N$.

We get that
$$\varphi^\varepsilon\left(
\ln \varphi-
\displaystyle{
1\over{c_{r+12\varepsilon}}}\right)
$$
represents a negative (nonpositive) kernel for $\ac_r$.

Thus $\varphi^{\varepsilon/2}
\left[\ln \varphi+
\left( 12 / \left( r-12
\left( \varepsilon /2\right) -
\left( \varepsilon /2\right) -1
\right) \right) 
\right]$
is negative for $\ac_{r+12
{\varepsilon\over{2}}}$,
and hence the square
$$\varphi^{\varepsilon}
\left[\ln \varphi+
\displaystyle{12\over{\displaystyle r-12
{\varepsilon \over{2}}-
{\varepsilon\over{2}}-1
}}
\right]^2=
\varphi^{\varepsilon}
\left[\ln \varphi
\displaystyle{12\over{\displaystyle r-12
{\varepsilon \over{2}}-
{\varepsilon\over{2}}-1
}}
\right]^2
$$
is positive for $\ac_{2r-12\varepsilon}$.

Consequently, the kernels
$$\displaystyle{c_t-12\varepsilon\over{
c_t}}\varphi^{\varepsilon}
\left[\ln \varphi-
\displaystyle{12\over{\displaystyle t-12
\varepsilon-1}}
\right]^2$$
are positive for $\ac_{2t+13\varepsilon}$.

Now we note the trivial calculus formulae
$$G_\varepsilon=
\displaystyle{S(\varepsilon)-\Id \over{
\varepsilon}}=\int_0^1S'(\varepsilon v)\,\d v,$$
$$G_{\varepsilon'}=
\displaystyle{S(\varepsilon')-\Id \over{
\varepsilon'}}=\int_0^1S'(\varepsilon' v)\,\d v.$$

The above equalities hold pointwise, that is, when evaluating the 
corresponding kernels
on points in $\H^2$.
Hence $$\eqalign{
G_\varepsilon-G_{\varepsilon'}&=
\int_0^1(S'(\varepsilon v)-S'(\varepsilon'v))\,\d v\cr
\noalign{\medskip}
&=\int_0^1(\varepsilon v-\varepsilon' v)
\int_0^1 S''(p(\varepsilon v)+(1-p)\varepsilon 'v)\,\d p \,\d v\cr
\noalign{\medskip}
&=(\varepsilon-\varepsilon')
\int_0^1\int_0^1 v S''(\alpha (v,p))\,\d p\,\d v,}$$
where $\alpha (v,p)=p(\varepsilon v)+(1-p)\varepsilon'v\leq
\max(\varepsilon,\varepsilon')$.

We have proved  that $S''(\alpha(v,p))$ is represented by a
positive kernel $R$, from which one has to subtract a quantity $Q$
(which is precisely $\displaystyle 
{\const\over {c_{v-12\alpha(v,p)}}}
\varphi^{\alpha(v,p)}$).

As such, by integration we obtain
$$G_\varepsilon-G_{\varepsilon'}=(\varepsilon-\varepsilon')[\rc-{\cal Q}],$$
where $\rc$ represents a positive kernel
for $\ac_{2t- 12\min(\varepsilon,\varepsilon')}$.
Moreover, $Q$ is a positive element in
$\ac_{2t}+ 12\min(\varepsilon,\varepsilon')$, and
   ${\cal Q}$ is bounded by $c\cdot\Id$, where
   $c$ is a universal constant.

Assume that $\varepsilon\geq \varepsilon'$; then in the sense
of inequalities in $\ac_{2r+1}$ we have that
$$G_\varepsilon-G_{\varepsilon'}\geq (\varepsilon-\varepsilon')(-{\cal Q}).$$

Since $0\leq {\cal Q}\leq c\cdot\Id  \cdot \Id $, we have
that $0\geq -{\cal Q}\geq c\cdot\Id 
\cdot -\Id $ (in $\ac_{2t+1}$). Consequently,  in $\ac_{2t+1}$, 
we have that
$$G_\varepsilon-G_{\varepsilon'}\geq(\varepsilon-\varepsilon')(-c).$$
Therefore, the following inequality holds in $\ac_{2t+1}$:
$$G_\varepsilon+\varepsilon c\geq G_{\varepsilon'}+\varepsilon'c.$$
If we take into account that $G_\varepsilon$ was
negative and replace $G_\varepsilon$ by
$H_\varepsilon=-G_\varepsilon$, then we get that
in $\ac_{2t+1}$ 
$$(-G_\varepsilon)-\varepsilon c\leq(-G_{\varepsilon'})-\varepsilon' c ,$$
i.e., that if $\varepsilon\geq \varepsilon'$,
$$H_{\varepsilon}-\varepsilon c\leq H_{\varepsilon'}-\varepsilon c.$$

We have consequently proved that, in  $\ac_{2t+1}$,
the kernels
$$H_\varepsilon(\overline{z},\xi)=-G_\varepsilon(\overline{z},\xi)=
-\displaystyle{\displaystyle{
c_t- 12\varepsilon\over{
c_t}}\varphi^\varepsilon+\Id 
\over{\varepsilon}}$$ are positive and they increase
(when $\varepsilon$ decreases to zero, modulo an infintesimal term)
to $-\left( c_t' / c_t\right) 
+\ln \varphi(\overline{z},\xi)$.\hfill\qed\vskip6pt

{\sm Lemma} 5.11.\kern.3em {\it Let $M\subseteq B\left( H\right) $ be a type
$\twoone$ factor and assume that
$\left( H_{n}\right) _{n\in \N}$ is an
increasing family of positive operators in $M$. Let
$${\cal D}\left( X\right) =
\left\{ \xi \in H \mid
\sup \left\langle H_{n}\xi ,\xi \right\rangle <\infty\right\} $$
and assume that ${\cal D}\left( X\right) $ is weakly dense in $H$.
Then ${\cal D}\left( X\right) $ is affiliated with $M$,
and 
$\left\langle X\xi ,\xi \right\rangle
= \sup_{n}\left\langle H_{n}\xi ,\xi \right\rangle $,
$\xi \in {\cal D}\left( X\right) $,
defines an operator affiliated with $M$.}\vskip6pt

{\it Proof.} Clearly ${\cal D}\left( X\right) 
=\left\{ \xi \mid \sup_{n}\left\| \smash{H_{n}^{1/2}}\xi \right\| 
<\infty \right\}$
and as such ${\cal D}\left( X\right) $ is a subspace, because if 
$\left\| \smash{H_{n}^{1/2}}\xi \right\| \leq A$,
$\left\| \smash{H_{n}^{1/2}}\eta \right\| \leq B$ for all $n$
then $\left\| \smash{H_{n}^{1/2}}\left( \xi +\eta \right) \right\| \leq A+B$.

Moreover, ${\cal D}\left( X\right) $ is clearly invariant under
$u' \in {\cal U}\left( M'\right) $,
and hence ${\cal D}\left( X\right) $ is affiliated to $M$.

The quadratic linear form $q_{X}\left( \xi \right) =\sup_{n} 
\left\langle H_{n}\xi ,\xi \right\rangle $
is weakly lower semicontinuous, thus
by [Co2],
$q_{X}$ defines a positive unbounded operator $X$,
affiliated with $M$,
with domain ${\cal D}\left( X\right) $.\hfill\qed\vskip6pt

{\sm Corollary} 5.12.\kern.3em
{\it The following holds:
$$
T_{\ln \varphi}^{t}+{c_{t}'\over c_{t}}\cdot \Id
=\Lambda\left( 1\right) .
$$
}\vskip6pt

{\it Proof.} Let
$H_{\varepsilon}=
-G_{\varepsilon}+K\cdot \Id+C\varepsilon $, where
$G_{\varepsilon}$ are as in Lemma~5.10.
Then by definition $X=-\Lambda\left( 1\right) +K\cdot \Id$
coincides with the supremum of $H_{\varepsilon}$ on
${\cal S}_{0}=\bigcup S_{\Delta^{\varepsilon}}^{\ast}$.

On ${\cal S}_{0}$, which is a core for
$T_{\ln \varphi}^{t}+K\cdot \Id$,
the same holds for $T_{\ln \varphi}^{t}+K\cdot \Id$.

Thus $T_{\ln \varphi}^{t}|_{{\cal S}_{0}}\subseteq X$,
hence $\overline{T_{\ln \varphi}^{t}}\subseteq X$
and so $T_{\ln \varphi}^{t}=X=\Lambda \left( 1\right) $
by [MvN].\hfill\qed\vskip6pt


\vfill\eject

%


\centerline{\S 6. {\ninebf The cyclic cocycle associated to the deformation}}\vskip12pt

In [Ra] we introduced a cyclic cocycle $\Psi_t$, which lives on the
algebra $\bigcup_{s<t}\widehat {\ac}_s$,
and we proved a certain form of
nontriviality for this cocycle.

We recall first the definition of the cocycle $\Psi_t$ and then we 
will show the
nontriviality of $\Psi_t$  by using a quadratic form deduced from 
the operator introduced
in Corollary~3.6 and Lemma~4.4.
The main result of this section will be the following:\vskip6pt

{\sm Theorem} 6.1.\kern.3em
{\it Let $t>1$, let
  $ {{\cal B}}_t=
\bigcup_{s<t}\widehat {\ac}_s$, which is a weakly dense subalgebra
of $\ac_{t}$, and let $R_t$ be defined on ${\cal B}_t$ {\rm (}with values in 
${\cal B}_t${\rm )}
by the formula
$$\langle R_tk,l\rangle _{L^2(\ac _t)}
=-{1\over2}\tau _{\ac_{t}}({\cal C}_t(k,l^*)),\qquad
k,l\in {\cal B}_t$$
{\rm (}that is,
$R_t$ implements the Dirichlet form $\tau_t({\cal C}_t(k,l^*)${\rm ).}

Let $(\nabla R_t)(k,l)=R_t(k,l)-kR_tl-R_t(k)l$,
which belongs to ${\cal B}_t$, if $k,l\in {\cal B}_t$, and
  let $\Psi_t$ be the cyclic cocycle
associated with the deformation  {\rm [Ra]}
$$\Psi_t(k,l,m)=\tau_{\ac_{t}}([{\cal C}_t(k,l)-
(\nabla R_t)(k,l)]m),\qquad k,l,m\in {\cal B}_t.$$

  Let $\Lambda_0$ be the operator, on  the weakly dense
{\rm (}non-unital\/{\rm )} subalgebra
$\dc_t^0\subseteq\ac_{t}$, introduced in Theorem {\rm 4.6,} by requiring
that $\Lambda _0(k)$ is the derivative
at $0$ of the operator represented in $\ac_t$ by
  the kernel $\varphi ^{\varepsilon }(\overline
{z},\xi)k(\overline {z},\xi)$. Thus $\Lambda _0(k)(\overline {z},\xi)$
is formally $k(\overline
{z},\xi)
\ln\varphi (\overline {z},\xi)$.

Let $\chi_t(k,l)=\langle \Lambda _0k,l^*\rangle _{L^2(\ac_{t})}-
\langle k,\Lambda _0(l^*)\rangle _{L^2(\ac_{t})}$
be the antisymmetric form associated with $\Lambda _0$.
Then
$$\Psi_t(k,l,m)={ c_t'\over{ c_t}}\tau_{\ac_{t}}(klm)+\chi_t(k*_tl,m)+
\chi_t(l*_tm,k)+\chi_t(m*_tk,l)$$
for $k,l,m\in \dc_t^0$.}\vskip6pt

We will split the proof of this result into several steps:
First we prove some properties about $\Lambda_0$
and its formal adjoint $\Lambda ^+$.
We start with the definition of $\Lambda ^+$.
The first lemma collects the definition
and basic properties of $\Lambda ^+$.\vskip6pt

{\sm Lemma} 6.2.\kern.3em
{\it \raggedright
Let $f$ be a bounded measurable function that is
$\PSLtwoZ$-\break
equivariant.

We define $\Lambda ^+(T_f^t)=T^t_{f\ln\varphi }$. Then
$\Lambda ^+$ has the following properties:

{\rm 1)} Assume in addition that
  $f\ln \varphi (\overline {z},z)$
is  a bounded function.  Then\break
  $\left( \Lambda _0\mid\dc_t\right) ^*\subseteq \Lambda ^+$ and
$$\eqalign{\tau_{\ac_t} (\Lambda _0(k)(T_f^t)^*)=\tau_{\ac_t} 
(k\Lambda ^+(T_f^t))=
&\displaystyle {1\over \area F} \int_{F}k(\overline {z},z)f(z)
\ln\varphi (\overline {z},z)\,\d\nu_0(z).}$$

{\rm 2)} For $k,l$ in $\dc_t^0$, we have that
$$\tau (k\Lambda ^+(T_f^t)l)=\tau (\Lambda_0 (l*k)T_f^t).$$
}\vskip6pt

{\it Proof.} The proof of this proposition is obvious, since the integrals
are absolutely summable. For part~2 we remark that
$k\Lambda ^+(T_f^t)l$ has symbol
$$c_t\int_{\H}k(\overline {z},\eta)
[f(\overline {\eta},\eta)\ln \varphi (\overline {\eta},\eta)]
\varphi (\overline {\eta},\xi)\,\d\nu_t(\eta).$$

Hence by summability, the trace is
$${c_t\over \area F}\mathop{\int\kern-6.pt \int}
\limits_{F\times \H}k(\overline {z},\eta)l(\overline {\eta},z)
f(\overline {\eta},\eta)\ln\varphi (\overline {\eta},\eta)\,\d\nu_t(z,\eta),$$
which is exactly
$\tau\left(\Lambda (l*k)T^t_{f\ln\varphi }\right)$. This completes the
proof.\hfill\qed\vskip6pt

Recall that in Section 2 we introduced
  the densely defined operator ${\cal T}_{\ln d }$ on $L^2(\ac_t)$,
  given by the formula
$$\langle {\cal T}_{\ln d }k,l\rangle ={c_t\over \area F}
\mathop{\int\kern-6.pt \int}\limits_{F\times \H}k(\overline {z},\eta)
\overline{l(\overline {z},\eta)}
\ln d (\overline {z},\eta) 
\left| d (\overline {z},\eta)\right| ^{2t}\,\d\nu_0(t,\eta),$$
  which is well defined for $k,l$ in
algebra $\hat{\cal B}_t$. We note  that ${\cal T}_{\ln d }$ acts like a
Toeplitz operator on $L^2(\ac_{t})$, with
symbol $\ln d $. In the next lemma we establish the relation between 
the operator
${\cal T}_{\ln d }$ and the operator $R_t$.\vskip6pt

{\sm Lemma} 6.3.\kern.3em
{\it The operator $R_t$, defined by the property
$$\langle R_tk,l\rangle =-{1\over2}\tau ({\cal C}_t(k,l)),$$ has the
  following simple expression in terms
of ${\cal T}_{\ln d }$:
$$R_t=-{1\over2}{ c_t'\over{ c_t}}-\langle {\cal T}_{\ln d}k,l\rangle ,\qquad 
k,l\in \widehat {{\cal B}}_t.$$}\vskip6pt

{\it Proof.} Indeed, we have that, for
  $k,l\in  {{\cal B}}_t\subseteq \widehat {\ac}_t$,
$${\cal C}_t(k,l)=
{ c_t'\over{ c_t}}(k*_tl)+
c_t\int_{\H}k(\overline {z},\eta)l(\overline {\eta},\xi)
[\overline{z},\eta,\overline{\eta},\xi]^t
\ln [\overline{z},\eta,\overline{\eta},\xi]\,\d\nu_t(\eta).$$

If we make $\xi=z$ in the above expression and then
  integrate over $F$ to get the trace of
${\cal C}_t(k,l)$, we get
$$\eqalign{\tau({\cal C}_t(k,l))
&= { c_t'\over{ c_t}} \tau(k*_tl)+{ c_t\over{ \area F}}
\mathop{\int\kern-6.pt \int}
\limits_{F\times \H}k(\overline {z},\eta)\overline{l(\overline {z},\eta)}
\left| d (\overline {z},\eta)\right| ^{2t}
\ln d^2\,\d\nu_t(\eta)\cr
\noalign{\medskip}
&=\displaystyle { c_t'\over{ c_t}}\tau(k*_tl)+
2\langle {\cal T}_{\ln \varphi }k,l\rangle .}$$

This completes the proof.\hfill\qed\vskip6pt

In the next proposition we prove a relation
  between $\Lambda_0 +\Lambda ^+$
and the other terms (remark that $\Lambda ^+$
  is not necessarily
   the adjoint of $\Lambda_0 $, but rather
we define $\Lambda ^+(T_f^t)=T^t_{f\ln\varphi }$
whenever  possible).\vskip6pt

{\sm Proposition} 6.4.\kern.3em
{\it For all $k,l$ in $\dc_t^0$ we have
$$\langle \Lambda _0k,l^*\rangle +\langle k,\Lambda _0l^*\rangle =
\tau(k T^t_{f\ln\varphi }l)+ \tau(l T^t_{f\ln\varphi }k)-
2\langle T_{\ln\varphi }k,l\rangle .
\eqno{(6.1)}$$

Consequently if we define {\rm ``$\Re\Lambda _0$'' (}formwise\/{\rm )}
by the relation
$$\langle (\Re \Lambda _0)k,l\rangle =
{1\over2}(\langle \Lambda _0k,l\rangle +\langle k,\Lambda _0(l)\rangle),$$
then
$$\langle (\Re \Lambda _0)k,l^*\rangle =
{1\over2}\tau(kTl+lTk)+\langle R_tk,l^*\rangle +{1\over2}{ c_t'\over{ c_t}}
\langle k,l^*\rangle _{L^2(\ac_{t})}
\eqno{(6.2)}$$}\vskip6pt

{\it Proof.} We prove first the relation (6.1). For
$k,l\in \dc_t^0$,  we have that
$$\langle \Lambda _0k,l^*\rangle _{L^2(\ac_{t})}
+\langle k,\Lambda _0l^*\rangle _{L^2(\ac_{t})}$$
is equal to
$${ c_t\over{ \area F}}
\mathop{\int\kern-6.pt \int}\limits_{F\times \H}
[\ln \varphi (\overline {z},\eta)+\ln \varphi (\overline {\eta},z)]k
(\overline {z},\eta)l(\overline {\eta},z) d(\overline {z},\eta)^{2t}
\,\d\nu_0(z,\eta).$$
Since
$$\ln \varphi (\overline {z},\eta)+\ln \varphi (\overline {\eta},z)=
\ln \varphi (\overline {z},z)+\ln\varphi (\overline {\eta},\eta)-
\ln[ d(\overline {z},\eta)]^2,$$
we get the relation (6.1).

Dividing by $2$ we get
$${1\over{2}}[\langle \Lambda _0k,l^*\rangle+\langle k,\Lambda _0(l)^*\rangle ]=
{1\over{2}}[\tau(kT^t_{\ln\varphi }l)+\tau(lT^t_{\ln\varphi }k)]-
\langle T_{\ln d}k,l^*\rangle _{L^2(\ac_{t})}.$$

The definition of $R_t$ and the previous lemma complete the
proof.\hfill\qed\vskip6pt

Recall that in \S 4, we proved that for all $k,l$ in $\dc_t^0$,
$${\cal C}_t(k,l)={ c_t'\over{ c_t}}\Id +kT^t_{\ln\varphi }l+
\Lambda _0(kl)-\Lambda _0(k)l-
k\Lambda _0(l).
\eqno{(6.3)}$$

We want to use (6.3) to find an expression for
$${\cal C}_t(k,l)-(\Delta R_t)(k,l)$$
by taking the trace of the product of  $m\in \dc_t^0$ with the 
previous expression.\vskip6pt

{\it Notation}. We denote $T=T^t_{\ln \varphi}$ and let
$\langle \Sym\nolimits_\varphi k,l\rangle ={1\over2}[\tau(kTl^*)+\tau(l^*Tk)]$,
for $k,l\in \dc_t^0$. Hence $$\tau_{\ac_t}(\Sym\nolimits_\varphi(k)l)=
{1\over2}[\tau(kTl)+\tau(lTk)].$$\vskip6pt

  In this terminology the relation
  in Proposition~6.4 becomes
$$\langle (\Re \Lambda _0)k,l^*\rangle =
\langle \Sym\nolimits_\varphi k,l\rangle +\langle R_tk,l^*\rangle +{1\over{2}}
{ c_t'\over{ c_t}}\langle k,l^*\rangle .$$
Note that in the relation above, the scalar product refers to the scalar
product on $L^2(\ac_t)$. Moreover, the following relations hold true:
$$\tau(\Sym\nolimits_{\varphi }(kl)m)={1\over2}[\tau(klTm)+\tau (mTkl)],
\eqno{(6.4)}$$
$$\tau((\Sym\nolimits_{\varphi }k)lm)=\tau(\Sym\nolimits_{\varphi }(k)(lm))=
{1\over2}[\tau(klTm)+\tau (lmTk)],
\eqno{(6.5)}$$
$$\tau(k(\Sym\nolimits_{\varphi }(l))m)=\tau(\Sym\nolimits_{\varphi }(mk)l)=
{1\over2}[\tau(lTmk)+\tau (mkTl)].
\eqno{(6.6)}$$\vskip6pt

{\sm Lemma} 6.5.\kern.3em
{\it For all $k,l,m$ in $\dc_t^0$ we have that
$$E=\Sym\nolimits_{\varphi }(kl)-(\Sym\nolimits_{\varphi }k)l-
k(\Sym\nolimits_{\varphi }l)+
kT^tl=0.$$

To check this, one has to verify that
$\tau(Em)=0$ for all $m$ in $\dc_t^0$.
}\vskip6pt

{\it Proof.} We have to check  that the expression
$$\tau(klTm)+\tau (mTkl)-\tau(kTlm)-\tau (lmTk)-\tau(lTmk)-\tau (mkTl)+
2\tau(kTlm)$$
vanishes. But
$$\eqalign{\tau(mTkl)&=\tau (lmTk),\cr
\noalign{\medskip}
\tau(klTm)&=\tau (lTmk).}$$

After cancelling the above terms we are left to check that
$$-\tau(kTlm)-\tau (mkTl)+2\tau(kTlm)$$
is equal to zero, which is obvious since $k,l,m\in \dc_t^0$.
This completes the proof.\hfill\qed\vskip6pt

We now decompose $\tau (\Lambda _0(k)l)$ in the following way:
$$\tau (\Lambda _0(k),l^*)=
\langle (\Re \Lambda _0)(k),l\rangle +\im
\langle (\Im \Lambda _0)(k),l\rangle ,$$
where
$$\langle (\Im \Lambda _0)(k),l^*\rangle =(1/2\im)[
\langle \Lambda _0(k),l^*\rangle -\langle k,\Lambda _0(l^*)\rangle ].
\eqno{(6.7)}$$

We now can proceed to the proof of Theorem~6.1.\vskip6pt

{\it  Proof of Theorem} 6.1. We have
$${\cal C}_t(k,l)={ c_t'\over{ c_t}}kl+kTl+\Lambda _0(kl)-\Lambda _0(k)l-
k\Lambda _0(l).$$

Hence by taking
the
scalar product with an $m$ in $\ac_{t}$, that is, computing
$\tau ({\cal C}_t(k,l)m)$, we obtain
$$\eqalign{\displaystyle \tau ({\cal C}_t(k,l)m)&= { c_t'\over{ c_t}}
\tau(klm)+\tau (kTlm)+\tau(\Lambda _0(kl)m)-\tau(\Lambda _0(k)lm)-
\tau (\Lambda _0(k)mk)\cr
\noalign{\medskip}
\displaystyle &={ c_t'\over{ c_t}}
\tau(klm)+\tau ([kTl]m)+\langle \Re \Lambda _0(kl),m^*\rangle \cr
\noalign{\medskip}
&\qquad -
\langle \Re \Lambda _0(k),(lm)^*\rangle -
\langle \Re \Lambda _0(l),(mk)^*\rangle \cr
\noalign{\medskip}
&\qquad +\im \langle \Im \Lambda _0(kl),m^*\rangle -
\im\langle \Im \Lambda _0(k),(lm)^*\rangle -
\im\langle \Im \Lambda _0(l),(mk)^*\rangle .}$$

By using the relation
$$\langle \Re \Lambda _0(k),l^*\rangle =
\langle R_tk,l^*\rangle +{1\over{2}}{ c_t'\over{ c_t}}\langle k,l^*\rangle +
\langle \Sym\nolimits_{\varphi }k,l^*\rangle ,$$
we obtain that $\tau ({\cal C}_t(k,l)m)$ is equal to
$$\tau([(\Delta R_t)(k,l)]m)$$
plus the following terms:
$$\langle \Sym\nolimits_{\varphi }(kl),m^*\rangle -
\tau (\Sym\nolimits_\varphi (k)lm)- \tau (\Sym\nolimits_\varphi(l)mk),
\eqno{(6.8)}$$
plus the terms
$${ c_t'\over{ c_t}}\tau(klm)+\left(
{ 1\over{ 2}}{ c_t'\over{ c_t}}-
{ 1\over{ 2}}{ c_t'\over{ c_t}}-
{ 1\over{ 2}}{ c_t'\over{ c_t}}\right)\tau (klm),
\eqno{(6.9)}$$
plus the terms
$$\im \langle \Im \Lambda _0kl,m^*\rangle -
\im \langle \Im \Lambda _0(k),(lm)^*\rangle -
\im \langle \Im \Lambda _0(l),(mk)^*\rangle .
\eqno{(6.10)}$$

The terms in (6.8) add up to  zero, as it was  proved in Lemma~6.5. 
The terms in
(6.9) add up to  ${ 1\over{ 2}}\left( c_t' / c_t\right) \tau(klm)$.

Since $\chi_t(k,l)={ 1\over{ 2}}[\langle\Lambda _0k,l^*\rangle- \langle 
k,\Lambda _0l^*\rangle]=\im
\langle \Im \Lambda _0(k),l^*\rangle ,$
we obtain, by adding the terms  from (6.8), (6.9), (6.10),  that
$$\tau ({\cal C}_t(k,l)m)= \tau([(\Delta R_t)(k,l)]m)+
{ 1\over{ 2}}{ c_t'\over{ c_t}}\tau(klm)+
\chi_t(kl,m)-\chi_t(k,lm)- \chi_t(l,mk).
$$
Thus $$\Psi_t(k,l,m)=\displaystyle { 1\over{ 2}}{ c_t'\over{ c_t}}
+
\chi_t(kl,m)-\chi_t(k,lm)- \chi_t(l,mk).
\eqno{\qed}$$\vskip6pt

{\sm Lemma} 6.6.\kern.3em {\it Let $t>1$.
Assume that $k$, $l$ are such that
$k=k_{1} \ast_{t} k_{2}$,
$l=l_{1} \ast_{t} l_{2}$,
$k_{i},l_{i}\in{\cal D}_{t}^{0}$.
Then 
$$
\tau \left( k_{2}{\cal L}_{t}\left( l\right) k_{1}\right)
+\tau \left( l_{2}{\cal L}_{t}\left( k\right) l_{1}\right)
+\tau \left( {\cal C}_{t}\left( k,l\right) \right)
=-{c_{t}'\over c_{t}}\tau \left( k\ast_{t}l\right).
$$
}\vskip6pt

{\it Proof.} Recall that
$$
{\cal L}_{t}=\left( \Lambda_{0}-{c_{t}'\over c_{t}}\cdot
\Id\right) 
-{1\over 2}\left\{ T,\,\cdot \,\right\}.
\eqno{(6.11)}
$$
Also $\tau \left( {\cal C}_{t}\left( k,l\right) \right)
=\left( c_{t}' / c_{t}\right) \tau \left( k\ast_{t}l\right)
+2\left\langle {\cal T}_{\ln d}k,l\right\rangle$.
Also
$$
\tau \left( \Lambda_{0}\left( k\right) \cdot l\right) 
+\tau \left( k\Lambda_{0}\left( l\right) \right) 
=\tau \left( kTl \right) 
+\tau \left( lTk \right) 
-2\left\langle {\cal T}_{\ln d}k,l\right\rangle .
\eqno{(6.12)}
$$
Moreover,
$$
\tau \left( kTl \right) 
=\tau \left( \Lambda_{0}\left( l\ast_{t}k\right) \right) 
=\int_{F}\ln \varphi \cdot \left( l\ast_{t}k\right) \left( \bar{z},z\right) 
\,\d\nu_{0}\left( z\right) .
\eqno{(6.13)}
$$
Hence, by (6.11), 
$$
\tau \left( k_{2}{\cal L}_{t}\left( l\right) k_{1}\right) 
=\tau \left( k_{2}\Lambda_{0}\left( l\right) k_{1}\right) 
-{c_{t}'\over c_{t}}\tau \left( kl \right) 
-{1\over 2}\tau \left( k_{2}Tlk_{1} \right) 
-{1\over 2}\tau \left( k_{2}lTk_{1} \right) 
$$
and
$$
\tau \left( l_{2}{\cal L}_{t}\left( k\right) l_{1}\right) 
=\tau \left( l_{2}\Lambda_{0}\left( k\right) l_{1}\right) 
-{1\over 2}\left( -\tau \left( l_{2}Tkl_{1} \right) 
-\tau \left( l_{2}kTl_{1} \right) \right) 
-{c_{t}'\over c_{t}}\tau \left( kl \right) .
$$
So by (6.11), (6.12),
$$
\eqalign{
&\tau \left( k_{2}{\cal L}_{t}\left( l\right) k_{1}\right) 
+\tau \left( l_{2}{\cal L}_{t}\left( k\right) l_{1}\right) 
+\tau \left( {\cal C}_{t}\left( k,l\right) \right) 
\cr
&\qquad =\tau \left( \Lambda_{0}\left( k\right) l\right) 
+\tau \left( \Lambda_{0}\left( l\right) k\right) 
-2{c_{t}'\over c_{t}}\tau \left( kl \right) 
-\tau \left( kTl \right) 
-\tau \left( lTk \right) 
+\tau \left( {\cal C}_{t}\left( k,l\right) \right) 
\cr
&\qquad =\tau \left( kTl \right) 
+\tau \left( lTk \right) 
-2{c_{t}'\over c_{t}}\tau \left( kl \right) 
-\tau \left( kTl \right) 
-\tau \left( lTk \right) 
\cr
&\qquad \qquad -2\left\langle {\cal T}_{\ln d}k,l\right\rangle 
+2\left\langle {\cal T}_{\ln d}k,l\right\rangle 
+{c_{t}'\over c_{t}}\tau \left( kl \right) 
\cr
&\qquad =-{c_{t}'\over c_{t}}\tau \left( kl \right) .
}\eqno{\qed}
$$

\vfill\eject



\centerline{\S 7. {\ninebf A dual solution; closability
  of $\Lambda$ }}\vskip12pt


In this section we analyze the Hilbert-space dual of the operator $\Lambda (k),
k\in D_0^t$, introduced in Section 4. This is achieved by
analyzing the derivative of the one-parameter family of completely
positive maps $\chi_{s,t}\colon  \ac_t\rightarrow \ac_s$, $1<s\leq t$,
defined as follows:
$$\chi_{s,t}(k)=
S_{\Delta^{(s-t)/12}}k \left(S_{\Delta^{(s-t)/12}}\right)^*,\qquad
k\in \ac_t.$$
Recall that
the
Hilbert space $L^2(\ac_t)$ is naturally identified with the Hilbert
space of all kernels $k=k(\overline{z},\xi)$ on $\H\times \H$ that
are diagonally $\Gamma$-equivariant, $\Gamma=\PSLtwoZ$.
  The kernels are also required
to be square-summable with respect to the the measure
$[| d(z,\xi)|]^{2t} \,\d\nu_0(z) \,\d\nu_0(\xi)$ on $F\times \H$.
(Recall that $F$ is a fundamental domain for $\Gamma$ in $\H$.)

Consider the Hilbert space $L_t$ of all
measurable functions on $\H\times F$ that are square-summable with respect to
the measure $d^{2t}\,\d\nu_0\times \d\nu_0$. This space
is obviously identified with a space of $\Gamma$-invariant (diagonally)
functions on $\H\times \H$, square-summable over $F\times\H$.

We let $\pc$ be the orthogonal projection from $L^t\ \hbox{into}\ L^2(\ac _t)$.
Let $\Phi$ be a measurable (diagonally) $\Gamma$-equivariant function on
$\H\times \H$. With the above identification let $M_{\Phi}$
  be the (eventually unbounded) operator
on $L_t$, defined by multiplication with $\Phi$ on $L_t$. 
Correspondingly, there is
a Toeplitz operator $\tc_{\Phi}=\pc \mc_{\Phi} \pc$, densely defined 
on $L^2(\ac _t)$.

For example, the map $\Lambda$, constructed in  Section 4, is $\tc_{\Phi}$ with
$\Phi=\ln \varphi$. In Section 6, Lemma~6.3, we have proved that
   the operator $R_t$ defined
by $$\left\langle R_t k, l\right\rangle _{L^2(\ac_t)}=
\displaystyle -{1\over 2} {\d\;\over {\hbox{d}s}}\tau_{\ac_s}(k*_s
l^*),$$
defined  for  $k, l$ in an algebra, is exactly $-\tc _{\ln 
d}-{1\over 2}\left( c'_t
/ c_t\right) \Id $.

Let also $P_t$ be the orthogonal projection from $L^2(\H, \d\nu_t)$ onto $H_t$.
Recall that the formula for $\pi_t$ has a trivial extension
to a projective unitary representation, $\widetilde {\pi}_t$
(given by the same formula as $\pi_t$), on functions on $L^2(\H,\d\nu_t)$.
Moreover, $P_t=P_t\widetilde{\pi}_t
=\widetilde{\pi}_tP_t=P_t\widetilde{\pi}_tP_t$.
Let $\widetilde{\ac}_t\subseteq B(L^2(H,\nu_t))$
be the commutant of $\widetilde{\pi}_t(\Gamma)$. By [A],
this is a type $\twoinfty$ factor, such that $L^2(\widetilde{\ac}_t)$
is canonically identified with $L_t$. Consequently, at least, for
$k$ in $L^2(\widetilde{\ac}_t)\cap \widetilde{\ac}_t$, it
makes sense to consider $\pc(k)=P_tkP_t$.\vskip6pt

{\sm Lemma} 7.1.
{\it Let $\pc_t$ be the orthogonal projection from $L_t$
{\rm (}identified with $L^2(\widetilde{\ac}_t)${\rm )} into $L^2(\ac_t)$.

Then $\pc_t(k)$ is given by the formula $P_t k P_t$, which
is well defined for $k\in L^2(\widetilde{\ac}_t)\cap \widetilde{\ac}_t$
and then extended by continuity. For such a $k$, the kernel of
$(\pc_tk)(\overline{z},\xi)$, $z,\xi $ in $\H $, is given by
the formula
$$(\pc_t)(\overline{z},\xi)=
\displaystyle c_t^2\zxi^t
\displaystyle \int \kern-6pt \int\limits_{\H^2}
\displaystyle {k(\eta_1,\eta_2)
\over{
(\overline{z}-\eta_1)^t(\overline{\eta_1}-\eta_2)^t(
\overline{\eta_2}-\xi)^t}}\,\d \nu_t(\eta_1)\,\d \nu_t(\eta_2).$$}\vskip6pt

{\it Proof.} One can check imediately that the map 
$\pc(k)=P_tkP_t$,
for $k$ in $L^2(\widetilde{\ac}_t)$, defines
an orthogonal projection on $L^2(\ac_t)$.

The formula for ${\cal P}(k)$ follows by writing
down the corresponding kernels, and
it holds as long as $k$ is in $L^2(\widetilde{\ac}_t).$\hfill\qed\vskip6pt

In the next lemma we will prove that the Toeplitz operators of
  symbols  $\ln \varphi$ and $\ln \overline{\varphi}$ have dense domain
(in $L^2(\ac_t)$) and that they are adjoint to each other.\vskip6pt

{\sm Lemma} 7.2.\kern.3em
{\it Let $\varphi(\overline{z},\xi)=
{1\over {12}}\ln
[
\overline{\Delta(z)}\Delta(\xi)\zxi^{12}]$
{\rm (}as in Section {\rm 4).} Let $\widetilde{\dc}_t^0$ be   the union of
$S_{\Delta^\varepsilon}\widehat{\ac}_{t-3-\varepsilon}
S^{*}_{\Delta^\varepsilon}$
with respect to $\varepsilon>0$, $t-3-\varepsilon>1$.

Then $\Dom(\mc_{\varphi})\cap H_t$
contains the weakly dense subalgebra $ \widetilde{\dc}_t^0$.
Consequently $\widetilde{\dc}_t^0$ is contained in the domain of
$\pc \mc_{\varphi}\pc$, which is the Toeplitz operator $\tc_{\varphi}$.
}\vskip6pt

  Before beginning the proof of the lemma, we note the following consequence:\vskip6pt

{\sm Corollary} 7.3\kern.3em
{\it
The operator $\Lambda$ introduced in {\rm \S4
   (}restricted to  $\widetilde{\dc}_t^0${\rm )}  coincides
with $\mc_{
(1/12)\ln\varphi}k$, acting
on the same domain.
Morever, the operators $ \tc_{\ln \overline{\varphi}}$
and $\tc_{\ln {\varphi}}$    are densely
defined and $ \tc_{\ln \overline{\varphi}}\subseteq
\left(\tc_{\ln {\varphi}}\right)^*$,
$ \tc_{\ln {\varphi}}\subseteq \left(\tc_{\ln \overline{\varphi}}\right)^*$.
  Consequently, these operators are closable
in $L^2(\ac_t)$.
}\vskip6pt

{\it Proof.}
  The fact that $\Lambda (k)$ is equal to $\mc_{
(1/12)\ln\varphi}k$ for
$k$  in $\dc^0_t$ is a consequence of the fact that
$\mc_{\varphi}(k)$ is the $L^2$-valued derivative at $0$  of
  the differentiable,  $L_t$-valued
function
   $\varepsilon \rightarrow  \mc_{\varphi^\varepsilon}(k)$. This is based
on the arguments in the
proof (below) of Lemma 7.2. Hence, $
\tc_{\ln {\varphi}}(k)$ is the derivative at $0$ of the differentiable,
   $L^2(\ac_t)$-valued function
   $\varepsilon \rightarrow  \tc_{\varphi^\varepsilon}$.\hfill\qed\vskip6pt

{\it Proof of Lemma} 7.2. We have to check  that
for $k$ in $\widetilde{\dc}_t^0$ having the expression
$k=S_{\Delta^\varepsilon} k_1 S^{*}_{\Delta^\varepsilon}$, with
$k_1\in \widehat{\ac }_{t-2-\varepsilon}$ (so that up to a constant
$k(\overline{z},\xi)=\varphi^{\varepsilon}k_1(\overline{z},\xi)$, $z,
\xi\in\H)$, the integral
$$
\mathop{\int \kern-6pt\int }\limits_{F\times \H}
\left| \ln\varphi (z,\xi)\right| ^2 \left| 
\overline{\Delta(z)}\Delta(\xi)(
\overline{z}-\xi)^{12}\right| ^{\varepsilon} 
\left| k_1(z,\xi)\right| d^{2t}\,\d\nu_0(z)\eqno(7.1)$$
is (absolutely) convergent.

Since $k_1$ belongs to $\widehat{\ac }_{t-2-\varepsilon}$ we may free
up a small power $\alpha $ of $d$, so that the integral
$$
\mathop{\int \kern-6pt\int }\limits_{F\times \H}
\left| k_1(\overline{z},\xi)
\right| d^{2t-4-\alpha}
\,\d\nu_0(z,\xi)$$
is still convergent.
We  proved in \S 5 that
for any $\varepsilon<\varepsilon '$,
   there exists a positive  constant $C_{\varepsilon,
\varepsilon '}$, such that
$$\left| \ln\Delta(z)\Delta^\varepsilon(z)\right| \leq C_{\varepsilon,
\varepsilon '}\displaystyle {\Re z \over (\Im z)^2}
\e^{-\varepsilon'   \Im z}.$$
When evaluating the integral in (7.1), we will have to find an estimate for
  each of the terms that arise
by writing
$$\ln \varphi(\overline z, \xi)=\ln \overline{\Delta(z)}+
\ln \Delta({\xi}) +12\ln ({\overline{z}-
\xi}).$$
After taking the square,  we see that it remains to prove that the  integrals
  containing the following quadratic terms are
finite:

$$\left| \ln \Delta(z)\right| ^2 
\left| \Delta(z)\right| ^\varepsilon d^{\alpha},$$
  $$\left| \ln \Delta(\xi)\right| ^2 
\left| \Delta(\xi)\right| ^\varepsilon d^{\alpha},$$
  $$\left| \ln \zxi\right| ^2 
\left| \Delta^\varepsilon (\xi)\right| . $$
We analyze, for example, the term involving $\left| \ln \Delta(z)\right| ^2$.
By using Corollary 5.6, we note that the integral is consequently bounded by
$$\mathop{\int \kern-6pt\int }\limits_{F\times \H}
\displaystyle{
(\Re z)^2\over{
(\Im z)^6}}
\e^{-\varepsilon' \Im z}
\e^{-\varepsilon \Im \xi}
\left| \overline{z}-\xi \right| ^{12 \varepsilon}\cdot
(d(z,\xi))^{2t}
\left| k_1(\overline{z}, \xi)\right| ^2 \,\d\nu_0(z).$$
We write $(d(z,\xi))^{2t}=(d(z,\xi))^{2(t-3)}\cdot d(z,\xi)^{6}$
to get that the above integral is bounded by
$$\mathop{\int \kern-6pt\int }\limits_{F\times \H}
\e^{-\varepsilon '\Im z}
\e^{-\varepsilon \Im \xi}
\displaystyle{
(\Im \xi)^6\over{
|\overline{z}-\xi|^{12-12\varepsilon} }}\cdot
(d(z,\xi))^{2(t-3)}
\left| k_1(z,\xi)\right| ^2 \,\d\nu_0(z,\xi).$$
Because of the term $\e^{-\varepsilon \Im \xi}$, by
eventually multiplying with a constant, we can neglect the
term ${(\Im \xi)}^6$.

Thus we are led to analyze the following integral:
$$\mathop{\int \kern-6pt\int }\limits_{F\times \H}
(\Re z)^2
\e^{-\varepsilon '\Im z}
\e^{-\varepsilon \Im \xi}
\displaystyle {1\over {|\overline{z}-\xi}|^{12-12\varepsilon}}
(d(z,\xi))^{2(t-3)}
\left| k_1(z,\xi)\right| ^2 \,\d\nu_0(z).$$
Because $(z,\xi)\in F\times \H$, it follows that there is
a constant $C$ such that $\left| \overline{z}-\xi\right| \geq C$,
for $z,\xi\in F\times \H$.
Also $(\Re z)^2 /
|\overline{z}-\xi|^2$ is bounded from above on this region.

Thus the above integral is bounded by a constant times
$$\mathop{\int \kern-6pt\int }\limits_{F\times \H}
(d(z,\xi))^{2(t-3)}
\left| k_1(\overline{z},\xi)\right| ^2 \,\d\nu_0(z,\xi),$$
which is finite if $k_1\in \ac_{t-3-\varepsilon}$.

The terms with $\left| \ln \zxi\right| $ are solved
by absorbing $\left| \ln \zxi\right| $ into some power
of $\left| \overline{z}-\xi\right| $.

Clearly $\pc \mc_{\ln \overline{\varphi}}$ has the same
domain
as
$\pc \mc_{\ln {\varphi}}$. This is precisely the vector space of all
  $k\in L^2(\ac_t)$ such that
$\left| k(\ln\varphi)\right| ^2
=\left| \smash{k\overline{(\ln\varphi)}}\right| ^2$ is summable
on $F\times \H$, with respect to the measure $d^t \,\d\nu_0\times \d\nu_0$.
This completes the proof.\hfill\qed\vskip6pt

We introduce the following definition which will be used in
the dual solution for the cohomology problem, corresponding to ${\cal C}_t.$\vskip6pt

{\sm Definition} 7.4.\kern.3em
{\it  Let $\chi_{s,t}\colon  \ac_t\rightarrow \ac_s$
be defined by the formula
$$\chi_{s,t}{k}=
S_{\Delta^{(t-s)/12}}^*
k S_{\Delta^
  {(t-s)/12}},$$
for $k$ in $ac_t$. Here $s\leq t$.}\vskip6pt

In the next proposition we analyze the relation between the derivative
of $\chi_{s,t}$ at $s=t$, $s\nearrow t$, with the derivative
of $\theta_{s',t}$, at $s'=t$, $(s'\searrow t)$, introduced
in  Section 4.\vskip6pt

{\sm Definition} 7.5.\kern.3em
{\it For $t>1$, we let $\dc_t^+$ be the algebra 
consisting of all
$k$ in $\ac_s$ that, for some $s<t$, are of the form
$S^*_{\Delta^\varepsilon}k_1S_{\Delta^\varepsilon}$, for some
$\varepsilon>0$ such that $s+\varepsilon<t$ and $k_1\in \ac_{s+\varepsilon}$.

Clearly $\dc_t^+$ is a weakly dense, unital subalgebra of $\ac_t$.}\vskip6pt

{\sm Lemma} 7.6.\kern.3em
{\it Fix $t>1$. Assume that $k$ in $\dc_t^+$ has
the expression $k=S^*_{\Delta^
\varepsilon} k_1 S_{\Delta^\varepsilon}\mkern-1mu$,
$k_1\in \ac_s$, $\varepsilon>0$, $\varepsilon+s<t$.
Then
$$k=\tc_{\overline{\varphi^\varepsilon}}(k_1)\tc_{ \overline{
\overline{\Delta^{\varepsilon}(z)}
\Delta^{\varepsilon}(\xi)\zxi^{
12\varepsilon} }}
(k_1).$$
}\vskip6pt

{\it Remark.} Note that by putting the variables $\overline{z},\xi$, 
we indicated
that $k_1$ is multiplied by a function, that contrary to $k_1$, is
antianalytic in the second variable and analytic in the first.
Thus $\tc_{\overline{\varphi^\varepsilon}}$
corresponds to a Toeplitz operator with  an ``antianalytic'' symbol.\vskip6pt

{\it Proof of Lemma} 7.6. Let $T(k)=S^*_{\Delta^\varepsilon}
  k S_{\Delta^\varepsilon}$. The statement
  follows
immediately from the fact
that the adjoint of $T$, as a map on
$L^2(\ac_t)$, is $l\rightarrow\break
S_{\Delta^\varepsilon} l S^*_{\Delta^\varepsilon}.$\hfill\qed\vskip6pt

In the next proposition we clarify the relation between the operator 
$Y_tk$ defined
as $\left.{\d\;\over{\d s}}\chi_{s,t}(k)
\right|_{s=t;\;s\nearrow t} $ and the operator
$X_t$ introduced in Section 4.

First we recall that the ``real part'' associated with the
deformation is given by the
Dirichlet form ${\cal E}_s(k,l)={\d\;\over{\d s}}
  \tau_{\ac_s}(k\ast_sl)=
\tau_{\ac_s}
c_s(k,l).$\vskip6pt

{\sm Definition} 7.7.\kern.3em
{\it  Recall {\rm (}from Section {\rm 6)}
that the real part of the cocycle ${\cal C}_t$
  is the operator $R_t$ given by
by the formula
$$
\left\langle R_t k,l^*\right\rangle 
=-{1\over{2}} \left.\displaystyle {\d\;\over{\d s}}\cdot
\tau_{\ac_s} (k*_sl)\right|_{{{s=t}}\atop {{s\searrow t}}}=-{1\over 2 }
{\cal E}_s(k,l). $$
This holds for all $k,l$ in $\bigcup_{r<t}
L^2(\ac_r)$, where $L^2(\ac_r)$
is identified with a vector subspace of $L^2(\ac_t)$ via the symbol map
$\Psi_{t,r}$.

Moreover, in Section {\rm 6} we proved that  $R_t$ has the following expression:
$$R_t=\tc_{\ln d}- {1\over{2}} \displaystyle {c_t' \over{c_t}}\Id .$$
}\vskip6pt

In the next proposition we construct the dual object for the
  generator used in Section 4.\vskip6pt

{\sm Proposition} 7.8.\kern.3em
{\it For any $k$ in $\dc_t^+ \subseteq \ac_t$,
the limit  $$Y_t(k)= \left. \displaystyle {\d\;\over{\d s}}\Psi_{t,s}(
\chi_{s,t}{k})\right|_{{{s=t}}\atop {{s\nearrow t}}} $$
exists in $L^2(\ac_t)$.
Moreover, we have that
$$Y_t=-\left(\tc_{(1/12)\ln\varphi}\right)^*
-\displaystyle{c_t'\over{c_t}}\Id +2R_t.$$
  The adjoint  $(\tc_{(1/12)\ln \varphi})^* $ is obtained
by first restricting $\tc_{(1/12)\ln \varphi}$ to $\widetilde{\dc}_t^0$ and
  then taking the adjoint.
}\vskip6pt

{\it Proof.} Indeed, $\chi_{s,t}(k)$ may be identified
with the {\it Toeplitz operator {\rm (}on $L^2(\ac_s)$})
with symbol
  $$\overline{\varphi^{(t-s)/12}} = \overline 
{\left[ \overline{\Delta^{(t-s)/12}}(z)
\Delta^{(t-s)/12}(\xi)\zxi^{(t-s)}\right] }.$$

Thus $\chi_{s,t}(k)$ is (modulo a multiplicative constant)
$$\pc_s\left[ \mc_{\overline{\varphi^{(t-s)/12}}} k\right] ,$$
and hence $\Psi_{t,s}\chi_{s,t}(k)$ is
$$\Psi_{t,s}\pc_s\left[ \mc_{\overline{\varphi^{(t-s)/12}}} (k)\right] $$

  The derivative at $s=t$  consequently involves two
components.

  One component is the derivative
$\left. {\d\;\over{\d s}}
\Psi_{t,s}
(\pc _s k)\right|_{s\nearrow t}$, which gives the
summand corresponding to $R_t$, i.e., $-\left( c_t' / c_t\right) +2R_t$.

The other component is
${\d\;\over{\d s}}
\pc_t(\mc_{\overline{\varphi^{(t-s)/12}}} (k))$,
which gives the multiplication-by-$\varphi$ part. Indeed, recall that
$k\ \hbox{belongs\ to\ }\ \dc_t^+ \subseteq \ac_t$, and hence
$k$ is of the form
$S^*_{\Delta^\varepsilon}k_1S_{\Delta^\varepsilon_0}$ for some
$\varepsilon_0>0$ such that $s+\varepsilon_0<t$ and $k_1\in 
\ac_{s+\varepsilon_0}$.
But then
$$\pc_t(\mc_{\overline{\varphi^{(t-s)/12}}} (k))=
\pc_t(\mc_{\overline{\varphi^{(t-s)/12}}} 
(\pc_s(\mc_{\overline{\varphi^{\varepsilon_0}}}(k_1)))).$$
Since $\overline\varphi$ plays the role of an antianalytic symbol, 
it follows that this
is further equal to
$$\pc_t \left( \mc_{\overline{\varphi^{[(t-s)/12+\varepsilon_0]}}} 
(k_1)\right).$$
The derivative (in the $s$ variable) of
  $s\rightarrow \overline{\varphi^{[(t-s)/12+\varepsilon_0]}}$ at $s=t$
  exists, by the method in Lemma 7.2 in $L_t$, and it is equal
to $$-{1\over 12}\ln \overline\varphi\cdot\overline{
\varphi^{\varepsilon_0}}\cdot k_1.$$
Thus, in the Hilbert space $L^2(\ac_t)$, we have that
$$\eqalign{\displaystyle{\d\;\over{\d s}}
\pc_t(\mc_{\varphi^{(t-s)/12}}  (k))
&=-\pc_t( \mc_{(1/12)(\ln 
\overline{\varphi})\overline{\varphi^{\varepsilon_0}}} (k_1))\cr
\noalign{\medskip}
&=-\pc_t( \mc_{(1/12)
\ln 
\overline{\varphi}}(\pc_{t+12\varepsilon_0}\overline{\varphi^{\varepsilon_0}} 
(k_1)))\cr
\noalign{\medskip}
&=-\pc_t( \mc_{(1/12)\ln \overline{\varphi}}(k)).}$$
This completes the proof.\hfill\qed\vskip6pt

We use the above arguments to prove that also the operator  $Y_t
={\d\;\over{\d s}}\chi_{s,t}$
implements a coboundary for ${\cal C}_t$.\vskip6pt

{\sm Lemma} 7.9.\kern.3em {\it
For $k,l$ in $\dc_t^+$, we have that
$\left. {\d\;\over{\d s}}\chi_{s,t}(k *_t
\varphi^{t-s}(\overline{z},\xi) *_t l)
\right|_{s\nearrow t}$
is equal to $Y_t(k*_tl)-k*_t\Lambda(1)*_t l$.
}\vskip6pt

{\it Proof.} Since $k,l$ are in $\dc_t^+$, there exists $\varepsilon_0>0$ and
  there are $k_1, l_1\in \ac_{t+12 \varepsilon_0}$ such that
  $k=
S^*_{\Delta^{\varepsilon_0}} k_1
S_{\Delta^{\varepsilon_0}}$,
$l=
S^*_{\Delta^{\varepsilon_0}} l_1
S_{\Delta^{\varepsilon_0}}$. This gives that
$$k*_t\Lambda(1) *_tl=
S^*_{\Delta^{\varepsilon}}
[k_1*_{t+\varepsilon}
\varphi^{\varepsilon}
\ln\varphi *_{t+\varepsilon}l_1]
S_{\Delta^{\varepsilon}},$$
where by $\varphi^{\varepsilon}\ln \varphi$ we
understand the unbounded operator defined
in \S 3, corresponding to
$$\varphi(\overline{z},\xi)^\varepsilon\ln \varphi(\overline{z},\xi).$$

As in the proof of Proposition 4.5, when computing this derivative, we
have a trivial summand plus a more complicated summand, corresponding
to the symbol
$$
\lim\limits_{s\nearrow t}\Psi_{t,s}\left[\pc_s
\overline{\varphi^{t-s}}
\left [
k*_t \displaystyle{
\varphi^{t-s}-\Id \over{
t-s}}*_t l \right]\right]
.$$

Because of the assumptions, the inside term
$$k*_t \displaystyle{
\varphi^{t-s}-\Id \over{
t-s}} *_t l$$
is equal to
$$\eqalign{
&S_{\Delta^{\varepsilon_0}}
\left[k_1*_{t+\varepsilon_0}
\varphi^{\varepsilon_0}
\displaystyle{
\varphi^{t-s}-\Id \over{
t-s}}*_{t+\varepsilon_0}l_1\right]S_{\Delta^{\varepsilon_0}}^*\cr
\noalign{\medskip}
&\qquad =\pc_t\left(
\overline{\varphi^{\varepsilon_0}}
\left(
k_1*_{t+\varepsilon_0}
\varphi^{\varepsilon_0}\displaystyle{
\varphi^{t-s}-\Id \over{
t-s}}*_{t+\varepsilon_0}l_1
\right)
\right) .}\eqno(7.2)$$

But the  methods in the proof of  the density of  the domain of
$\mc_{\ln\varphi}$ may also be used to prove that
$$\varphi^{\varepsilon_0}\left(
\displaystyle{
\varphi^{t-s}-\Id \over{
t-s}}\right)$$ converges, as $s\nearrow t$, in
$L^2(\ac_{t+12\varepsilon_0})$ to
$-\varphi^{\varepsilon_0}\ln \varphi$.

Since $\Psi_{t,s} \pc_s$ converges strongly to the identity, and
the norm of $\Psi_{t,s} \pc_s$ as an operator from
$L^2(\widehat{\ac}_s)$ into $L^2(\ac_t)$ is bounded
by $1$, it follows that the expression in (7.2) converges to
$$ \pc_t(\overline{\varphi}_{\varepsilon_0}[
k_1*_{t+\varepsilon_0}\varphi^{\varepsilon_0}(-\ln \varphi)
*_{t+\varepsilon_0}l_1
]),$$
which is  $k *_t(-\Lambda(1))*_t l$.\hfill\qed\vskip6pt

Similarly we have the following lemma.\vskip6pt

{\sm Lemma} 7.10.\kern.3em {\it
For $k,l$ in $\dc_t^+$ we have that
$$\left.\displaystyle {\d\;\over{\d s}}
\left[ \chi_{t,s}(k) *_s \chi_{t,s}(l)\right] 
\right|_{s\nearrow t}=
Y_t(k)*_tl+{\cal C}_t(k,l)+k*_tY_t(l).
$$}\vskip6pt

{\it Proof.} Again this derivative has three
summands. The first summand is
$$\lim_{s\nearrow t}
\displaystyle{
\chi_{t,s}(k) *_s \chi_{t,s}(l)-
\chi_{t,s}(k) *_t \chi_{t,s}(l)
\over{
t-s}}.$$
The same type of argument as in Proposition 4.5 gives that
this is ${\cal C}_t(k,l)$.

{}From the remaining two summands, the only one that is complicated is
$$\chi_{s,t}(k) *_t \displaystyle{\chi_{s,t}(l)-l\over{t-s}}.$$

Because for $l$ in $\dc_+^t$ we have that
$
\left( \chi_{s,t}(l)-l \right) / \left( t-s\right) $
converges in $L^2(\ac_t)$ to $Y_t l$, and
since $\chi_{s,t}(k)$ is bounded in $L^2(\ac_t)$
as $s\nearrow t$, it follows that this term converges
too, to $k*_tY_tl$.

The remaining term trivially converges to
$ Y_t k*_t l$.

This completes the proof of the lemma.\hfill\qed\vskip6pt

As
  a corollary  we obtain the following result.\vskip6pt

{\sm Proposition} 7.11.\kern.3em
{\it
Let $\dc_t^+$ be as in Definition {\rm 7.5.}  Assume $t>3$; then
for all $k,l$ in $\dc_t^+$, we have that
$$Y_t(k*_t l)-Y_t k*_t l- k*_t Y_t l-
k*_t {1\over{12}} \ln \varphi *_tl={\cal C}_t(k,l).$$

Here by ${1\over{12}}\varphi(\overline{z},\xi)$ we
understand $\Lambda(1)=\mc_{\ln\varphi}(1)$, the
operator constructed in Corollary~{\rm 3.6} and in Lemma~{\rm 4.4.}
}\vskip6pt

{\it Proof.} Indeed the identity
$$\chi_{s,t}(k *_t \varphi^{t-s}*_t l)=
\chi_{s,t}(k) *_s \chi_{s,t}(l)$$ is obvious, valid for all
$k,l$ in $\dc_t^+$, $s\leq t$.

By differentiation, and using the two previous lemmas, we get our result.\hfill\qed\vskip6pt

{\it Remark.} Recall that in Proposition 6.4 we proved that if
$\Lambda_0(k)=\mc_{\ln \varphi}k$ (for $k$ in
$\dc_0^t$), then, denoting $S=\Sym\nolimits_{ \varphi}$, we have
$$\left\langle
\Lambda_0(k),l\right\rangle
+\left\langle k,\Lambda_0(l)\right\rangle
=2\left\langle
S
k,l\right\rangle
+2\left\langle
R_tk,l\right\rangle
+\displaystyle{
c_t'\over{c_t}}\left\langle
k,l\right\rangle .$$

If $k$ belonged to $\widetilde{\dc}_t^0$ (which is the domain
of $\Lambda_0$) and also to the domain of $Y_t$,
which is $\dc_t^+$, then the above relation could be
rewritten as
$$\Lambda_0+\Lambda_0^*=2S+2R_t+
\displaystyle{
c_t'\over{c_t}}\eqno(7.3).$$
Recall that  $S=\Sym\nolimits_{ \varphi}$ is the operator defined by
$$\left\langle
Sk,l\right\rangle
=\tau_{\ac_t}(kTl^*)+\tau_{\ac_t}(l^*Tk).$$

But on the intersection of the domains  we have
(from Proposition 7.8)
that
$$Y_t=-\Lambda_0^*+2R_t-
\displaystyle{
c_t'\over{c_t}}\Id .\eqno(7.4)$$
Consequently
$$
\Lambda_0^*=-Y_t+2R_t-
\displaystyle{
c_t'\over{c_t}}\Id .
$$
Thus, by (7.3), for  $k$ in $\widetilde{\dc}_t^0\cap\dc_t^+$  we get that
$$\Lambda_0=2S
+Y_t+
2\displaystyle{
c_t'\over{c_t}}\Id ,$$
and hence that $$X_t=\Lambda_0- \displaystyle{
c_t'\over{c_t}}\Id =2S+Y_t+\displaystyle{
c_t'\over{c_t}}\Id  ,\eqno(7.5)$$
where equality holds on $\widetilde{\dc}_t^0\cap\dc_t^+$.

Now we compare the way $X_t$, $Y_t$ implement
a coboundary for ${\cal C}_t(k,l)$.
Recall the notation
$(\nabla \Phi)(k,l)=\Phi(k,l)-k\Phi(l)-\Phi(k)l$.

Thus we have proved that
$$\nabla X_t(k,l)={\cal C}_t(k,l)-kT^t_{\ln\varphi}\varphi,\qquad 
k,l\ \hbox{in}\
\widetilde{\dc}_t^0 ,\eqno(7.6)$$
$$\nabla Y_t(k,l)={\cal C}_t(k,l)-k\Lambda(1)l,\qquad 
k,l\ \hbox{in}\ \dc_t^+.
\eqno(7.7)$$

Now if $k,l$ were in $\widetilde{\dc}_t^0\cap\dc_t^+$, it
would follow, by substituting (7.5) into (7.6),
that
$$2\nabla S_t(k,l)+\nabla Y_t(k,l)-
\displaystyle{
c_t'\over{c_t}} k *_t l={\cal C}_t(k,l)-kT_{\ln\varphi}^t l. \eqno(7.8)$$

By using (7.7) in (7.8) we get
$$2(\nabla S_t)(k,l)- k\Lambda(1)l-
\displaystyle{
c_t'\over{c_t}} (k*_tl)=-kT_{\ln \varphi}^t l,
$$
and thus for $k,l$ in $\widetilde{\dc}_t^0\cap\dc_t^+$ we would
get that
$$k\left[T^t_{\ln \varphi}-\left(\Lambda(1)+\displaystyle{
c_t'\over{c_t}}\right)\right]l=
2\nabla S_t(k,l)\eqno(7.9)$$ for all $k,l$ in $\widetilde{\dc}_t^0\cap\dc_t^+$.
But recall that
$$\left\langle
S_t(k),l^*\right\rangle
=\tau(kTl^*+l^*Tk).$$
This corresponds, at least formally, to the fact that
$S_tk=kT+Tk$
and hence $(\nabla S_t)(k,l)$ is $2k T_{\ln \varphi}^t l$.

Thus (7.9) would imply directly that
$$T_{\ln \varphi}^t =\Lambda(1)+\displaystyle{
c_t'\over{c_t}}$$
  if $\widetilde{\dc}_t^0\cap\dc_t^+$ is nonzero.
\vfill \eject

\centerline {\bf Appendix}\vskip12pt
\centerline {\ninebf A more general coboundary for ${\cal C}_t$}\vskip12pt

In this appendix, we want to construct a more general solution for a
coboundary (which is necessarily unbounded, see [Ra]) for
${\cal C}_t$.
This will be constructed out of a measurable function $g$ that
has the same $\Gamma$-invariance properties as $\ln \Delta(z)$.
By this construction we will lose the complete positivity properties
of the solution.

Recall that $L_t$ consists of all kernels $k$ or $\H\times \H$,
that are diagonally $\Gamma$-invariant and square
summable on $F\times \H$, with respect to  the
measure $d^t\,\d\nu_0 \times \d \nu_0$. Also recall
that the elements in $L_t$ are canonically identified with
operators in the $\twoinfty$ factor of all
operators that commute with $\widetilde {\pi_t}(\Gamma)$, acting
on $L^2(\H,\d \nu_t)$.\vskip6pt

{\sm Proposition} A1.\kern.3em {\it Let $g$ be a measurable
function $\H$ such that the bivariable
function $\theta$ in $\H\times \H$ defined by
$$\theta(z,\xi)=\overline{g(z)}+g(\xi)+\ln\zxi$$
is $\Gamma$-invariant. {\rm (}It is this point which makes the problem 
solvable, by this method,
only for $\PSLtwoZ$.{\rm )}

Let $\dc_g,\widetilde{\dc}_g$ consist of all $k$ in $L^2(\ac_t)$
{\rm (}respectively $L_t${\rm )} such that $k\cdot \theta$ still
belongs to $L_t$. Let $\mc_\theta$ be the
{\rm (}unbounded\/{\rm )} operator with domain $\widetilde{\dc}_g$,  of
multiplication by $\theta$.
Let $\tc_{\theta}=\pc_t\mc_{\theta}|_{L^2(\ac_t)}$ and let
$T^t_{\theta}$ be the Toeplitz operator with symbol
$\theta(z,z)=\Re g(z)+\ln (\overline{z}-z)$.

Let $k,l$ be in ${\cal D}_g$ such that $k,l$ also belong
  to the domain of ${\cal C}_t(k,l)$.
Then
$$\mc_\theta (k*_tl)-
\mc_\theta k*_tl-k*_t \mc_\theta l+M_\theta(k,l)={\cal C}_t(k,l),
$$
where $M_\theta(k,l)$ is a bimodule map, equal to $k*_t T_{\theta}^t*_t l$,
if $T_{\theta}^t$ exists.

Consequently, by taking $ P_t$ on the left- and right-hand sides,
  the same will hold true for
$\tc_\theta={\cal P}_t \mc_\theta|_{L^2(\ac_t)}$.
}\vskip6pt

{\it Proof.}
Indeed ${\cal C}_t(k,l)$ is given by the kernel
$$\eqalign{{\cal C}_t(k,l)(\overline {z},\xi)&=\displaystyle
{ c_t'\over{ c_t}}(k*_tl)(\overline {z},\xi)\cr
\noalign{\medskip}
&\qquad +\displaystyle c_t\int_{\H}k(\overline {z},\eta)l(\overline {z},\xi)
[\overline{z},\eta, \overline{\eta},\xi]^t\ln
[\overline{z},\eta, \overline{\eta},\xi]\,\d \nu_0(\eta).}
$$
On the other hand
$$\mc_\theta (k*_tl)-
\mc_\theta k*_tl-k*_t \mc_\theta l$$
has the kernel
$$\displaystyle c_t\int_{\H}k(\overline {z},\eta)l(\overline {z},\xi)
[\overline{z},\eta, \overline{\eta},\xi]^t(\theta(z,\xi)-
\theta(z,\eta)-\theta (\eta,\xi))
\,\d \nu_0(\eta).$$

Since $\theta(z,\xi)-\theta(z,\eta)-\theta (\eta,\xi)$ is
equal to $\theta (\eta,\eta)$, it follows that
$${\cal C}_t(k,l)-[\mc_{\theta}(k*_tl)-
\mc_{\theta}k*_t l-  k*_t\mc_{\theta}l]$$
is given by the kernel
$$\zxi^t\int_{\H}
\displaystyle{
k(\overline{z},\eta) l(\overline{\eta},\xi)
\over{
\zueta^t\etxi^t}}
\left (\theta(\eta,\eta)
+\displaystyle{
c_t'\over{c_t}}\right)\,\d \nu_t(\eta),$$
which indeed corresponds to
$T_\theta^t+
\left( 
c_t' / c_t\right) \cdot\Id \ $, as
long  as we can make sense of the unbounded  Toeplitz  operator $T_\theta^t$.

A dual version could be obtained if we consider
$$\left\langle
2R_t k, l^*\right\rangle _{L^2(\ac_t)}=-
\left.\displaystyle{\d\;\over{\d s}}\tau_{\ac_s} (k *_s l^*)\right|_{s=t},$$
which is in other terms
$$2R_t=-\displaystyle{
c_t'\over{c_t}}-2 \tc_{\ln d}.$$

One can check immediately that
$$\eqalign{
&[{\cal C}_t(k,l)-(\nabla 2R_t)(k,l)](\overline{z},\xi)\cr
\noalign{\medskip}
&\qquad =\displaystyle{
c_t'\over{c_t}}\tau(k*_t l)+ c_t
\int _{\H}k(\overline{z},\eta)l(\overline{\eta},\xi)
\ln[\overline{z},\eta,\overline{\eta},\xi][\overline{z},\eta
,\overline{\eta},\xi]^t\,\d \nu_0(\eta)\cr
\noalign{\medskip}
&\qquad \qquad +
\int_{\H}k(\overline{z},\eta)l(\overline{\eta},\xi)
[-\ln( d(z,\xi))^2+\ln( d(z,\eta))^2+\ln(\d(\eta, \xi))^2]\,\d \nu_0(\eta)\cr
\noalign{\medskip}
&\qquad \qquad -\left(
-\displaystyle{
c_t'\over{c_t}}+
\displaystyle{
c_t'\over{c_t}}+
\displaystyle{
c_t'\over{c_t}}\right)\tau(k*_t l)\cr
\noalign{\medskip}
&\qquad =
c_t\int_{\H} k(\overline{z},\eta)l (\overline{\eta},\xi)
[\overline{z},\eta
,\overline{\eta},\xi]^t\cr
\noalign{\medskip}
&\qquad \qquad \cdot \{
\ln(
[\overline{z},\eta
,\overline{\eta},\xi])-
2\ln  d(z,\xi)+
2\ln  d(z,\eta)+
2\ln  d(\eta,\xi)\}\,\d \nu_0(\eta)\cr
\noalign{\medskip}
&\qquad =-c_t
\int_{\H} k(\overline{z},\eta)l (\overline{\eta},\xi)
[\overline{z},\eta
,\overline{\eta},\xi]\cdot
\ln
[ \overline{
\overline{z},\eta
,\overline{\eta},\xi}]\,\d \nu_0(\eta).}
$$
Then consider $g$ such that
$$\theta(z,\xi)=g(z)+\overline{g(\xi)}+\ln (z-\overline{\xi})$$
is $\Gamma$-invariant.

The same argument as above gives that for $k,l$ in $\dc(\mc_{\theta})$ s.t.\
$k*_tl\in \dc(\mc_{\theta})$ and $k,l$ in $\Dom (R_t)$, $k*_tl$ in
$\Dom (R_t)$
we have that
$$\mc_\theta (k*_tl)-
\mc_\theta k*_tl-k*_t \mc_\theta l+
k*_t T_{\theta}^t*_t l$$
is equal to ${\cal C}_t(k,l)-2\nabla R_t(k,l)$.

Finally remark that we have proved that for $k,l$ in $\dc_t^+$,
which is the vector space of all $k$ that are of
the form $\pc_t(k_1\overline{\varphi^\varepsilon})$,
the expression
$$\mc_{\overline{\varphi}} (k*_tl)-
\mc_{\overline{\varphi}} k*_tl-
k*_t\mc_{\overline{\varphi}} l+
k*_t\Lambda(1)*_tl-{\cal C}_t(k,l)$$
is orthogonal to $\pc_t({L^t})$ (in other words, if we apply $P_t$ to
the left and right we get $0$).

If we could extend the above relation to all
$k,l$ in $\Dom(\mc_{\overline \varphi})\cap \Dom R_t$ such
  that $k*_tl$ belongs to the same domain, then the above relation,
by the considerations at the end of \S 7, would imply that
$\Lambda(1)-\left( 
c_t' / c_t\right) \cdot\Id$ coincides (on an affiliated domain) with
$T^t_{\ln \varphi}$.

  Note that this corresponds formally to the fact that
$$T^t_{\ln \varphi}=T^t_{\ln (\overline{\eta}-\eta)}+
T^t_{\overline{\ln \Delta}}+
T^t_{\ln \Delta}.$$
On the other hand
$T^t_{\ln (\overline{\eta}-\eta)}=
{\cal P}_{\ln (\overline z-\xi)}-\left( 
c_t' / c_t\right) \cdot\Id$,
while $T^t_{\overline{\ln \Delta}}$ is clearly, on its domain,
$S^t_{\overline{\ln \Delta}}$ (and similarly for $T^t_{\ln \Delta}$).

  If the domains had nonzero intersection,
  one could directly conclude that
$$S^t_{\overline{\ln \Delta}}+S^t_{{\ln \Delta}}+
{\cal P}_{\ln (\overline z-\xi)}=\Lambda(1).$$

\vfill\eject



\ninepoint

\centerline{REFERENCES}\vskip12pt

\leftskip24pt

\parindent=-24pt

[Acc]\kern.5em L. Accardi, R.L. Hudson,
{\it Quantum stochastic flows and nonabelian cohomology},
Quantum Probability and Applications, V 
(Heidelberg, 1988)
(L. Accardi and W. von 
Waldenfelds, eds.),
Lecture Notes in Math., vol. 1442, Springer, Berlin, 1990,
pp. 54--69.

[A]\kern.5em Atiyah, M. F.
{\it Elliptic operators, discrete groups and von Neumann algebras},
  Colloque ``Analyse et Topologie'' en l'Honneur de Henri Cartan (Orsay, 1974),
 Asterisque, No. 32-33,
Soc. Math. France, Paris, 1976, pp. 43--72.

[AS]\kern.5em            M.F. Atyiah, W. Schmidt,
{\it  A geometric construction
of the discrete series for semi\-sim\-ple Lie groups}, Invent. Math. {\bf 42}
(1977), 1--62.

[Be]\kern.5em F.A. Berezin, {\it General Concept of Quantisation},
Comm. Math. Phys, {\bf 40}(1975), 153--174.

[Ba]\kern.5em V. Bargmann, {\it Irreducible unitary representations of the
Lorenz group}, Annals of Mathematics,  {\bf 48} (1947), 568-640.

[BG]\kern.5em Barge, J.; Ghys, E.
{\it Surfaces et cohomologie born\'ee},
Invent. Math. 92 (1988), no. 3, 509--526.

[CF] \kern.5em Chebotarev, A. M. ; Fagnola, F.
{\it Sufficient conditions for conservativity of minimal
  quantum dynamical semigroups},
J. Funct. Anal. {\bf 153} (1998), 382--404.

[CE]\kern.5em Christensen, E; Evans, D. E.
{\it Cohomology of operator algebras and quantum dynamical semigroups},
J. London Math. Soc, {\bf 20} (1979), 358--368.

[CH] \kern.5em Cohen, P. Beazley,; Hudson, R. L.,
{\it Generators of quantum stochastic flows and cyclic cohomology},
Math. Proc. Cambridge Philos. Soc. {\bf 123} (1998),  345--363.

[Co] \kern.5em A. Connes, {\it Sur la Theorie Non Commutative de
l'Integration}, Springer Verlag, 725.

[Co1]\kern.5em A. Connes, {\it Noncommutative Geometry},
  Academic Press, Harcourt Brace \& Company Publishers, 1995.

[Co2]
\kern.5em A. Connes,
{\it On the spatial theory of von Neumann algebras},
J. Funct. Anal. {\bf 35} (1980), 153--164. 

[CoCu] \kern.5em Connes, A.; Cuntz, J., {\it
Quasi homomorphismes, cohomologie cyclique et positivit\'e
  [Quasihomomorphisms, cyclic homology and
positivity]},
Comm. Math. Phys. {\bf 114} (1988),  515--526.

[Dav] \kern.5em Davies, E. B.
{\it Uniqueness of the standard form of the generator
  of a quantum dynamical semigroup},
Rep. Math. Phys. {\bf 17} (1980),  249--264.

[Dy] \kern.5em Dykema, K.,
{\it Interpolated free group factors},
Pacific J. Math. {\bf 163} (1994), no. 1, 123--135.

[GKS] \kern.5em Gorini, V.; Kossakowski, A.; Sudarshan, E. C. G.,
{\it Completely positive dynamical semigroups of $N$-level systems},
J. Mathematical Phys. {\bf 17} (1976),  821--825.

[GS] \kern.5em Goswami, D.; Sinha, K. B.,
{\it Hilbert modules and stochastic dilation of a quantum
  dynamical semigroup on a von Neumann algebra}   Comm. Math. Phys. 
{\bf 205} (1999),
   377--403.

[GHJ] \kern.5em      F. Goodman, P. de la Harpe, V.F.R. Jones, {\it Coxeter
Graphs and Towers of Algebras}, Springer Verlag, New York, Berlin,
Heidelberg, 1989.

[Jo] \kern.5em Jolissaint, P.; Valette, A.,
{\it Normes de Sobolev et convoluteurs born\' es  sur $L\sp 2(G)$},
Ann. Inst. Fourier (Grenoble) {\bf 41} (1991),  797--822.

[Ho] \kern.5em Holevo, A. S.
{\it Covariant quantum dynamical semigroups: unbounded generators}, in
  Irreversibility and causality (Goslar,
1996), 67--81,
Lecture Notes in Phys., {\bf 504},
Springer, Berlin, 1998.

  [HP] \kern.5em Hudson, R. L.; Parthasarathy, K. R.,
{\it Quantum Ito's formula and stochastic evolutions},
Comm. Math. Phys. {\bf 93} (1984),  301--323.

  [Li] \kern.5em Lindblad, G.
{\it On the generators of quantum dynamical semigroups},
Comm. Math. Phys. {\bf 48} (1976),  119--130.

  [MS] \kern.5em Mohari, A.; Sinha, K. B.,
{\it Stochastic dilation of minimal quantum dynamical semigroup},
Proc. Indian Acad. Sci. Math. Sci. {\bf 102} (1992), 159--173.

[MvN]\kern.5em F. J. Murray; J. von Neumann, {\it On ring of Operators, IV},
  Annals of Mathematics, {\bf 44} (1943), 716-808.

[Mi] \kern.5emMiyake, T.,
{\it Modular Forms},
Springer-Verlag, Berlin-New York, 1989.

[Puk] \kern.5emPuk\'anszky, L.,
{\it The Plancherel formula for the universal covering group of}
 ${\rm SL}(R,2)$, Math. Ann.
{\bf 156} (1964), 96--143.

[Ra] \kern.5em R\u adulescu, F.,
{\it The $\Gamma$-equivariant form of the Berezin quantization
  of the upper half plane},
Mem. Amer. Math. Soc. {\bf 133} (1998), no. 630,

[Ra1] \kern.5em R\u adulescu, F.,
{\it On the von
Neumann Algebra of Toeplitz Operators with
Automorphic Symbol}, in Subfactors,
Proceedings of the Taniguchi Symposium on Operator
Algebras, eds.\ H. Araki, Y. Kawahigashi, H. Kosaki,  World Scientific,
  Singapore-New Jersey,1994.

[Ra2] \kern.5em R\u adulescu, F.,
{\it Quantum dynamics and Berezin's deformation quantization},
  in Operator algebras and quantum field
theory (Rome, 1996), 383--389,
Internat. Press, Cambridge, MA, 1997.

[Sal] \kern.5em  P. Sally, {\it Analytic
  Continuation of the Irreducible
Unitary Representations of the Universal
  Covering Group}, Memoirs A. M. S. , 1968.

[Sau]  \kern.5em Sauvageot, J.-L.,
{\it Quantum Dirichlet forms, differential calculus and semigroups}, in
Quantum probability and applications, V (Heidelberg, 1988), 334--346,
Lecture Notes in Math., {\bf 1442},
Springer, Berlin, 1990.

[ShS]  \kern.5em Shapiro, H. S.; Shields, A. L.
{\it On the zeros of functions with finite Dirichlet integral and 
some related function spaces},
Math. Z., {\bf 80} (1962), 217--229.

[Vo]  \kern.5em Voiculescu, D.,
{\it Limit laws for random matrices and free products},
Invent. Math. {\bf 104} (1991),  201--220.


\bye